\documentclass[11pt]{article}

\usepackage{apacite}
\usepackage{color}
\usepackage[dvipsnames]{xcolor}
\usepackage{natbib}
\setcitestyle{authoryear,round}

\usepackage{epsfig,lscape}
\usepackage{amssymb,amsfonts,amsmath,amsthm}
\usepackage{rotating}
\usepackage{bm}
\usepackage{bbm}
\usepackage[hypertexnames=false]{hyperref}
\usepackage{mathtools}

\usepackage{threeparttable}
\usepackage{booktabs}
\usepackage{subcaption}
\usepackage[left=1in,right=1in,top=.9in,bottom=1.1in]{geometry}

\def\me{\mathrm e}

\def\dif{\mathrm d}

\def\N{\mathrm{N}}

\def\T{ {\mathrm{\scriptscriptstyle T}} }

\def\a{\alpha}
\def\b{\beta}
\def\g{\gamma}
\def\w{\omega}

\def\bbR{\mathbb R}

\def\argmin{\mathrm{argmin}}

\def\one{\mathbf{I}}
\def\two{\mathbf{II}}

\newenvironment{prf}
{\noindent \textbf{Proof.}}{\hfill $\Box$ \vspace{.1in}}

\newtheorem{thm}{Theorem}
\newtheorem{lem}{Lemma}
\newtheorem{pro}{Proposition}
\newtheorem{cor}{Corollary}
\newtheorem{ass}{Assumption}

\theoremstyle{definition}

\theoremstyle{definition}

\allowdisplaybreaks
\begin{document}

\renewcommand\thefootnote{\fnsymbol{footnote}}

\begin{titlepage}
    \begin{center}
        {\Large Understanding Accelerated Gradient Methods: Lyapunov Analyses and Hamiltonian Assisted Interpretations}

        \vspace{.15in} Penghui Fu\footnotemark[1] and Zhiqiang Tan\footnotemark[1]

        \vspace{.1in}
        \today
    \end{center}

    \footnotetext[1]{Department of Statistics, Rutgers University. Address: 110 Frelinghuysen Road,
Piscataway, NJ 08854. E-mails: penghui.fu@rutgers.edu, ztan@stat.rutgers.edu.}

    \paragraph{Abstract.} We formulate two classes of first-order algorithms more general than previously studied
    for minimizing smooth and strongly convex or, respectively, smooth and convex functions.
    We establish sufficient conditions, via new discrete Lyapunov analyses,
    for achieving accelerated convergence rates which match Nesterov's methods in the strongly and general convex settings.
    Next, we study the convergence of limiting ordinary differential equations (ODEs) and point out currently notable gaps between the convergence properties of
    the corresponding algorithms and ODEs. Finally, we propose a novel class of discrete algorithms, called the Hamiltonian assisted gradient method,
    directly based on a Hamiltonian function and several interpretable operations,
    and then demonstrate meaningful and unified interpretations of our acceleration conditions.
    %\vspace{-.2in}

\end{titlepage}

\renewcommand\thefootnote{\arabic{footnote}}

\section{Introduction}

Optimization plays a vital role in machine learning, statistics, and many other fields.
%In particular, first-order\footnotemark[2] methods are desirable
%in handling large-scale problems due to low per-iteration costs and easy adaption to parallel and distributed computing.
In the optimization literature, there exists a striking phenomenon
where after suitable modifications of a first-order method,\footnote{First-order methods refer to methods using function values and gradients only, whereas second-order methods additionally rely on the Hessian matrices or their approximations.}
the convergence guarantee can be improved,
often attaining the complexity lower bound, with a similar computational cost as before.
Such an acceleration has been widely studied
since the seminal work of \cite{nesterov1983method},
which improves gradient descent for minimizing smooth convex functions.
Examples include constrained optimization \citep{nesterov2018lectures},
mirror descent with a non-Euclidean norm \citep{krichene2015accelerated}, composite optimization with a proximable function \citep{beck2009fast},
primal-dual splitting \citep{chambolle2011first}, stochastic
gradient methods \citep{zhang2017stochastic, allen2018katyusha}, and others.

Despite extensive research, the scope and mechanism of acceleration remains to be fully understood,
even including the original acceleration of gradient descent for smooth convex optimization.
\cite{nesterov1988approach} established the improved convergence through
a particular technique, called estimation sequences, which is instrumental in the early study of accelerated methods
but does not offer a heuristic explanation for when and how acceleration can be achieved.

Recently, Nesterov's acceleration has been studied from various perspectives
while circumventing the estimation-sequence technique, conceptually or analytically.
For example, geometric formulations are proposed by coupling gradient descent and mirror descent \citep{allen2017linear},
and by averaging two minimizers of an upper and a lower quadratic bound for the objective function \citep{bubeck2015geometric,drusvyatskiy2018optimal}.
Another useful approach is to relate discrete algorithms to their continuous limits, which are ordinary differential equations (ODEs)
obtained by letting the stepsize in the discrete algorithms tend to zero.
Continuous ODEs are usually more tractable than their discrete counterparts, and can be studied by exploiting
a rich set of analytical tools from dynamical systems and control theory.
The analyses and properties of continuous ODEs can in turn provide insights about the behavior of discrete algorithms \citep{qian1999momentum,su2016differential,yang2018physical,sun2020high,shi2021understanding}.
Conversely, ODEs can be directly formulated and then their numerical discretizations are studied
\citep{wibisono2016variational,wilson2021lyapunov}.
It has been argued that acceleration can be attributed to suitable discretizations from certain ODEs,
such as the symplectic integrator \citep{shi2019symplectic,francca2020conformal,muehlebach2021optimization} and Runge--Kutta integrator \citep{zhang2018runge}.

A prominent strategy emerging from the recent literature, motivated by analysis of ODEs based on Lyapunov functions, is
the use of Lyapunov analysis to establish convergence properties
including accelerated rates for discrete algorithms.
See \cite{bansal2019potential} and \cite{aspremont2021acceleration} for overviews of Lyapunov-based proofs for gradient methods.
A central step in Lyapunov analysis is to construct an appropriate Lyapunov function (aka energy functional),
which satisfies a contraction inequality to ensure the desired convergence.
But this task is nontrivial especially for discrete algorithms.
As mentioned above, considerable progress has been made 
by building and exploiting connections between the corresponding ODEs and discrete algorithms
to facilitate the construction of Lyapunov functions for analyzing discrete algorithms in a systematic manner
\citep{wibisono2016variational,wilson2021lyapunov, shi2021understanding}.

The ODE-based approach, however, currently does not fully account for different behaviors among discrete algorithms,
in particular, whether Nesterov's acceleration is achieved.
For example,  for minimizing strongly convex functions,  both Nesterov's accelerated gradient method (NAG-SC) and Polyak's heavy-ball method
admit the same limiting (low-resolution) ODE, but only Nesterov's method is known
to achieve accelerated convergence.
A more elaborate approach  has been proposed by \cite{shi2021understanding} using high-resolution ODEs which
are defined by retaining certain terms which would otherwise vanish as the stepsize tends to 0.
In this approach, a continuous Lyapunov function is first constructed to analyze a high-resolution ODE
and then translated into a discrete Lyapunov function to analyze the original discrete algorithm.\footnote{Currently, the translation of continuous Lyapunov functions into discrete ones seems to still require ad hoc considerations
when using high-resolution ODEs as well as low-resolution ODEs.}
While this approach succeeds in providing a Lyapunov analysis to demonstrate the acceleration of NAG-SC,
there remain some gaps.
In fact, the high-resolution ODEs for NAG-SC and heavy-ball converge at the same rates,
as shown in Proposition \ref{pro:high-res-ODE NAG and HB}.
In a similar manner, as shown in Proposition \ref{pro:high-res-ODE-C}, the same convergence rates are achieved by
the high-resolution ODEs for a class of algorithms generalizing Nesterov's accelerated gradient method (NAG-C),
but only a subclass of algorithms are known to achieve accelerated convergence as NAG-C
for minimizing general convex functions.

We take a direct approach to studying the scope and mechanism of Nesterov's acceleration
for discrete algorithms. Our main contributions can be summarized as follows.\vspace{-.05in}
\begin{itemize} \addtolength{\itemsep}{-.05in}

\item We formulate two classes of algorithms which are more general than previously studied
and establish sufficient conditions for when the algorithms achieve Nesterov's acceleration as well as non-accelerated convergence,
in both the strongly and general convex settings.
Our proofs are developed by constructing new discrete Lyapunov functions applicable to the general classes. See Section \ref{sec:outline}
for a comparison of our and existing Lyapunov analyses.

\item We study the low-resolution and high-resolution ODEs
derived from our classes of algorithms and compare the conditions for when acceleration is achieved by
the discrete algorithms and ODEs.
Our comparison points to currently notable gaps between the convergence properties of the corresponding algorithms and ODEs.
See Section \ref{sec:ODE-interpretation} for a more detailed summary.

\item We propose a novel class of discrete algorithms, Hamiltonian assisted gradient method (HAG), directly
based on a Hamiltonian function, and demonstrate that the conditions from our convergence results can be
interpreted through HAG, with regard to the momentum and
gradient correction terms, in a meaningful and unified manner in both the strongly and
general convex settings. See Section \ref{subsec:interpretation-HAG} for a more detailed summary.
\end{itemize}  \vspace{-.05in}
In the Supplement, we also present numerical results to illustrate different performances of the algorithms
falling inside versus outside the sufficient conditions for accelerated convergence in our theoretical results.
All technical details are deferred to the Supplement.

\textbf{Notation.}
We largely adopt the notation in \cite{nesterov2018lectures}. For a smooth and convex function $f$, consider the unconstrained minimization problem
\begin{equation}\label{eq:main_optimization}
    \min_{x\in\bbR^n} f(x).
\end{equation}
Denote by $x^*$ one of the minimizers of $f$ and $f^*=f(x^*)$. When $f$ is strongly convex, $x^*$ is unique.
For $m\ge 1$, let $\mathcal{F}^m $ be the set of functions that are convex and $k$ times continuously differentiable on $\bbR^n$. Moreover, we define $\mathcal{F}^m_L \subset \mathcal{F}^m $ which further requires that $f$ is $L$-smooth, that is,  the gradient $\nabla f$ is $L$-Lipschitz continuous, $\|\nabla f(y) - \nabla f(x)\| \leq L\|y - x\|$ for $x, y \in \bbR^n$, where $L>0$ is the Lipschitz constant and $\|\cdot\|$ denotes the Euclidean norm. Similarly, for $\mu>0$ we define $\mathcal{S}^m_\mu \subset\mathcal{F}^m$ which further requires that $f$ is $\mu$-strongly convex, that is, $f(y)\geq f(x) + \langle \nabla f(x),y-x\rangle + \frac{\mu}{2}\|y-x\|^2$ for $x, y \in \bbR^n$. Let $\mathcal{S}^m_{\mu,L}=\mathcal{F}^m_L\cap \mathcal{S}^m_{\mu}$. For additional properties of smoothness and strong convexity, readers are referred to, for example, Appendix A in \cite{aspremont2021acceleration}.

Given two sequences $\{a_k\} $ and $\{b_k\} $, %with $k$ indexing the algorithm iterations,
we write $a_k=O (b_k)$ if there exist constants $C>0$ and $K\ge 1$,
such that $a_k\leq C b_k$ for $k\geq K$.
Similarly, we write $a_k=\Omega (b_k)$ if $a_k\geq Cb_k$ for $k\geq K$.
We write $a_k=\Theta(b_k)$ if both $a_k=\Omega(b_k)$ and $a_k=O(b_k)$,
and write $a_k \sim b_k $ if $\lim_{k} a_k/b_k =1$.
%Similar notations are used for $g(t)$ and $h(t)$ being functions of continuous $t$.
More generally, for a set of pairs $(g,h)$, we write $g = O(h)$ or $g \lesssim h$ if there exists a constant $C>0$
such that $g \le C h$ for $(g,h)$.
Similarly, we write $g = \Omega (h)$ or $g \gtrsim h$ if $g \ge C h$ for $(g,h)$.
We write $g \asymp h$ if both $g \lesssim h$ and $g \gtrsim h$.
The constant $C$ may depend only on how the set of $(g,h)$ is defined from the context. For example,
for $a_k=O (b_k)$ above, $(g,h)$ represents $(a_k, b_k)$ for $k \ge K$, and
$C$ and $K$ may depend only on the sequences $\{a_k\} $ and $\{b_k\} $.
For NAG-C (\ref{eq:NAG-C}), the bound (\ref{eq:optimal bound-c}) indicates that there exists a constant $C>0$ such that
$f(x_k)-f^* \le C \frac{\|x_0-x^*\|^2}{s k^2}$ for any $L>0$, $f \in \mathcal{F}^1_L$, $x_0\in \bbR^n$, $0 < s \le 1/L$, and $k\ge 1$.
For Theorem~\ref{thm:converge-c}, the bound (\ref{eq:optimal bound-c}) is restricted to $0 < s \le C_0/L$ and $k\ge K$.
The constants $C$, $C_0$, and $K$ depend only on the algorithm parameters.

\section{Acceleration for strongly convex functions}\label{sec:acc-sc}

\subsection{Review of Nesterov’s acceleration}

A basic method for solving (\ref{eq:main_optimization}) is the gradient descent (GD):
\begin{equation}\label{eq:GD}
    x_{k+1} = x_k - s \nabla f(x_k),
\end{equation}
with an initial point $x_0$ and a stepsize $s$.
In this section, we study the case where $f\in \mathcal{S}^1_{\mu,L}$ is smooth and strongly convex. See Section \ref{sec:acc-c} for the case
where $f$ is smooth and convex.

For $f\in \mathcal{S}^1_{\mu,L}$, GD can be exponentially convergent in $k$, but the convergence rate depends on the condition number, namely $L/\mu$.
In fact, the GD iterates (\ref{eq:GD}) satisfy that  for $0<s\leq 2/(\mu+L)$ \citep{nesterov2018lectures},
$$
f(x_k)-f^* \leq \frac{L\|x_0-x^*\|^2}{2}\left(1-\frac{2\mu s}{1+\mu/L}\right)^k.
$$
When $s=1/L$ (in this case, the value of $\mu$ is not required for implementing GD), the rate is $(\frac{1-\mu/L}{1+\mu/L})^k\sim (1-\frac{\mu}{L})^{2k}$ for $\mu/L\to 0$. When $s=2/(L+\mu)$ (in this case, the value of $\mu$ is required), the rate improves to $(\frac{1-\mu/L}{1+\mu/L})^{2k}\sim(1-\frac{\mu}{L})^{4k}$. For both choices of the stepsize, the iteration complexity for $f(x_k)-f^*\leq \epsilon$ is $O(\frac{L}{\mu}\log(\frac{1}{\epsilon}))$, which has a linear dependency on $L/\mu$.

%The vanilla GD is not optimal in terms of the dependency of its convergence rate on $L/\mu$ for minimizing $f\in\mathcal{S}^1_{\mu,L} $.
\cite{nesterov1988approach} proposed an accelerated gradient method (NAG-SC) of the following form with an additional extrapolation step:
\begin{subequations}\label{eq:NAG-SC}
    \begin{align}
    y_{k+1} &= x_k - s \nabla f(x_k), \label{subeq:NAG-SC-y}\\
    x_{k+1} &= y_{k+1} + \sigma(y_{k+1}-y_k), \label{subeq:NAG-SC-single}
    \end{align}
\end{subequations}
with $x_0=y_0$ and the momentum coefficient $\sigma = \frac{1-\sqrt{q}}{1+\sqrt{q}}$, where $q =\mu s$, a shorthand to be used throughout this paper.
Equivalently, (\ref{eq:NAG-SC}) can be expressed in a single-variable form:
\begin{equation}\label{eq:NAG-SC-single}
    x_{k+1} = \underbrace{\vphantom{\frac{1-\sqrt{q}}{1+\sqrt{q}} }x_k - s \nabla f(x_k)}_{\text{gradient descent}} + \underbrace{\frac{1-\sqrt{q}}{1+\sqrt{q}}  (x_{k}-x_{k-1})}_{\text{momentum}}  - \underbrace{\frac{1-\sqrt{q}}{1+\sqrt{q}} \cdot s(\nabla f(x_{k})-\nabla f(x_{k-1}))}_{\text{gradient correction}},
\end{equation}
with $x_0$ and $x_1 = x_0 - \frac{2s\nabla f(x_0)}{1+\sqrt{q}}$.
Compared with GD, the iterate (\ref{eq:NAG-SC-single}) involves two additional terms, called the momentum and gradient correction.
For $0<s\leq 1/L$, the NAG-SC iterates satisfy $f(x_k)-f^* = O((1-\sqrt{q})^k)$.
When $s=1/L$, the bound reduces to $O((1-\sqrt{\frac{\mu}{L}})^k)$, and the iteration complexity is lowered to $O(\sqrt{\frac{L}{\mu}}\log(\frac{1}{\epsilon}))$, with a square-root dependency on $L/\mu$. \cite{drori2022oracle} established that a lower bound for minimizing general $f\in\mathcal{S}^1_{\mu, L}$ is $f(x_k)-f^*=\Omega((1-\sqrt{\frac{\mu}{L}})^{2k})$, where $\{x_k\}$ are iterates of any black-box\footnote{Black-box means that no prior knowledge of $f$ (e.g., $f$ is quadratic) is available except for the class $f$ belongs to. For $\mathcal{F}^1_L$, the available information is only the Lipschitz constant $L$.} first-order method. Therefore, NAG-SC is optimal up to a constant factor of 2 in terms of the iteration complexity.

There are several first-order methods \textit{exactly} reaching the lower bound $(1-\sqrt{\frac{\mu}{L}})^{2k}$, for instance, the information-theoretic exact method (ITEM) \citep{taylor2022optimal} and triple-momentum method (TMM) \citep{van2017fastest}.
While ITEM involves time-dependent coefficients, TMM is defined with time-independent coefficients \citep{aspremont2021acceleration}:
\begin{subequations}\label{eq:TMM}
    \begin{align}
        y_{k+1} &= x_{k} - s \nabla f(x_{k}),\label{eq:TMM-y}\\
        z_{k+1} &= \sqrt{q}\left(x_k-\frac{1}{\mu}\nabla f(x_k)\right) + (1- \sqrt{q}) z_k,\label{eq:TMM-z} \\
        x_{k+1} &= \frac{2 \sqrt{q}}{1+\sqrt{q}} z_{k+1} + \left(1-\frac{2 \sqrt{q}}{1+\sqrt{q}}\right)y_{k+1}, \label{eq:TMM-x}
    \end{align}
\end{subequations}
with $x_0=z_0$.
See Supplement \ref{subsec:reformulation of NAG-SC and TMM}
for an equivalent form of (\ref{eq:TMM}), where $(y_{k+1}, x_{k+1})$ are defined independently of $\{ z_k \}$, and then
$z_{k+1}$ is defined as
\begin{align}
z_{k+1} = \frac{1+\sqrt{q}}{2\sqrt{q}} x_{k+1} + \left(1-\frac{1+\sqrt{q}}{2\sqrt{q}}\right)y_{k+1},\label{eq:TMM-z-auxiliary}
\end{align}
which is of a similar form as (\ref{eq:NAG-sc-z}) below for NAG-SC.
%a weighted average of $x_{k+1}$ and $y_{k+1}$, with weights $ (1+\sqrt{q})/(2\sqrt{q})$ and $(-1+\sqrt{q})/(2 \sqrt{q})$.
In this sense, the sequence $\{z_k\}$ in TMM is auxiliary.
Nevertheless, $\{z_k\}$ plays a vital role in the existing analysis of TMM.\
\cite{van2017fastest} showed that $\{ z_k\}$ achieves the lower bound, i.e.,
$f(z_k)-f^*=O((1-\sqrt{\frac{\mu}{L}})^{2k})$ when $s=1/L$.
It seems to be an open question whether $\{ x_k\}$ also achieves the lower bound.

In the Lyapunov analysis of NAG-SC in \cite{bansal2019potential}, a similar auxiliary sequence $\{z_k\}$ as above is introduced by defining
\begin{equation}
    z_{k+1} = \frac{1+\sqrt{q}}{\sqrt{q}} x_{k+1} + \left(1-\frac{1+\sqrt{q}}{\sqrt{q}}\right)y_{k+1}, \label{eq:NAG-sc-z}
\end{equation}
and the NAG-SC iterates in (\ref{eq:NAG-SC}) together with (\ref{eq:NAG-sc-z}) are equivalently reformulated as
\begin{subequations}\label{eq:reformulate-NAG-SC}
    \begin{align}
        y_{k+1} &= x_{k} - s \nabla f(x_{k}),\label{subeq:reformulate-NAG-SC-y}\\
        z_{k+1} &= \sqrt{q}\left(x_k-\frac{1}{\mu}\nabla f(x_k)\right) + (1- \sqrt{q}) z_k, \label{subeq:reformulate-NAG-SC-z}\\
        x_{k+1} &= \frac{ \sqrt{q}}{1+\sqrt{q}} z_{k+1} + \left(1-\frac{ \sqrt{q}}{1+\sqrt{q}}\right)y_{k+1}, \label{subeq:reformulate-NAG-SC-x}
    \end{align}
\end{subequations}
with $x_0=z_0$.
See Supplement \ref{subsec:reformulation of NAG-SC and TMM} for a derivation, provided for completeness.
Note that TMM and NAG-SC differ only in how $x_{k+1}$ is defined in (\ref{eq:TMM-x}) and
(\ref{subeq:reformulate-NAG-SC-x}).
\cite{bansal2019potential} constructed a new Lyapunov function to simplify the convergence proof for NAG-SC,
but only directly showed that $\{ y_k\}$
achieves the bound $f(y_k)-f^*=O((1-\sqrt{\frac{\mu}{L}})^k)$ when $s=1/L$.
By some additional arguments, similar convergence bounds can also be deduced for $\{z_k\}$ and $\{ x_k\}$.

For comparison, we also mention the heavy-ball method \citep{polyak1964some}, defined as\footnote{
Following \cite{shi2021understanding},
algorithm (\ref{eq:HB}) is a slight modification of the original method in \cite{polyak1964some} where $s=4/(\sqrt{L}+\sqrt{\mu})^2$ and the momentum coefficient $\sigma = (1-\sqrt{q})^2$. If $s$ is small, the two coefficients $(1-\sqrt{q})^2$ and $\frac{1-\sqrt{q}}{1+\sqrt{q}}$ are close.
The original heavy-ball method, with the specific $s$, achieves the accelerated convergence rate $(\frac{1-\sqrt{\mu/L}}{1+\sqrt{\mu/L}})^{2k}$
for $f\in\mathcal{S}^2_{\mu,L}$ \citep{polyak1987intro}, 
but it may fail to converge for some $f\in\mathcal{S}^1_{\mu,L}$ \citep{lessard2016analysis}.}
\begin{equation}\label{eq:HB}
    x_{k+1} = x_k - s\nabla f(x_k) + \sigma (x_k-x_{k-1}) ,
\end{equation}
with $\sigma = \frac{1-\sqrt{q}}{1+\sqrt{q}}$ as in (\ref{eq:NAG-SC}). Note that (\ref{eq:HB})  differs from NAG-SC (\ref{eq:NAG-SC-single}) only in the absence of gradient correction.
The heavy-ball method may not achieve accelerated convergence for general $f\in \mathcal{S}^1_{\mu,L}$.
In fact, the iterates (\ref{eq:HB}) satisfy $f(x_k)-f^* = O((1-C_1 \sqrt{q})^k)$, but only provably for the stepsize $0<s \le C_0\mu/L^2$,
with some  constants $C_0,C_1>0$ \citep{shi2021understanding}, where $\mu/L^2$ can be much smaller than $1/L$.
Hence the convergence bound for the heavy-ball method is still
of order $O\left((1-C\frac{\mu}{L})^k\right)$ for some constant $C>0$, corresponding to the same iteration complexity as GD.

\subsection{Main results}

We formulate a broad class of algorithms including NAG-SC and TMM as special cases and
establish (simple and interpretable) sufficient conditions for when the algorithms in the class achieve Nesterov's acceleration,
which is defined as reaching an objective gap of $O((1-C\sqrt{\frac{\mu}{L}})^k)$ at iteration $k$,
for a constant $C>0$ and a suitable stepsize $s$, usually in the order $s \asymp 1/L$.
Our work does not aim to find a sharp value of $C$ or
address the question of whether these algorithms \textit{exactly} achieve the complexity lower bound corresponding to $C=2$ \citep{drori2022oracle}. In Supplement Section \ref{sec:experiment}, we present numerical results to illustrate different performances from specific algorithms
where our sufficient conditions are either satisfied or violated.

To unify and extend NAG-SC and TMM, we consider the following class of algorithms:
\begin{subequations}\label{eq:extended-NAG-SC}
    \begin{align}
        y_{k+1} &= x_{k} - \eta s \nabla f(x_{k}),\label{subeq:extended-NAG-SC-y}\\
        z_{k+1} &= \nu \sqrt{q}\left(x_k-\frac{1}{\mu}\nabla f(x_k)\right) + (1-\nu\sqrt{q}) z_k, \label{subeq:extended-NAG-SC-z}\\
        x_{k+1} &= \frac{\tau \sqrt{q}}{1+\sqrt{q}} z_{k+1} + \left(1-\frac{\tau \sqrt{q}}{1+\sqrt{q}}\right)y_{k+1},\label{subeq:extended-NAG-SC-x}
    \end{align}
\end{subequations}
with $x_0=z_0$, where $\eta, \nu, \tau \ge 0$ are three parameters which may depend on $q$.
Then NAG-SC (\ref{eq:NAG-SC}) or (\ref{eq:reformulate-NAG-SC}) is recovered by setting $(\eta,\nu,\tau)=(1,1,1)$, and TMM (\ref{eq:TMM}) is recovered by $(\eta,\nu,\tau)=(1,1,2)$. In this way, the two algorithms differ only in the choice of $\tau$.
Equivalently, (\ref{eq:extended-NAG-SC}) can be put into a single-variable form in terms of $\{x_k\}$, similarly to (\ref{eq:NAG-SC-single}) for NAG-SC:
\begin{equation}\label{eq:supp-extended-NAG-SC-single}
    \begin{split}
        x_{k+1} &= x_k - \frac{\nu(\tau+\zeta\eta\sqrt{q})}{1+\sqrt{q}} s\nabla f(x_k) + \frac{\zeta(1-\nu\sqrt{q})}{1+\sqrt{q}}(x_k-x_{k-1}) \\
        &\quad - \frac{\zeta\eta(1-\nu\sqrt{q})}{1+\sqrt{q}} s(\nabla f(x_k)-\nabla f(x_{k-1})),
    \end{split}
\end{equation}
with $x_0$ and $x_1=x_0-\frac{\zeta\eta+\nu\tau}{1+\sqrt{q}}s\nabla f(x_0)$, where $\zeta = 1+(1-\tau)\sqrt{q}$, a shorthand throughout this paper.
The coefficients in the three terms for $ s\nabla f(x_k) $, $x_k-x_{k-1}$, and $s(\nabla f(x_k)-\nabla f(x_{k-1}))$
are implicitly constrained due to the translation from (\ref{eq:extended-NAG-SC}),
as discussed after Corollary~\ref{cor:symmetrized conditions for extended-NAG-SC-single}.
See also Lemma~\ref{lem:matching of two forms of extended-NAG-SC} for a partial converse from a single-variable form $\{x_k\}$ to the three-variable form (\ref{eq:extended-NAG-SC})
under some conditions on the coefficients in the single-variable form.

The following result gives sufficient conditions for convergence of algorithm (\ref{eq:extended-NAG-SC})
in the scenario where  $\eta$, $\nu$, and $\tau$ are absolute constants, free of $q$. In all our results, the constants $C_0$ and $C_1$ can be explicitly calculated given algorithm parameters
by our proofs in the Supplement.

\begin{thm}\label{thm:converge-sc-const}
     Let $f\in \mathcal{S}^1_{\mu,L}$. Assume that $\eta=\eta_0$, $\nu=\nu_0$, and $\tau=\tau_0$ for some constants $\eta_0$, $\nu_0$, and $\tau_0$, free of $q$.
     (i) There exist constants $C_0,C_1>0$, depending only on $(\eta_0,\nu_0,\tau_0)$, such that for $0<s \le C_0\frac{\mu}{L^2}$ and $k\ge 1$, the iterates of (\ref{eq:extended-NAG-SC}) satisfy
     \begin{equation}\label{eq:exponential bound in q}
        f(x_k)-f^* = O \left( L\|x_0-x^*\|^2(1-C_1 \sqrt{\mu s})^k  \right),
    \end{equation}
    provided that one of the following conditions holds: \vspace{-.05in}
        \begin{itemize} \addtolength{\itemsep}{-.05in}
        \item[(ia)] $\nu_0,\tau_0>0$, $\nu_0\not=\tau_0$, and $0 \le \eta_0 < \nu_0\tau_0/2$;
        \item[(ib)] $\nu_0=\tau_0 >2$ and $\eta_0 = \tau_0^2/2$.
    \end{itemize} \vspace{-.05in}
    (ii) There exist constants $C_0,C_1>0$, depending only on $(\eta_0,\nu_0,\tau_0)$, such that for $0<s \le C_0 \frac{1}{L}$, the iterates of (\ref{eq:extended-NAG-SC}) satisfy (\ref{eq:exponential bound in q}), provided that one of the following conditions holds: \vspace{-.05in}
        \begin{itemize} \addtolength{\itemsep}{-.05in}
        \item[(iia)] $\nu_0,\tau_0>0$, $\nu_0\not=\tau_0$, and $\eta_0\geq\nu_0\tau_0/2$;
        \item[(iib)] $\nu_0= \tau_0 \geq 2$ and $\eta_0 > \tau_0^2/2$;
        \item[(iic)] $1<\nu_0= \tau_0 <2$ and $\eta_0 > \tau_0$;
        \item[(iid)] $0<\nu_0= \tau_0 \leq 1$ and $\eta_0\geq\tau_0$.
    \end{itemize} \vspace{-.05in}

\end{thm}

The convergence bound (\ref{eq:exponential bound in q}) exhibits two types of dependency on $\mu/L$,
determined by how large the stepsize $s$ is allowed.
As mentioned earlier, a similar phenomenon occurs in the comparison of the heavy-ball and NAG-SC methods.
Under the conditions in Theorem~\ref{thm:converge-sc-const}(i), algorithm (\ref{eq:extended-NAG-SC}) achieves a usual convergence bound
$O((1-C \frac{\mu}{L})^k)$ for $s \asymp \frac{\mu}{L^2}$, resulting in an iteration complexity $O(\frac{L}{\mu}\log(\frac{1}{\epsilon}))$
with a linear dependency on $L/\mu$.
Under the conditions in Theorem~\ref{thm:converge-sc-const}(ii),
algorithm (\ref{eq:extended-NAG-SC}) reaches an accelerated convergence bound
$O((1-C \sqrt{\frac{\mu}{L}})^k)$ for $s \asymp \frac{1}{L}$, resulting in an iteration complexity $O(\sqrt{\frac{L}{\mu}} \log(\frac{1}{\epsilon}))$
with a square-root dependency on $L/\mu$.
In particular, TMM with $(\eta_0,\nu_0,\tau_0)=(1,1,2)$ is covered by condition (iia),
whereas NAG-SC with $(\eta_0,\nu_0,\tau_0)=(1,1,1)$ is covered by condition (iid).
It is interesting that for both TMM and NAG-SC, the choice of $(\eta_0,\nu_0,\tau_0)$ lies on the boundaries of the acceleration regions
identified in Theorem~\ref{thm:converge-sc-const}(ii).\footnote{The TMM choice, $\eta_0=\nu_0=1$ and $\tau_0=2$, satisfies $\nu_0 \not= \tau_0$ and $\eta_0 = \nu_0\tau_0/2$, lying on the boundary of condition (iia) in Theorem~\ref{thm:converge-sc-const}. The NAG-SC choice, $\eta_0 =\nu_0=\tau_0=1$, satisfies $\nu_0= \tau_0= 1$ and $\eta_0 = \tau_0$, lying on the boundary of condition (iid) in Theorem~\ref{thm:converge-sc-const}.}

Next, we study convergence of algorithm (\ref{eq:extended-NAG-SC}) in a more general scenario
where $\eta$, $\nu$ and $\tau$ are analytical functions of $\sqrt{q}$, free of negative exponents like $q^{-1/2}$,
as stated in Assumption \ref{ass:regular func of sqrt-q}.
In particular, $\eta$, $\nu$ and $\tau$ which are polynomials of $\sqrt{q}$ satisfy Assumption~\ref{ass:regular func of sqrt-q}.

\begin{ass}\label{ass:regular func of sqrt-q}
    There exist non-negative, analytic functions $\tilde \eta(\cdot)$, $\tilde \nu(\cdot)$ and $\tilde \tau(\cdot)$ in a neighborhood of $0$ such that $\eta = \tilde \eta(\sqrt{q})$, $\nu = \tilde \nu(\sqrt{q})$ and $\tau=\tilde \tau(\sqrt{q})$. Then $\eta$, $\nu$ and $\tau$ admit (convergent) Taylor expansions:
    \begin{equation} \label{eq:param-expan}
        \eta = \sum_{i=0}^\infty \eta_i (\sqrt{q})^{i}, \quad \nu = \sum_{i=0}^\infty \nu_i (\sqrt{q})^{i}, \quad \tau = \sum_{i=0}^\infty \tau_i (\sqrt{q})^{i},
    \end{equation}
    when $0\leq q\leq q_0$ for some constant $q_0>0$.
\end{ass}

The following result gives sufficient conditions on convergence of algorithm (\ref{eq:extended-NAG-SC}), in terms of only the constant coefficients,
$(\eta_0,\nu_0,\tau_0)$, in the expansions (\ref{eq:param-expan}).

\begin{thm}\label{thm:converge-sc-0}
    Suppose that $f\in \mathcal{S}^1_{\mu,L}$ and Assumption~\ref{ass:regular func of sqrt-q} holds.
    (i) For $0<s \le C_0 \frac{\mu}{L^2}$ and $k\ge 1$, the iterates of (\ref{eq:extended-NAG-SC}) satisfy (\ref{eq:exponential bound in q}),
    with $C_0,C_1 >0$ depending only on $(\tilde \eta, \tilde \nu, \tilde \tau)$,
    provided that condition (ia) in Theorem \ref{thm:converge-sc-const} holds.
    (ii) For $0<s \le C_0 \frac{1}{L}$ and $k\ge 1$, the iterates of (\ref{eq:extended-NAG-SC}) satisfy (\ref{eq:exponential bound in q}),
    with $C_0,C_1>0$ depending only on $(\tilde \eta, \tilde \nu, \tilde \tau)$, provided that condition (iia) in Theorem \ref{thm:converge-sc-const} holds.
\end{thm}

Theorem \ref{thm:converge-sc-0} can be understood similarly as Theorem \ref{thm:converge-sc-const}.
However, we stress that Theorems~\ref{thm:converge-sc-const} and \ref{thm:converge-sc-0} are partially overlapped, but distinct from each other.
The conclusions from Theorem~\ref{thm:converge-sc-const} under condition (ia) or (iia) can be obtained as special cases
of Theorem~\ref{thm:converge-sc-0}, where $\eta$, $\nu$, and $\tau$ are constant in $q$, i.e.,
$\eta_i = \nu_i=\tau_i =0$ for $i \ge 1$.
But the conclusions from  Theorem~\ref{thm:converge-sc-const} under condition (ib) or (iib)--(iid) are not covered by Theorem \ref{thm:converge-sc-0},
and the conclusions from Theorem \ref{thm:converge-sc-0} in the case of $\eta$, $\nu$, and $\tau$ depending on $q$ are
not covered by Theorem \ref{thm:converge-sc-const}.
For example, NAG-SC is covered only by Theorem~\ref{thm:converge-sc-const},
whereas TMM is covered by both Theorems~\ref{thm:converge-sc-const} and \ref{thm:converge-sc-0}.

In the case of
$\nu_0 =\tau_0 >0$ and $\eta_0 \ge \tau_0^2/2$ for $\eta$, $\nu$, and $\tau$ depending on $q$, which is not addressed by Theorem \ref{thm:converge-sc-0},
the following result gives sufficient conditions on convergence of algorithm (\ref{eq:extended-NAG-SC}),
involving the linear coefficients $(\eta_1,\nu_1,\tau_1)$ in the expansions (\ref{eq:param-expan}).

\begin{thm}\label{thm:converge-sc-1}
    Suppose that $f\in \mathcal{S}^1_{\mu,L}$ and Assumption~\ref{ass:regular func of sqrt-q} holds with $0<\nu_0=\tau_0$ and $\eta_0\geq \tau_0^2/2$.
    (i) For $0<s \le C_0 \frac{\mu}{L^2}$ and $k\ge 1$, the iterates of (\ref{eq:extended-NAG-SC}) satisfy (\ref{eq:exponential bound in q}), with constants $C_0,C_1>0$ depending only on $(\tilde \eta, \tilde \nu, \tilde \tau)$,  provided that the following holds:
       \begin{itemize}\addtolength{\itemsep}{-.05in}
        \item[(ia)] $\eta_0= \tau_0^2/2$, $\nu_1-\tau_1<\tau_0\left(\frac{\tau_0}{2}-1\right)$, and $2\eta_1<\nu_1\tau_1 + \frac{\tau_0^2}{2}\left(\frac{5}{2}\tau_0-2\right)$.
    \end{itemize} \vspace{-.05in}
    (ii) For $0<s \le C_0 \frac{1}{L}$ and $k\ge 1$,  the iterates of (\ref{eq:extended-NAG-SC}) satisfy (\ref{eq:exponential bound in q}), with constants $C_0, C_1>0$ depending only on $(\tilde \eta, \tilde \nu, \tilde \tau)$,  provided that one of the following conditions holds: \vspace{-.05in}
    \begin{itemize}\addtolength{\itemsep}{-.05in}
        \item[(iia)] $\eta_0=\tau_0^2/2$, $\nu_1-\tau_1<\tau_0\left(\frac{\tau_0}{2}-1\right)$, and $2\eta_1\geq \nu_1\tau_1 + \frac{\tau_0^2}{2}\left(\frac{5}{2}\tau_0-2\right)$;
        \item[(iib)] $\eta_0>\tau_0^2/2$ and $\nu_1-\tau_1<(\tau_0-1)\tau_0-\eta_0$;
        \item[(iic)] $\eta_0>\tau_0^2/2$ and $(\tau_0-1)\tau_0-\eta_0 < \nu_1-\tau_1<\eta_0-\tau_0$.
    \end{itemize} \vspace{-.05in}
\end{thm}

Similarly, Theorems \ref{thm:converge-sc-const} and \ref{thm:converge-sc-1} are also partially overlapped, but distinct from each other.
Theorem \ref{thm:converge-sc-1} deals with more general parameters $\eta$, $\nu$, and $\tau$, possibly depending on $q$.
But in the case of constant parameters (hence $\eta_1=\nu_1=\tau_1=0$), the sufficient conditions in Theorem \ref{thm:converge-sc-1}
is more restrictive than those in Theorem \ref{thm:converge-sc-const}, i.e., if
a sufficient condition stated in Theorem \ref{thm:converge-sc-1} holds, then
at least one of the sufficient conditions in Theorem \ref{thm:converge-sc-const} must also hold.
This is because Theorem \ref{thm:converge-sc-const} is deduced by exploiting the constant assumption on $(\eta,\nu,\tau)$,
which provides information on all coefficients
$(\eta_i,\nu_i,\tau_i)$ for $i \ge 0$, not just for $i=0, 1$, in the expansions (\ref{eq:param-expan}).

All the preceding results are applicable to algorithm (\ref{eq:extended-NAG-SC}) or equivalently
its single-variable form (\ref{eq:supp-extended-NAG-SC-single}),
but the coefficients in (\ref{eq:supp-extended-NAG-SC-single}) appear to be highly structured.
To facilitate comparison and interpretation (see Sections \ref{sec:ODE-interpretation} and \ref{sec:HAG-interpretation}), we translate the conclusion of Theorem \ref{thm:converge-sc-0} in the case of $\nu_0\not=\tau_0$ directly
in terms of the leading coefficients in a single-variable form.
The parameterization of (\ref{eq:extended-NAG-SC-single-2}) in terms of $(c_1,c_2,c_3)$ is motivated by the HAG algorithm,
as shown in (\ref{eq:extended-NAG-SC-single-2b}).

\begin{cor}\label{cor:symmetrized conditions for extended-NAG-SC-single}
    For $f\in \mathcal{S}^1_{\mu,L}$, consider the following algorithm:
    \begin{equation}\label{eq:extended-NAG-SC-single-2}
        \begin{split}
            x_{k+1} &= x_k - (c_0+R_1) s\nabla f(x_k)  + (1-c_1\sqrt{q} + R_2)(x_k-x_{k-1}) \\
            &\quad -(c_2\sqrt{c_0}-\frac{c_0}{2}+R_3)s(\nabla f(x_k)-\nabla f(x_{k-1})),
        \end{split}
    \end{equation}
    with $x_1=x_0-h_1 s\nabla f(x_0)$, where $c_0,c_1,c_2>0$ and
    $R_1=O(\sqrt{q})$, $R_2=O(q)$, $R_3=O(\sqrt{q})$ and $h_1$ are analytic functions of $\sqrt{q}$ around $0$.
    (i) For $0<s \le C_0 \mu/L^2$ and $k\ge 1$,
    the iterates of (\ref{eq:extended-NAG-SC-single-2}) satisfy (\ref{eq:exponential bound in q}), with $C_0,C_1>0$ depending only on $(c_0,c_1,c_2)$
    and $(R_1,R_2,R_3)$, provided $c_0/4\leq c_2^2< c_0$ and $c_1^2>4c_0$.
    (ii) For $0<s \le C_0 /L$ and $k\ge 1$,
    the iterates of (\ref{eq:extended-NAG-SC-single-2}) satisfy (\ref{eq:exponential bound in q}), with $C_0,C_1>0$ depending only on $(c_0,c_1,c_2)$
    and $(R_1,R_2,R_3)$,
    provided $c_2^2\geq c_0$ and $c_1^2>4c_0$.
\end{cor}

%The conclusions (i) and (ii) in Corollary \ref{cor:symmetrized conditions for extended-NAG-SC-single} are translated from
%(i) and (ii) in Theorem \ref{thm:converge-sc-0} respectively.

The single-variable form (\ref{eq:extended-NAG-SC-single-2}) is more general than induced by (\ref{eq:extended-NAG-SC}).
For example, unlike TMM and NAG-SC, the heavy-ball method (\ref{eq:HB}) does not fall in the class of algorithms (\ref{eq:extended-NAG-SC}).
But all the three methods can be put into (\ref{eq:extended-NAG-SC-single-2}) as follows. \vspace{-.05in}
\begin{itemize} \addtolength{\itemsep}{-.05in}
\item TMM:   $(c_0,c_1,c_2) = (2,3,\sqrt{2})$ with $c_1^2>4c_0$;

\item NAG-SC: $(c_0,c_1,c_2) = (1,2, 3/2)$ with $c_1^2=4c_0$;

\item Heavy-ball: $(c_0,c_1,c_2) = (1,2, 1/2)$ with $c_1^2=4c_0$.
\end{itemize}
As mentioned earlier, TMM is covered by both Theorems~\ref{thm:converge-sc-const} and \ref{thm:converge-sc-0} (and hence
Corollary \ref{cor:symmetrized conditions for extended-NAG-SC-single}),
whereas NAG-SC is covered by Theorem~\ref{thm:converge-sc-const}
but not Theorem \ref{thm:converge-sc-0} or Corollary \ref{cor:symmetrized conditions for extended-NAG-SC-single}.

The translation from Theorem \ref{thm:converge-sc-0} about  (\ref{eq:extended-NAG-SC}) to
Corollary \ref{cor:symmetrized conditions for extended-NAG-SC-single} about (\ref{eq:extended-NAG-SC-single-2})
relies on the fact that if $c_1^2 > 4c_0$ and $c_2^2 \ge c_0/4$, then
the iterates $\{x_k\}$ from the single-variable form (\ref{eq:extended-NAG-SC-single-2}) can be equivalently put into the three-variable form (\ref{eq:extended-NAG-SC}) for some parameters $\eta$, $\nu$, and $\tau$,
which satisfy $\nu_0,\tau_0>0$, $\nu_0\not=\tau_0$, and $\eta_0\geq 0$ in the expansions (\ref{eq:param-expan}). See Lemma \ref{lem:matching of two forms of extended-NAG-SC} in the Supplement.
By similar reasoning as in the proof of Lemma \ref{lem:matching of two forms of extended-NAG-SC},
when $c_1^2=4c_0$, (\ref{eq:extended-NAG-SC-single-2}) may still be put into (\ref{eq:extended-NAG-SC}) with $\nu_0=\tau_0$,
depending on additional information about  $R_1$, $R_2$ and $R_3$.
In that case, convergence of (\ref{eq:extended-NAG-SC-single-2}) can be deduced from Theorem \ref{thm:converge-sc-1}, where
 $\nu_0=\tau_0$ and the linear coefficients $(\eta_1,\nu_1,\tau_1)$ are involved,
or from Theorem~\ref{thm:converge-sc-const} with $\nu_0=\tau_0$ under the assumption that
$\eta=\eta_0$, $\nu=\nu_0$, and $\tau=\tau_0$ are constants, free of $q$.
However, when $c_1^2<4c_0$, (\ref{eq:extended-NAG-SC-single-2}) can no longer be put into (\ref{eq:extended-NAG-SC}).
This reveals an inherent constraint which is satisfied by the single-variable form (\ref{eq:supp-extended-NAG-SC-single})
derived from algorithm (\ref{eq:extended-NAG-SC}).

\section{Acceleration for general convex functions}\label{sec:acc-c}

\subsection{Review of Nesterov's acceleration}
Consider the unconstrained minimization problem (\ref{eq:main_optimization}) with $f\in\mathcal{F}^1_L$.
For stepsize $0<s\leq 1/L$, the GD iterates $\{f(x_k)\}$ is non-increasing and satisfy $f(x_k)-f^*\leq \|x_0-x^*\|^2/(2ks)$.\footnote{\label{fn:GD with large s}For stepsize $1/L <s< 2/L$, the GD iterates $\{f(x_k)\}$ is still non-increasing but satisfies a slightly different bound $O\left(\frac{\|x_0-x^*\|^2}{ks(2-Ls)}\right)$ for the objective gap \citep{nesterov2018lectures}.} However, the $O(1/k)$ rate is not optimal.
\cite{nesterov1988approach} proposed an accelerated gradient method in a similar form as NAG-SC (\ref{eq:NAG-SC}):
\begin{equation}\label{eq:NAG-C}
    \begin{split}
        y_{k+1} &= x_k - s \nabla f(x_k), \\
        x_{k+1} &= y_{k+1} + \sigma_{k+1} (y_{k+1}-y_k),
    \end{split}
\end{equation}
with $x_0=y_0$ and $\sigma_{k+1}=k/(k+3)$.
Compared with NAG-SC (\ref{eq:NAG-SC}), the momentum coefficient $\sigma_{k+1}$ varies with $k$, instead of a fixed value depending on the strong-convexity parameter $\mu$.
For $0<s\leq 1/L$, the NAG-C iterates $\{f(x_k)\}$ may not be non-increasing but satisfy
\begin{equation}\label{eq:optimal bound-c}
    f(x_k)-f^*=O\left(\frac{\|x_0-x^*\|^2}{s k^2}\right).
\end{equation}
Compared with GD, the iteration complexity for $f(x_k)-f^*\leq \epsilon$ is reduced from $O(L/\epsilon)$ to $O(\sqrt{L/\epsilon})$ when taking $s=1/L$. Notably, the convergence bound in (\ref{eq:optimal bound-c}) matches with the lower bound of black-box first-order methods for minimizing functions in $\mathcal{F}^1_L$ when $n$ (the dimension of $x$) is relatively large compared with $k$ \citep{nesterov2018lectures}.

The single-variable form in terms of $\{x_k\}$ for NAG-C (\ref{eq:NAG-C}) is
\begin{equation}\label{eq:NAG-C-single}
    x_{k+1} = \underbrace{\vphantom{\frac{k}{k+3}}x_k - s \nabla f(x_k)}_{\text{gradient descent}} + \underbrace{\frac{k}{k+3} (x_{k}-x_{k-1})}_{\text{momentum}}  - \underbrace{\frac{k}{k+3}\cdot s(\nabla f(x_{k})-\nabla f(x_{k-1}))}_{\text{gradient correction}},
\end{equation}
with $x_0$ and $x_1=x_0 - s \nabla f(x_0)$.
The iterate (\ref{eq:NAG-C-single}) involves an additional momentum term and a gradient correction term
similarly as in (\ref{eq:NAG-SC-single}) for NAG-SC, but with time-dependent coefficients.
\cite{shi2021understanding} studied NAG-C by relating (\ref{eq:NAG-C-single}) to a high-resolution ODE (see Section \ref{sec:ODE-c}) and
obtained a new result on the squared gradient norm:
for stepsize $0<s \le  1/(3L)$,
\begin{equation}\label{eq:inverse cubic rate}
    \min_{0\leq i \leq k} \|\nabla f(x_i)\|^2 = O\left(\frac{\|x_0-x^*\|^2}{s^2 (k+1)^3} \right) .
\end{equation}
The inverse cubic rate (\ref{eq:inverse cubic rate}) cannot be obtained directly from (\ref{eq:optimal bound-c}).
Moreover, \cite{shi2021understanding} extended the coefficients in the momemtum and gradient correction terms in (\ref{eq:NAG-C-single}) to $k/(k+r+1)$ and
$(k/(k+r+1))\beta s$ for any $r \ge 2$ and $\beta>1/2$ and showed that (\ref{eq:optimal bound-c}) and (\ref{eq:inverse cubic rate}) remain valid.

The momentum coefficient $\sigma_{k+1}$ in (\ref{eq:NAG-C}) can be defined in more flexible forms than $k/(k+r+1)$ with $r\ge 2$
to achieve the optimal rate (\ref{eq:optimal bound-c}). A popular scheme is to set $\sigma_{k+1} = (\alpha_k-1)/\alpha_{k+1}$, where $\{\alpha_k\}$ is a scalar sequence to be chosen.
It is known \citep{beck2017first} that for any sequence $\{\alpha_k\}$ satisfying $\alpha_k=\Omega(k)$ and
\begin{equation}\label{eq:recursive condition}
  \alpha_{k+1}(\alpha_{k+1}-1)\leq \alpha^2_k,
\end{equation}
the corresponding algorithm (\ref{eq:NAG-C}) achieves the optimal bound (\ref{eq:optimal bound-c}). The recursive condition (\ref{eq:recursive condition}) implies that $\a_{k+1}\leq (1+\sqrt{1+4\alpha_k^2})/{2}\leq \a_k+1$, and hence $\a_k=\Theta(k)$.
For the accelerated proximal gradient method or FISTA
\citep{beck2009fast}, $\{\alpha_k\}$ is defined recursively as $\alpha_{k+1}=(1+\sqrt{1+4\alpha_{k}^2})/2$ with $\alpha_0=1$,
i.e., (\ref{eq:recursive condition}) holds as equality.
It can also be verified that $\{\a_k=(k+r)/r\}$ satisfies (\ref{eq:recursive condition}) for any $r\geq 2$,
with the corresponding $\sigma_{k+1} = k/(k+r+1)$.

\subsection{Main results}

We formulate a broad class of algorithms including NAG-C (\ref{eq:NAG-C}) and existing variations \citep{beck2017first,shi2021understanding}
and establish (simple and interpretable) sufficient conditions for when the algorithms in the class achieve
both the optimal bound (\ref{eq:optimal bound-c}) for the objective gap and the inverse cubic rate (\ref{eq:inverse cubic rate})
for the squared gradient norm, similarly as NAG-C. In Supplement Section \ref{sec:experiment}, we present numerical results to illustrate different performances from specific algorithms
where the algorithm parameters either satisfy or violate our sufficient conditions.

To unify and extend existing choices of the momentum and gradient correction terms related to NAG-C, we consider the following class of algorithms:
\begin{equation}\label{eq:extended-NAG-C}
    \begin{split}
        y_{k+1} &= x_k - \beta_k s \nabla f(x_k), \\
        x_{k+1} &= x_k - \gamma_k s \nabla f(x_k) + \sigma_{k+1}(y_{k+1}-y_k),
    \end{split}
\end{equation}
where $x_0=y_0$, $\sigma_{k+1}=(\alpha_k-1)/\alpha_{k+1}$, and $\{\alpha_k\}$, $\{\beta_k\}$ and $\{\gamma_k\}$ are three scalar sequences.
The equivalent single-variable form of (\ref{eq:extended-NAG-C}) is
\begin{equation}\label{eq:extended-NAG-C-single}
    \begin{split}
        x_{k+1} &= x_{k} - (\gamma_k+\sigma_{k+1}(\beta_k-\beta_{k-1})) s \nabla f(x_k) + \sigma_{k+1}(x_k-x_{k-1}) \\
        &\quad - \sigma_{k+1}\beta_{k-1} \cdot s(\nabla f(x_{k})-\nabla f(x_{k-1})),
    \end{split}
\end{equation}
starting from $x_0$ and $x_1 = x_0 - (\gamma_0+(\alpha_0-1)\beta_0/\alpha_1)s\nabla f(x_0)$.
The algorithms studied in \cite{beck2017first} correspond to the case of $\beta_k = \gamma_k=1$ and $\alpha_k$ of a general functional form
subject to certain inequality constraints.
The algorithms studied in \cite{shi2021understanding} correspond to $\gamma_k=1$, $\beta_k = \beta$ (a constant), and
$\alpha_k = (k+r)/r $ for $r \ge 2$.

Motivated by the Lyapunov analysis of NAG-C in \cite{su2016differential},
we introduce the following three-variable form of (\ref{eq:extended-NAG-C}).
Define $z_k = \alpha_k x_k + (1 - \alpha_k) y_k$ for $k\geq 0$, which is reminiscent of (\ref{eq:TMM-z}) and (\ref{eq:NAG-sc-z}) in the strongly convex setting.
Then (\ref{eq:extended-NAG-C}) can be equivalently reformulated as
   \begin{subequations}\label{eq:extended-NAG-C-xyz}
        \begin{align}
            y_{k+1} &= x_k - \beta_k s \nabla f(x_k), \label{subeq:extended-NAG-C-y} \\
            z_{k+1} &= z_k - \widetilde \alpha_k s \nabla f(x_k),\label{subeq:extended-NAG-C-z} \\
            x_{k+1} &= \frac{1}{\alpha_{k+1}}z_{k+1} + \left(1-\frac{1}{\alpha_{k+1}}\right) y_{k+1}, \label{subeq:extended-NAG-C-x}
        \end{align}
    \end{subequations}
starting from $x_0=z_0$, where $\widetilde \alpha_k=\beta_k\alpha_k + (\gamma_k-\beta_k)\alpha_{k+1}$.
See Supplement \ref{sec:tech for section c} for a proof.
By definition,  given $\{\beta_k\}$ and $\{\alpha_k >0\}$, there is a one-to-one correspondence between $\{ \gamma_k\}$ and $\{ \widetilde \alpha_k\}$:
$ \gamma_k = (\widetilde \alpha_k - \beta_k (\alpha_k - \alpha_{k+1}) )/ \alpha_{k+1}$.
Hence algorithm (\ref{eq:extended-NAG-C}) or (\ref{eq:extended-NAG-C-xyz}) can be considered to be directly parameterized by the three sequences $\{\alpha_k\}$, $\{\beta_k\}$ and $\{\widetilde\alpha_k\}$.

Our main result gives sufficient conditions on when algorithm (\ref{eq:extended-NAG-C})  or (\ref{eq:extended-NAG-C-xyz}) achieves accelerated convergence in both the objective gap and gradient norm similarly as NAG-C.
We discuss interpretations of the conditions from the perspective of HAG in Section \ref{sec:HAG-interpretation}.

\begin{thm}\label{thm:converge-c}
     Let $f\in\mathcal{F}^1_L$. There exist constants $C_0>0$ and $K\ge 1$, depending only on $\{\alpha_k\}$, $\{\beta_k\}$, and $\{\gamma_k\}$, such that
     for $0<s \le C_0 /L$ and $k\ge K$, the iterates $\{x_k\}$ from (\ref{eq:extended-NAG-C})  or (\ref{eq:extended-NAG-C-xyz})
     satisfy the optimal bound (\ref{eq:optimal bound-c}) for the objective gap and the inverse cubic rate (\ref{eq:inverse cubic rate})
     for the squared gradient norm, provided that the following conditions jointly hold.\vspace{-.05in}
    \begin{itemize} \addtolength{\itemsep}{-.05in}
        \item[(i)] $\lim_k \beta_k=\beta$ and $\lim_k \gamma_k=\gamma$ with $\beta>\gamma/2>0$;
        \item[(ii)] $\{\alpha_k>0\}$ satisfies that $\alpha_k=\Omega(k)$, $\lim_k \alpha_{k+1}/\alpha_k=1$, $\alpha_{k+1}(\alpha_{k+1}-1)\leq \alpha^2_k$;
        \item[(iii)] $\{\widetilde \alpha_k/\alpha_k\}$ is monotone (either non-increasing or non-decreasing) in $k$.
    \end{itemize} \vspace{-.05in}
\end{thm}

For Theorem~\ref{thm:converge-c},
the existence of the limits, $\lim_k \beta_k=\beta$, $\lim_k \gamma_k=\gamma$, and $\lim_k \alpha_{k+1}/\alpha_k=1$,
are introduced mainly to simplify the sufficient conditions and may be relaxed even with the same Lyapunov function in our proofs.
The monotonicity condition on $\{\widetilde \alpha_k/\alpha_k\}$ may also be relaxed, but an alternative Lyapunov function may be required.
In the setting of $\beta_k\equiv\beta$ and $\gamma_k\equiv \gamma$,
we have $\widetilde \alpha_k/\alpha_k = \beta + (\gamma-\beta)\alpha_{k+1}/\alpha_k$. Then the monotonicity condition is equivalent to
either requiring $\beta=\gamma$ without any additional constraint on $\{\alpha_k\}$ or, if $\beta\not=\gamma$,
requiring that $\{\alpha_{k+1}/\alpha_k\}$ is monotone.
The latter condition is satisfied by both the simple choice $\alpha_k = (k+r)/r $
and the iterative choice $\alpha_{k+1}=(1+\sqrt{1+4\alpha_{k}^2})/2$ with $\alpha_0=1$ in FISTA.

We also point out that the conditions in Theorem \ref{thm:converge-c} are general enough to allow that a limiting ODE may not exist for
algorithm (\ref{eq:extended-NAG-C}).
As described in Section \ref{sec:ODE-interpretation}, by taking $\Delta t = \sqrt{s} \to 0$,
NAG-C (\ref{eq:NAG-C-single}) admits a limiting ODE (\ref{eq:low-res-ODE-NAG-C}), where the coefficient $3/t$ results from the fact that
 $\sigma_{k+1}= k/(k+3)=1-3/k+O(1/k^2)$.
Similarly, the algorithm defined by substituting $\sigma_{k+1}= k/(k+r+1)$ for $k/(k+3)$ in (\ref{eq:NAG-C-single}) admits a limiting ODE (\ref{eq:low-res-ODE-C}).
For the existence of a limiting ODE for (\ref{eq:extended-NAG-C}), a necessary condition is that
$\lim_k k(1-\sigma_{k+1})$ exists. However, we provide an example where the conditions in  Theorem \ref{thm:converge-c} are satisfied,
but $\lim_k k(1-\sigma_{k+1})$ does not exist.

\begin{lem}\label{lem:counterexample for alphak}
    Consider the sequence $\{\alpha_k\}$ defined by alternating two rules:
    \begin{equation*}
        \alpha_k = \begin{cases}
        (k+r)/{r}, & \text{if } k \text{ is even},\\
        (1+\sqrt{1+4\alpha^2_{k-1}})/{2}. & \text{if } k \text{ is odd},
        \end{cases}
    \end{equation*}
    starting from $\alpha_0=1$. Then for any $r\geq 2$, condition (ii) in Theorem \ref{thm:converge-c} holds.
    But $\lim_k k(1-\sigma_{k+1})$ does not exist if $r>2$: along $\{k'\}=\{2k\}$ and $\{k''\}=\{2k+1\}$, we have
    \begin{equation*}
        \lim_{k'\to \infty} k'(1-\sigma_{k'+1}) = 3r/2, \quad \lim_{k''\to \infty} k''(1-\sigma_{k''+1}) = 2 + r/2.
    \end{equation*}
\end{lem}

Lemma~\ref{lem:counterexample for alphak} shows that accelerated convergence can be achieved by algorithm (\ref{eq:extended-NAG-C}), independently of
whether a limiting ODE exists. Hence convergence of discrete algorithms may not always be explained from the ODE perspective,
which is consistent with our findings in Section \ref{sec:ODE-interpretation}.

\section{ODE connection and comparison}\label{sec:ODE-interpretation}

As the stepsize $s$ vanishes to $0$ in a discrete algorithm, the limiting ODE (if exists) can be studied to understand
the behavior of its discrete counterpart \citep{su2016differential, shi2021understanding}.
We study convergence of the ODEs derived from the classes of algorithms in Sections~\ref{sec:acc-sc}--\ref{sec:acc-c},
and compare the conditions for when acceleration is achieved by the discrete algorithms and ODEs.
The comparison may help to understand the scope and mechanism of acceleration.

We briefly review how the limiting ODEs from NAG-SC (\ref{eq:NAG-SC-single}) and NAG-C (\ref{eq:NAG-C-single}) exhibit interesting differences from that of the vanilla gradient descent (\ref{eq:GD}), as discussed in \cite{su2016differential}.
On one hand, by taking $\Delta t= s$ and $x_k=X(t_k)=X(ks)$ for a continuous-time trajectory $X_t=X(t)$, the limit of gradient descent (\ref{eq:GD}) as $s\to 0$ is the gradient flow
\begin{equation}\label{eq:gradient flow}
    \dot{X}_t + \nabla f(X_t) = 0,
\end{equation}
with $X(0)=x_0$. It is known that when $f\in \mathcal{F}^1 $, $f(X_t)-f^*\leq \|x_0-x^*\|^2/(2t)$ and when $f\in \mathcal{S}^1_{\mu} $, $f(X_t)-f^*\leq \me^{-2\mu t}(f(x_0)-f^*)$. By relating $X_t$ to $x_k$ (i.e., taking $t=ks$ for small $s$), the former $O(1/t)$ rate translates into $O({1}/({ks}))$ which is exactly the discrete rate of gradient descent in the general convex setting. The second $O(\me^{-2\mu t})$ rate resembles $O((1-2\mu s)^k)$ for $s\approx 0$, which is similar to the discrete rate $O((1-\frac{2\mu s}{1+\mu/L})^k)$ for gradient descent in the strongly convex setting.

On the other hand, by taking $\Delta t=\sqrt s$ and $x_k=X(t_k)=X(k\sqrt s)$ for a continuous-time trajectory $X_t=X(t)$, the limit of NAG-C (\ref{eq:NAG-C-single}) as $s\to 0$ is
\begin{equation}\label{eq:low-res-ODE-NAG-C}
    \ddot{X}_t + \frac{3}{t}\dot{X}_t + \nabla f(X_t)=0,
\end{equation}
with initial conditions $X(0)=x_0$ and $\dot X(0)=0$.
The solution $X_t$ satisfies that for $f\in \mathcal{F}^1 $, $f(X_t)-f^*\leq 2\|x_0-x^*\|^2/t^2$. The limit of NAG-SC (\ref{eq:NAG-SC-single}) is
\begin{equation}\label{eq:low-res-ODE-NAG-SC}
    \ddot X_t + 2\sqrt{\mu} \dot X_t + \nabla f(X_t) = 0,
\end{equation}
with initial conditions $X(0)=x_0$ and $\dot X(0)=0$. The solution $X_t$ satisfies that when $f\in \mathcal{S}^1_{\mu} $, $f(X_t)-f^*\leq 2\me^{-\sqrt{\mu}t} (f(x_0)-f^*)$. By relating $X_t$ to $x_k$ (i.e., taking $t=k\sqrt{s}$), the former $O(1/t^2)$ rate translates into $O({1}/(sk^2))$, which is the discrete rate for NAG-C. The latter $O(\me^{-\sqrt{\mu}t})$ rate matches the discrete rate $O((1-\sqrt{\mu s})^k)$ for NAG-SC. In these cases, the convergence rate of the continuous-time trajectory $X_t=X(t)$ matches that of the discrete iterates $\{x_k\}$.

The overall findings of our study can be summarized as follows.
Although the convergence properties of discrete algorithms and those of limiting ODEs (low- or high-resolution) are related to each other,
there are currently notable gaps between the corresponding algorithms and ODEs
with respect to the configuration parameters studied,
except for the effect of the momentum term in the minimization of general convex functions \citep{su2016differential}.\vspace{-.05in}
\begin{itemize} \addtolength{\itemsep}{-.05in}

\item For minimization of strongly convex functions,
ODEs can be shown to converge under much weaker conditions
than the associated algorithms in terms of the momentum parameter.
%while comparable convergence rates are obtained.
In addition, the ODE convergence bounds do not directly inform the range of feasible stepsizes,\footnote{
This observation is also valid in the general convex setting. But in that case,
there does not seem to be the issue that allowing the stepsize in different ranges
leads to qualitatively different convergence rates, such as the linear or square-root dependency on $L/\mu$ in the strongly convex setting.
See the discussion of Proposition~\ref{pro:low-res-ODE-SC}.}
which are crucial in determining whether Nesterov's acceleration is achieved.

\item For minimization of either strongly convex or general convex functions,
the gradient correction parameter vanishes in the (low-resolution)  ODEs.
The high-resolution ODEs, although designed to reflect the gradient correction \citep{shi2021understanding},
can still converge at a similar rate even when the gradient correction parameter is nullified.
\end{itemize} \vspace{-.05in}
Our work focuses on comparison of the convergence properties between related algorithms and ODEs, and
does not address how algorithms can be studied by leveraging ideas (not necessarily results) from
studying ODEs, as shown in \cite{su2016differential} and \cite{shi2019symplectic, shi2021understanding} among others.
Our findings suggest that further work is needed to resolve the current gaps.

\subsection{Strongly convex setting} \label{sec:ODE-sc}

For strongly convex $f$, we compare convergence of algorithm (\ref{eq:extended-NAG-SC}) and related ODEs, which are derived from
a more general class of algorithms in the single-variable form (\ref{eq:extended-NAG-SC-single-2}),
with unspecified remainder terms $R_1$, $R_2$ and $R_3$ depending on $q$.
As discussed after Corollary~\ref{cor:symmetrized conditions for extended-NAG-SC-single}, algorithm (\ref{eq:extended-NAG-SC-single-2})
may be put into the three-variable form (\ref{eq:extended-NAG-SC}) only when $c_1^2 \ge 4c_0$.

By taking $\Delta t=\sqrt s$ and $x_k=X(t_k)=X(k\sqrt s)$, the limiting ODE of (\ref{eq:extended-NAG-SC-single-2}) is
\begin{equation}\label{eq:low-res-ODE-SC}
    \ddot X_t + c_1\sqrt{\mu} \dot X_t + c_0 \nabla f(X_t)=0,
\end{equation}
with initial conditions $X(0)=x_0$ and $\dot X(0)=0$. This can be viewed as a Newtonian equation of motion in a viscous medium in the potential field $c_0 f$. The damping coefficient is $c_1\sqrt{\mu}$, resulting from the momentum coefficient $1-c_1\sqrt{q}$ in (\ref{eq:extended-NAG-SC-single-2}).
For convenience, the three settings $c_1^2 > 4c_0$, $c_1^2 = 4 c_0$, or $c_1^2 < 4 c_0$ are referred to as
over-damping, critical-damping, or under-damping respectively.
However, the parameter $c_2$ associated with the gradient correction in (\ref{eq:extended-NAG-SC-single-2}) vanishes in the ODE (\ref{eq:low-res-ODE-SC}).
Nevertheless, for any $\mu$-strongly convex $f$ and any $c_0,c_1 >0$, we show that $f(X_t)$ converges to
 $f^*$ exponentially fast with a decaying rate proportional to $\sqrt{\mu}$. The condition, $c_0, c_1 > 0$, is much
weaker than previously realized, for example, $c_0 > 0$ and $c_0^2 < c_1 \le 4c_0^2$ in \cite{wilson2021lyapunov}, which covers ODE (\ref{eq:low-res-ODE-NAG-SC}) for NAG-SC with $(c_0, c_1) = (1, 2)$.

\begin{pro}\label{pro:low-res-ODE-SC}
    Let $f\in\mathcal{S}^1_\mu$. Then for any $c_0, c_1>0$, the solution $X_t$ to (\ref{eq:low-res-ODE-SC}) satisfies that $f(X_t)-f^* = O(\me^{-C\sqrt{\mu}t}(f(x_0)-f^*))$ for a constant $C=\frac{c_1-\sqrt{(c_1^2-4c_0) \vee 0}}{2}>0$. Otherwise (either $c_0\le 0$ or $c_1 \le 0$),
    $X_t$ may fail to converge to $x^*$ as $t\to\infty$.
\end{pro}

The dependency on $\mu$ in the bound above is improved to $\sqrt{\mu}$, compared with the $O(\me^{-2\mu t})$ bound for gradient flow (\ref{eq:gradient flow}). A direct translation by $t_k=k\sqrt{s}$ suggests that a discrete bound of $O((1-C\sqrt{\mu s})^k)$ may be expected for suitable discretizations of (\ref{eq:low-res-ODE-SC}), as opposed to $O((1-C\mu s)^k)$ for gradient descent (\ref{eq:GD}).
However, there are notable gaps between such results suggested from ODEs
and our results as well as existing results for discrete algorithms.

First, the bound in Proposition \ref{pro:low-res-ODE-SC} only requires $c_0,c_1>0$,
which is much weaker than the conditions in our convergence bounds in Section \ref{sec:acc-sc}.
Corollary~\ref{cor:symmetrized conditions for extended-NAG-SC-single}, deduced from Theorem \ref{thm:converge-sc-0} for (\ref{eq:extended-NAG-SC}),
establishes a convergence bound $(1-C\sqrt{\mu s})^k$ for algorithm (\ref{eq:extended-NAG-SC-single-2})
only in the case of $c_1^2 > 4c_0$ (over-damping) and additional conditions on the relationship of $c_0$ and $c_2$.
As mentioned after Corollary~\ref{cor:symmetrized conditions for extended-NAG-SC-single},
convergence of algorithm (\ref{eq:extended-NAG-SC-single-2}) can also be deduced from Theorems \ref{thm:converge-sc-const} and \ref{thm:converge-sc-1}
in case of $c_1^2=4c_0$ (critical-damping) and additional information about $R_1$, $R_2$ and $R_3$.
However, our results or, to our knowledge, existing results do not address
convergence of (\ref{eq:extended-NAG-SC-single-2}) in the case of $c_1^2 < 4c_0$ (under-damping).
As suggested by Proposition \ref{pro:low-res-ODE-SC}, in that case,
the bound $O((1-C\sqrt{\mu s})^k)$ may potentially be achievable by some discretization of (\ref{eq:low-res-ODE-SC}).
But this remains to be resolved, as the discretization would need to be outside the class of algorithms (\ref{eq:extended-NAG-SC}) studied.

Second, even if the form of the bound $O((1-C\sqrt{\mu s})^k)$ may be valid,
Proposition \ref{pro:low-res-ODE-SC} does not directly inform the range of feasible stepsize $s$,
which is needed to determine whether Nesterov's acceleration is achieved in the strongly convex setting.
%given such a bound, whether a discrete algorithm achieves Nesterov's accelerated bound, $O((1-C \sqrt{\frac{\mu}{L}})^k)$,
%depends on how large the stepsize $s$ is allowed, in terms of the smoothness constant $L$.
(Nesterov’s acceleration in Corollary \ref{cor:symmetrized conditions for extended-NAG-SC-single} applies to only $f\in\mathcal{S}^1_{\mu,L}$,
but $f\in\mathcal{S}^1_\mu$ is allowed in Proposition \ref{pro:low-res-ODE-SC}.)
As in the discussion of Theorem~\ref{thm:converge-sc-const}, from the bound $O((1-C\sqrt{\mu s})^k)$,
allowing  $s \asymp 1/L$ leads to Nesterov's acceleration under some conditions on the algorithm parameters,
whereas taking $s \asymp \mu/L^2$ results in a non-accelerated bound under other conditions.
Finding the range of feasible stepsize $s$ (involving $L$) and the associated conditions is essential in studying convergence of discrete algorithms.
In contrast, the range of $s$ bears no effect in the convergence of ODEs, which are defined as $s\to 0$.
(This also underlies why Proposition \ref{pro:low-res-ODE-SC} applies to $f\in\mathcal{S}^1_\mu$, regardless of $L$.)
Hence there is a limitation in how convergence properties of discrete algorithms can be transferred from those of limiting ODEs.

To complement the preceding discussion, we also study the high-resolution ODEs for
NAG-SC (\ref{eq:NAG-SC-single}) and the heavy-ball method (\ref{eq:HB}),
which are proposed to reflect the gradient correction or lack of in the two methods \citep{shi2021understanding}.
The high-resolution ODEs are derived by retaining $O(\sqrt{s})$ terms
which would otherwise vanish in the limit of $s \to 0$.
By taking $\Delta t = \sqrt{s}$ and $x_k=X(t_k)=X(k\sqrt{s})$, the low-resolution ODEs for NAG-SC (\ref{eq:NAG-SC-single}) and the heavy-ball method (\ref{eq:HB}) are the same equation in (\ref{eq:low-res-ODE-NAG-SC}). The high-resolution ODE for NAG-SC is
\begin{equation}\label{eq:high-res-ODE-NAG-SC}
    \ddot X_t + 2\sqrt{\mu}\dot X_t  + \sqrt{s}\nabla^2 f(X_t)\dot X_t + (1+\sqrt{\mu s})\nabla f(X_t) =0,
\end{equation}
with initial conditions $X(0)=x_0$ and $\dot X(0)=-\frac{2\sqrt{s}}{1+\sqrt{\mu s}}\nabla f(x_0)$. The high-resolution ODE for the heavy-ball method is
\begin{equation}\label{eq:high-res-ODE-HB}
    \ddot X_t + 2\sqrt{\mu}\dot X_t  + (1+\sqrt{\mu s})\nabla f(X_t) =0,
\end{equation}
with the same initial conditions.
For NAG-SC, the gradient correction results in a Hessian term, $\sqrt{s}\nabla^2 f(X_t)\dot X_t$, in (\ref{eq:high-res-ODE-NAG-SC}).
For the heavy-ball method, there is no gradient correction and hence no Hessian term in (\ref{eq:high-res-ODE-HB}).
However,  the convergence rates of these two high-resolution ODEs are the same.
The following result is qualitatively similar to Theorems 1 and 2 in \cite{shi2021understanding},
but involves a sharper objective gap due to a new Lyapunov analysis.

\begin{pro}\label{pro:high-res-ODE NAG and HB}
    Let $f\in\mathcal{S}^2_\mu $. Then the solutions $X_t$ of both (\ref{eq:high-res-ODE-NAG-SC}) and (\ref{eq:high-res-ODE-HB}) satisfy that $f(X_t)-f^* \leq V_0 \me^{-\sqrt{\mu}t}$, where the constant $V_0=(1+\sqrt{\mu s})(f(x_0)-f^*)+\frac{1}{2}\|\sqrt{\mu}(x_0-x^*)-\sqrt{s}\frac{1-\sqrt{\mu s}}{1+\sqrt{\mu s}}\nabla f(x_0)\|^2$ for (\ref{eq:high-res-ODE-NAG-SC}) and $V_0=(1+\sqrt{\mu s})(f(x_0)-f^*)+\frac{1}{2}\|\sqrt{\mu}(x_0-x^*)-\frac{2\sqrt{s}}{1+\sqrt{\mu s}}\nabla f(x_0)\|^2$ for (\ref{eq:high-res-ODE-HB}).
\end{pro}

Unfortunately, although the stepsize and the gradient correction (or lack of) are explicitly incorporated,
the high-resolution ODEs do not capture the difference that NAG-SC achieves Nesterov's acceleration, but heavy-ball does not.
A possible explanation can also be seen from the earlier discussion.
Nesterov's acceleration is directly affected by the range of feasible stepsize $s$, involving $L$, for a discrete algorithm.
But convergence of high-resolution ODEs remains unrelated to the range of stepsize $s$
(hence Proposition \ref{pro:high-res-ODE NAG and HB} applies to $f\in\mathcal{S}^2_\mu$, regardless of $L$).

\subsection{General convex setting} \label{sec:ODE-c}

For general convex $f$, we compare convergence of algorithm (\ref{eq:extended-NAG-C-single}) and related ODEs
in the setting where $\gamma_k =\gamma$, $\beta_k = \beta$, and
$\alpha_k=\frac{k+r}{r}$.
By rescaling the stepsize $s$ to $s/\gamma$, the parameters $\gamma$ and $\beta$ can be reset to $1$ and $\beta/\gamma$ respectively.
In this setting, (\ref{eq:extended-NAG-C-single}) reduces to the sub-class in \cite{shi2021understanding}, with two parameters $r$ and $\beta/\gamma$:
\begin{equation}\label{eq:extended-NAG-C-shi}
    x_{k+1} = x_{k} - s \nabla f(x_k) + \sigma_{k+1} (x_k-x_{k-1})  - \sigma_{k+1} \frac{\beta}{\gamma}  \cdot s(\nabla f(x_{k})-\nabla f(x_{k-1})),
\end{equation}
where $\sigma_{k+1} = \frac{\alpha_k-1}{\alpha_{k+1}}=\frac{k}{k+r+1} $.
For a general sequence $\{\alpha_k\}$,
 algorithm (\ref{eq:extended-NAG-C-single}) may not admit a limiting ODE, as shown in Lemma \ref{lem:counterexample for alphak}.

By taking $\Delta t=\sqrt s$ and $x_k=X(t_k)=X(k\sqrt s)$, the limiting ODE for (\ref{eq:extended-NAG-C-shi}) is
\begin{equation}\label{eq:low-res-ODE-C}
    \ddot X_t + \frac{r+1}{t}\dot X_t + \nabla f(X_t)=0,
\end{equation}
with initial conditions $X(0)=x_0$ and $\dot X(0)=0$. This can be viewed as a Newtonian equation for a particle moving in the potential field $f$ with friction.
The damping coefficient is $\frac{r+1}{t}$, resulting from the momentum coefficient $\sigma_{k+1}$.
\cite{su2016differential} showed that the convergence of $X_t$ from (\ref{eq:low-res-ODE-C}) exhibits a phase transition at $r=2$, the choice in NAG-C (\ref{eq:NAG-C}). For $f\in\mathcal{F}^1$, if $r\geq 2$, then $f(X_t)-f^*=O({1}/{t^2})$. But if $r<2$, the rate $O({1}/{t^2})$ may fail as illustrated by counterexamples. The condition
for $O({1}/{t^2})$ convergence of ODE (\ref{eq:low-res-ODE-C}) is consistent with that for algorithm (\ref{eq:extended-NAG-C-shi})
to achieve $O(1/(sk^2))$ convergence in Theorem \ref{thm:converge-c}.\footnote{In case of $r\ge 2$, Nesterov's acceleration for algorithm (\ref{eq:extended-NAG-C-shi}) is applicable only to $f \in\mathcal{F}^1_L$,
whereas the $O({1}/{t^2})$ convergence of ODE (\ref{eq:low-res-ODE-C}) holds for $f\in\mathcal{F}^1$.
A similar distinction is mentioned in the discussion of Proposition \ref{pro:low-res-ODE-SC}.}
For $\alpha_k=\frac{k+r}{r}$, it can be directly verified that
the recursive condition $\alpha_{k+1}(\alpha_{k+1}-1)\leq \alpha_k^2$ holds if and only if $r\geq 2$.
Therefore, the condition on the momentum term for Nesterov's acceleration of (\ref{eq:extended-NAG-C-shi})
is well captured by the limiting ODE.\
Such a phenomenon is not observed in the strongly convex setting in Section \ref{sec:ODE-sc},
where regarding the momentum term, the condition for convergence of the limiting ODE (\ref{eq:low-res-ODE-SC}) is much weaker than
those currently established for algorithm (\ref{eq:extended-NAG-SC-single-2}) to achieve
a comparable convergence rate.

The ODE (\ref{eq:low-res-ODE-C}), however, does not reflect
the parameter $\beta/\gamma$ associated with the gradient correction in algorithm (\ref{eq:extended-NAG-C-shi}),
which is similarly observed in the low-resolution ODE (\ref{eq:low-res-ODE-SC}) for algorithm (\ref{eq:extended-NAG-SC-single-2}) in the strongly convex setting.
The gradient correction term, of order $O(\sqrt{s})$, can be incorporated in a high-resolution ODE \citep{shi2021understanding}.
By taking $t_0=\frac{(r+1)\sqrt{s}}{2}$ and $t_k=t_0+k\sqrt{s}$ for algebraic convenience, the high-resolution ODE for (\ref{eq:extended-NAG-C-single}) is
\begin{equation}\label{eq:high-res-ODE-C}
    \ddot X_t + \frac{r+1}{t}\dot X_t + \frac{\beta}{\gamma} \sqrt{s}\nabla^2 f(X_t)\dot X_t + \left(1+\frac{(r+1) \sqrt{s}}{2t}\right)\nabla f(X_t)=0,
\end{equation}
with initial conditions $X(t_0)=x_0$ and $\dot X(t_0)=-\sqrt{s}\nabla f(x_0)$.
The following result gives convergence properties of (\ref{eq:high-res-ODE-C}), with a new Lyapunov analysis in the case of $r=2$.
\cite{shi2021understanding} provided similar bounds only explicitly in the case of $r >2$.

\begin{pro}\label{pro:high-res-ODE-C}
Let $f\in\mathcal{F}^2$. If $r\geq 2$ and $\frac{\beta}{\gamma} >0$, then there exists a time point $t_1\geq t_0$,
with $t_1/\sqrt{s}$ depending only on $r$ and $\beta/\gamma$, such that the solution to (\ref{eq:high-res-ODE-C}) satisfies the bounds for $t> t_1$:
\begin{equation*}
    f(X_t)-f^* = O\left(\frac{V_{t_1}}{t^2}\right), \quad \inf_{t_1\leq u \leq t} \|\nabla f(X_u)\|^2 = O\left(\frac{\gamma V_{t_1}}{\beta\sqrt{s}(t^3-t_1^3)}\right).
\end{equation*}
In addition, if $r\geq 2$ and $\beta=0$, then the first bound still holds. Here $V_{t_1}$ is the value that the continuous Lyapunov function takes at $t_1$: $V_{t}= (t+C\sqrt{s})(t+ (\frac{r+1}{2}-\frac{\beta}{\gamma})\sqrt{s})(f(X_t)-f^*) + \frac{t+C\sqrt{s}}{t}\frac{1}{2}\|r(X_t-x^*)+t(\dot X_t + \frac{\beta\sqrt{s}}{\gamma}\nabla f(X_t))\|^2$, where $C=0$ if $r>2$ and $C=\frac{\beta}{\gamma}$ if $r=2$.
\end{pro}

If we replace $t$ by $t_k=t_0+k\sqrt{s}\sim k\sqrt{s}$, then the first bound above translates into
the accelerated bound $O({1}/(sk^2))$ on the objective gap as in (\ref{eq:optimal bound-c}) and
the second becomes $O({1}/{(s^2k^3)})$ on the squared gradient norm as in (\ref{eq:inverse cubic rate}).
These two bounds above are exactly the continuous analogs of the associated bounds for algorithm (\ref{eq:extended-NAG-C-shi}). However, compared with the condition
$\beta/\gamma > 1/2$ in Theorem \ref{thm:converge-c}, the condition in the continuous ODE setting is much weaker: $\beta/\gamma$ is allowed to be arbitrarily small, and even $0$ (no Hessian term) if only the $O(1/t^2)$ rate for the objective gap is desired.
The latter case corresponds to a discrete algorithm without the gradient correction, similar to the heavy-ball method (\ref{eq:HB}).
Therefore, even when a Hessian term is included as a continuous counterpart of the gradient correction, the convergence properties of the high-resolution ODE are currently not fully matched by those of algorithm (\ref{eq:extended-NAG-C-shi}).
This remains to be further studied as an open question, as also mentioned in \cite{shi2021understanding}, Section 5.1.

\section{Hamiltonian assisted interpretation} \label{sec:HAG-interpretation}

As observed in Section \ref{sec:ODE-interpretation},
the low-resolution and high-resolution ODEs do not fully capture the convergence properties of the algorithms studied.
Alternatively, we directly formulate a broad class of discrete algorithms, HAG, based on a Hamiltonian function with a position variable and a momentum variable,
and demonstrate that the conditions from our convergence results in Sections
 \ref{sec:acc-sc}--\ref{sec:acc-c} can be interpreted through HAG in a meaningful and unified manner
in both the strongly convex and general convex settings.
This development is motivated by a related formulation of Hamiltonian assisted Metropolis sampling (HAMS) \citep{song2021hamiltonian,song2022irreversible}.

\subsection{HAG: Hamiltonian assisted gradient method}

For $f\in\mathcal{F}^1_L$, consider the following unconstrained minimization problem:
\begin{equation*}
    \min_{x,u \in\bbR^n} H(x,u)= f(x) + \frac{1 }{2 }\|u\|^2,
\end{equation*}
where $H(x,u)$ can be interpreted as a Hamiltonian function (or total energy),
with $x$  a position variable and $u$ a momentum variable.
On one hand, the above problem is apparently equivalent to the original problem (\ref{eq:main_optimization}).
If $x$ and $u$ are updated separately, then there is no possible improvement, compared with solving  (\ref{eq:main_optimization}) directly.
On the other hand, we show how $x$ and $u$ can be updated in a coupled manner to derive a rich class of first-order algorithms, which include
representative algorithms from the classes (\ref{eq:extended-NAG-SC}) and (\ref{eq:extended-NAG-C-single}) studied in Sections \ref{sec:acc-sc}--\ref{sec:acc-c}.

\textbf{Gradient descent.} Given initial points $(x_0,u_0)$, consider minimizing (or decreasing) $H(x,u)$ subject to a linear constraint $x- \delta u = x_0 - \delta u_0$ for some $\delta >0$.
By substituting $u=\frac{x-x_0+ \delta u_0}{ \delta}$, the problem  with respect to $x$ becomes minimizing (or decreasing)
$\tilde f(x; x_0,u_0) = f(x) + \frac{1}{2 \delta^2}\|x-x_0+ \delta u_0\|^2$.
Consider updating $x_0$ by gradient descent, i.e., minimizing a quadratic surrogate function at $x_0$
(with surrogate smoothness parameter being $1$ for $f$) as follows:
\begin{equation*}
    \begin{split}
        x_1 &= \argmin_{x}\; \langle \nabla f(x_0),x-x_0\rangle + \frac{1}{2}\|x-x_0\|^2 + \frac{1}{2 \delta^2}\|x-x_0+ \delta u_0\|^2 ,
    \end{split}
\end{equation*}
and then updating $u_0$ by solving $x_1- \delta u_1 = x_0- \delta u_0$. This results in the update:
\begin{equation*}
    \begin{split}
        x_1 &= x_0 - \frac{ \delta^2}{1+\delta^2}\nabla f(x_0) - \frac{ \delta}{1+\delta^2}u_0, \\
        u_1 &= u_0 - \frac{ \delta}{1+\delta^2}\nabla f(x_0) - \frac{1}{1+\delta^2} u_0.
    \end{split}
\end{equation*}
The function $\tilde f(x; x_0,u_0)$ has a smoothness parameter $L_{\tilde f} = L + \frac{1}{\delta^2}$, and
the GD stepsize above is $s = \frac{1}{1 + \frac{1}{\delta^2}} = \frac{ \delta^2}{1+\delta^2}$.
If $0< \delta^2 \leq \frac{1}{L}$, then
$ L_{\tilde f} \le \frac{2}{\delta^2} < \frac{2}{s}$ and $s < \frac{2}{L_{\tilde f}}$,
so that the above GD update is well-behaved and $H$ is non-increasing (see the footnote \ref{fn:GD with large s}).

\textbf{Momentum negation and extrapolation.}
The above update can be repeated on $(x_1,u_1)$, but the iterates will be stuck in the same hyperplane $x- \delta u = x_0- \delta u_0$.
To resolve this issue, we negate $u$ at the beginning of each iteration. In other words, given $(x_0,u_0)$ we first negate $u_0$ to $-u_0$, construct a hyperplane passing $(x_0,-u_0)$ and then implement the above update. The negation of $u$ does not change the value of $H$, but keeps the hyperplane changing from iteration to iteration. In addition, we introduce an extrapolation step before the negation of $u_0$ for $\rho \ge 0$. This leads to the following update
given the initial points $(x_0,u_0)$:
\begin{equation*}
    \begin{split}
        x_1 &= x_0 - \frac{(1+\rho) \delta^2}{1+\delta^2}\nabla f(x_0) + \frac{ (1+\rho) \delta}{1+\delta^2}u_0, \\
        u_1 &= -u_0 - \frac{(1+\rho) \delta}{1+\delta^2}\nabla f(x_0) + \frac{ 1+\rho }{1+\delta^2}u_0.
    \end{split}
\end{equation*}
By setting $a=\frac{(1+\rho)\delta^2}{1+\delta^2}$ and $b=\frac{ 1+\rho }{1+\delta^2}$,
the above can be put in a simple form as
\begin{equation*}
    \begin{split}
        x_1 &= x_0 - a\nabla f(x_0) + \sqrt{ab}u_0, \\
        u_1 &= -u_0 - \sqrt{ab}\nabla f(x_0) + bu_0.
    \end{split}
\end{equation*}
The parameters $(\delta,\rho)$ can be determined from $(a,b)$ as
$\delta = \sqrt{ a/b }$ and $\rho= a+b-1 $.

\textbf{Gradient correction.}
The above updates of $x_1$ and $u_1$ are parallel, so that $u_1$ uses only the gradient information at $x_0$, but not at the newly updated $x_1$.
To exploit the gradient information at $x_1$, we incorporate a gradient correction into the update:
\begin{equation*}
    \begin{split}
        x_1 &= x_0 - a\nabla f(x_0) + \sqrt{ab}u_0, \\
        u_1 &= -u_0 - \sqrt{ab}\nabla f(x_0) + bu_0 - \phi(\nabla f(x_1)-\nabla f(x_0)),
    \end{split}
\end{equation*}
for some scalar $\phi$. Lastly, we allow $a$, $b$ and $\phi$ to vary from iteration to iteration and consider the following class of algorithms,
called Hamiltonian assisted gradient method (HAG):
\begin{equation}\label{eq:HAG}
    \begin{split}
        x_{k+1} &= x_k - a_k\nabla f(x_k) + \sqrt{a_k b_k}u_k, \\
        u_{k+1} &= -u_k - \sqrt{a_k b_k}\nabla f(x_k) + b_k u_k - \phi_k(\nabla f(x_{k+1})-\nabla f(x_k)),
    \end{split}
\end{equation}
starting from $(x_0,u_0)$. Notably, the overall contribution from gradients remains unchanged in the update for $u_{k+1}$, as the new gradient $\nabla f(x_{k+1})$ and old gradient $\nabla f(x_k)$ are introduced in opposite signs. The two gradients are re-weighted in the momentum update:
\begin{equation}\label{eq:update of momentum-reweighted gradients}
    u_{k+1} = (b_k-1) u_k - \underbrace{\big( (\sqrt{a_k b_k}-\phi_k)\nabla f(x_k) + \phi_k\nabla f(x_{k+1}) \big)}_{\text{re-weighted gradient}}.
\end{equation}
The larger $\phi_k$ is, the larger weight is assigned to the new gradient $\nabla f(x_{k+1})$.
%When $\phi_k = \sqrt{a_k b_k}$, the effect of $\phi$-term is just like replacing the old gradient $\nabla f(x_k)$ entirely by the new gradient $\nabla f(x_{k+1})$.

To facilitate comparison with algorithms in Sections \ref{sec:acc-sc}--\ref{sec:acc-c},
the HAG algorithm (\ref{eq:HAG}) can be put into a single-variable form involving $\{x_k\}$ only (see Supplement Section \ref{sec:tech for section HAG}):
\begin{equation}\label{eq:HAG-single}
    \begin{split}
        x_{k+1} &= x_k - \left(a_{k-1}\sqrt{\frac{a_k b_k}{a_{k-1}b_{k-1}}}+a_k\right)\nabla f(x_k) +  (b_{k-1}-1)\sqrt{\frac{a_k b_k}{a_{k-1}b_{k-1}}}(x_k-x_{k-1})  \\
        &\quad - \left(\phi_{k-1}\sqrt{a_k b_k}-a_{k-1}\sqrt{\frac{a_k b_k}{a_{k-1}b_{k-1}}}\right)(\nabla f(x_k)-\nabla f(x_{k-1})),
    \end{split}
\end{equation}
starting from $x_0$ and $x_1 = x_0-a_0\nabla f(x_0) + \sqrt{a_0 b_0}u_0$. The parameter $\{\phi_k\}$ only appears in the gradient correction in (\ref{eq:HAG-single}),
and in the case of $a_kb_k$ being constant in $k$,
the parameters $\{a_k\}$ and $\{b_k\}$ fully determine the gradient descent and the momentum terms respectively.
For convenience, the term of $x_k-x_{k-1}$ is still referred to as the momentum term in single-variable forms,
although a momentum variable $u_k$ is explicitly introduced in HAG.

\subsection{Interpretation from HAG}\label{subsec:interpretation-HAG}

HAG (\ref{eq:HAG}) or (\ref{eq:HAG-single}) are general enough to represent various algorithms studied in Sections \ref{sec:acc-sc}--\ref{sec:acc-c},
by choosing the parameters $\{a_k\}$, $\{b_k\}$ and $\{\phi_k\}$ accordingly.
We examine our convergence results from the HAG perspective,
and obtain meaningful interpretations of the conditions
for when Nesterov's acceleration is achieved in both strongly convex and general convex settings.\vspace{-.05in}
\begin{itemize} \addtolength{\itemsep}{-.05in}
\item The parameter $a_k$ acts as a re-scaled stepsize $s$. The parameter $b_k$, which controls the momentum
term in the single-variable form, has a leading constant 2,
and the gap $2-b_k$ is $\Omega(\sqrt{q})$ or $\Omega(\frac{1}{k})$ in the strongly or general convex setting respectively.
In terms of $(\delta_k,\rho_k)$ in the HAG derivation, this indicates that $\delta_k\sim \sqrt{\frac{a_k}{2}}$ and $1-\rho_k \sim 2-b_k$. Therefore, $\delta_k$ is of order $\sqrt{s}$, $\rho_k$ has a leading constant 1 (symmetric extrapolation),
and the gap  $1-\rho_k$ is $\Omega(\sqrt{q})$ or $\Omega(\frac{1}{k})$ in the strongly or general convex setting.

\item The parameter $\phi_k$ is of order $\sqrt{s}$ and controls the gradient correction such that in the re-weighted gradient for the momentum update
(\ref{eq:update of momentum-reweighted gradients}), the new gradient $\nabla f(x_{k+1})$ fully dominates
with its weight greater than the total weight, whereas the old gradient $\nabla f(x_k)$ has a negative weight.
The boundary case of a zero weight for $\nabla f(x_k)$ is also allowed in the strongly convex setting. Such heuristic interpretations are not feasible from the single-variable forms.
\end{itemize} \vspace{-.05in}
Currently, these interpretations are derived from the convergence results in Sections \ref{sec:acc-sc}--\ref{sec:acc-c}, where
our Lyapunov analyses are based on the three-variable forms (\ref{eq:extended-NAG-SC}) and (\ref{eq:extended-NAG-C-xyz}).
It remains an interesting question to develop a direct analysis of HAG and study further implications of HAG.

\textbf{Strongly convex setting.}
For $f\in\mathcal{S}^1_{\mu,L} $, we set
\begin{align*}
a_k\equiv a= (\frac{c_0}{2}+O(\sqrt{q}))s, \quad b_k\equiv b= 2-c_1\sqrt{q} + O(q), \quad \phi_k\equiv \phi= (c_2+O(\sqrt{q}))\sqrt{s},
\end{align*}
where $c_0, c_1, c_2 >0$. Then HAG (\ref{eq:HAG-single}) reduces to (\ref{eq:extended-NAG-SC-single-2}) exactly:
\begin{equation} \label{eq:extended-NAG-SC-single-2b}
    \begin{split}
        x_{k+1} &= x_k - 2a\nabla f(x_k) + (b-1)(x_k-x_{k-1}) - (\phi\sqrt{ab}-a)(\nabla f(x_k)-\nabla f(x_{k-1})) \\
        &= x_k - (c_0 + O(\sqrt{q}))s\nabla f(x_k) + (1-c_1\sqrt{q}+O(q))(x_k-x_{k-1}) \\
        &\quad - (c_2\sqrt{c_0}-\frac{c_0}{2} + O(\sqrt{q}))s(\nabla f(x_k)-\nabla f(x_{k-1})).
    \end{split}
\end{equation}
In the following, we interpret the conditions in Corollary \ref{cor:symmetrized conditions for extended-NAG-SC-single} through HAG.

The parameter $a=(\frac{c_0}{2}+O(\sqrt{q}))s$ plays the role of a re-scaled stepsize $s$. There is no condition on $c_0$ alone
but all conditions are stated in the ratios $c_1^2/c_0$ and $c_2^2/ c_0$.

The parameter $b= 2-c_1\sqrt{q} + O(q)$ controls the momentum term $(x_k-x_{k-1})$ in (\ref{eq:extended-NAG-SC-single-2b}).
The leading constant of $b$ is $2$, which ensures  the existence of a limiting ODE as $s\to 0$.
For the rate $(1-C\sqrt{\mu s})^k$, which gives either non-accelerated or accelerated convergence depending on
the range of $s$, the sufficient conditions in Corollary \ref{cor:symmetrized conditions for extended-NAG-SC-single} require $c_1^2> 4c_0$ (over-damping).\footnote{Note that $c_1^2=4c_0$ (critical-damping) may also be feasible, depending on higher-order information as in Theorems \ref{thm:converge-sc-const} and \ref{thm:converge-sc-1}. See the discussion after Corollary \ref{cor:symmetrized conditions for extended-NAG-SC-single}.}
With $q = \mu s$, this indicates a simple condition about $a$ and $b$ as $s\to 0$:
\begin{align*}
\frac{2-b} {\sqrt{2a}} = \frac{c_1\sqrt{q} }{\sqrt{c_0 s}} + O(\sqrt{s}) > 2 \sqrt{\mu} + O(\sqrt{s}),
\end{align*}
where $2 \sqrt{\mu}$ is exactly the coefficient of $\dot X_t$ in ODE (\ref{eq:low-res-ODE-NAG-SC}) for NAG-SC.
In terms of $(\delta,\rho)$ in the HAG derivation, the preceding discussion also gives
\begin{align}  \label{eq:HAG-delta-rho-SC}
\delta = \sqrt{a/b} = \frac{1}{2} \sqrt{c_0 s} + O(s),\quad
1 - \rho = 2-b-a > 2 \sqrt{\mu (c_0 s)} + O(s).
\end{align}
Hence $\delta$ is in the order of $\sqrt{s} $, similarly as $\Delta t$ in deriving the ODEs (\ref{eq:low-res-ODE-NAG-SC}) and (\ref{eq:low-res-ODE-SC}),
and $1- \rho = \Omega(\sqrt{\mu s})$ as $s\to 0$, i.e., $\rho$ is close to 1 (symmetric extrapolation) but with a gap being $\Omega( \sqrt{\mu s})$.

The parameter $\phi=(c_2+O(\sqrt{q}))\sqrt{s}$ controls the gradient correction term $(\nabla f(x_{k+1})-\nabla f(x_k))$ in (\ref{eq:extended-NAG-SC-single-2b}).
The re-wighted gradient in the momentum update (\ref{eq:update of momentum-reweighted gradients}) for $u_{k+1}$ becomes
\begin{align} \label{eq:reweighted-grad-SC}
      (\sqrt{c_0}-c_2+O(\sqrt{q}))\sqrt{s}\nabla f(x_k)+ (c_2+O(\sqrt{q}))\sqrt{s}\nabla f(x_{k+1})  .
\end{align}
For accelerated convergence in Corollary \ref{cor:symmetrized conditions for extended-NAG-SC-single},
the condition $c_2^2\geq c_0$ indicates that the new gradient $\nabla f(x_{k+1})$ fully dominates
with its weight, $c_2 \sqrt{s}$, no smaller than the total weight $\sqrt{c_0 s}$,
whereas the old gradient $\nabla f(x_k)$ has a zero or negative weight.
For non-accelerated (sub-optimal) convergence in Corollary \ref{cor:symmetrized conditions for extended-NAG-SC-single},
the condition $c_0/4\leq c_2^2 <c_0$ indicates that
both the new and old gradients
have positive weights, but $\nabla f(x_{k+1})$ contributes no less than $\nabla f(x_k)$.
We notice that a re-weighted gradient can also be identified in the single-variable form (\ref{eq:extended-NAG-SC-single-2}) or (\ref{eq:extended-NAG-SC-single-2b}):
\begin{align*}
& \quad (c_0 + O(\sqrt{q}))s\nabla f(x_k) + (c_2\sqrt{c_0}-\frac{c_0}{2} + O(\sqrt{q}))s(\nabla f(x_k)-\nabla f(x_{k-1}))\\
& =  (c_2\sqrt{c_0} + \frac{c_0}{2} + O(\sqrt{q}))s \nabla f(x_k) + (-c_2\sqrt{c_0} +\frac{c_0}{2} + O(\sqrt{q}))s \nabla f(x_{k-1}) .
\end{align*}
But there seems to be no meaningful interpretation. For example,
the condition $c_2^2\geq c_0$ indicates that the weight of $\nabla f(x_k)$ is no smaller than $\frac{3}{2} c_0$,
whereas the weight of $\nabla f(x_{k-1})$ is no larger than $-\frac{1}{2} c_0$,
although the total weight is $c_0 s$, independently of $c_2$.

\textbf{General convex setting.}
For $f\in\mathcal{F}^1_L $, we set
\begin{equation}\label{eq:config-HAG-C}
    b_k=1+\sigma_{k+2}, \quad a_k=\frac{c_0}{b_k}s, \quad \phi_k=c_2 \sqrt{s},
\end{equation}
where $c_0,c_2>0$ and, as before, $\sigma_{k+1}=\frac{\alpha_k-1}{\alpha_{k+1}}$ for some sequence $\{\alpha_k\}$.
Then HAG (\ref{eq:HAG-single}) reduces to
\begin{equation}\label{eq:HAG-C}
    \begin{split}
        x_{k+1} &= x_k - c_0 \left(\frac{1}{b_{k-1}}+\frac{1}{b_k}\right) s \nabla f(x_k) + \sigma_{k+1}(x_k-x_{k-1}) \\
        &\quad -\left(c_2 \sqrt{c_0}-\frac{c_0}{b_{k-1}}\right)s(\nabla f(x_{k})-\nabla f(x_{k-1})),
    \end{split}
\end{equation}
which falls in the class (\ref{eq:extended-NAG-C-single}) with $\gamma_k$ and $\beta_k$ varying in $k$.
We choose $a_k$ varying in $k$ such that $a_k b_k$ is constant in $k$, mainly to simplify the coefficients in (\ref{eq:HAG-single}).

Given $\alpha_k = \Omega(k)$ and $\frac{\alpha_{k+1}}{\alpha_k}=1+O(\frac{1}{k})$ (hence $\sigma_{k+1}=1+O(\frac{1}{k})$ and $b_k=2+O(\frac{1}{k})$),
the preceding HAG algorithm (\ref{eq:HAG-C}) is further simplified to
\begin{equation} \label{eq:extended-NAG-SC-single-2c}
    \begin{split}
        x_{k+1} &= x_k - \left(c_0 + O(\frac{1}{k})\right) s \nabla f(x_k) + \sigma_{k+1}(x_k-x_{k-1}) \\
        &\quad -\left(c_2 \sqrt{c_0}-\frac{c_0}{2} + O(\frac{1}{k}) \right)s(\nabla f(x_{k})-\nabla f(x_{k-1})),
    \end{split}
\end{equation}
which resembles algorithm (\ref{eq:extended-NAG-SC-single-2b}) in the strongly convex setting, in all the leading constants involved.
Algorithm (\ref{eq:extended-NAG-SC-single-2c}) can be put into (\ref{eq:extended-NAG-C-single})  with $\gamma_k\to \gamma= c_0$ and $\beta_k\to\beta = c_2 \sqrt{c_0}-\frac{c_0}{2}$. 
Note that (\ref{eq:HAG-C}) or (\ref{eq:extended-NAG-SC-single-2c}) 
is more complex than (\ref{eq:extended-NAG-C-shi}) with constant $\gamma_k$ and $\beta_k$.
In the following, we interpret conditions (i) and (ii) in Theorem \ref{thm:converge-c} through HAG.
The monotonicity condition (iii) is satisfied for the choice $\alpha_k = \frac{k+r}{r}$ discussed below (see Supplement \ref{sec:tech for section HAG} for a proof).

Similarly as in the strongly convex setting, the parameter $a_k = (c_0+O(\frac{1}{k}))s/2$ plays the role of a re-scaled stepsize $s$. Moreover, the parameter $b_k=1+\sigma_{k+2}$ controls the momentum term $(x_k-x_{k-1})$ in (\ref{eq:extended-NAG-SC-single-2c}), and has its leading constant being 2, which is necessary (but not sufficient by Lemma \ref{lem:counterexample for alphak})
for a limiting ODE to exist as $s\to 0$. For the choice $\alpha_k = \frac{k+r}{r}$
and $\sigma_{k+1} = \frac{k}{k+r+1} = 1-\frac{r+1}{k}+O(\frac{1}{k^2})$,
the recursive condition $\alpha_{k+1}(\alpha_{k+1}-1)\leq \alpha_k^2$ is equivalent to $r\geq 2$. By taking 
$\Delta t = \sqrt{c_0 s}$ and $t=k\Delta t=k\sqrt{c_0 s}$,
the above indicates a simple condition about $a_k$ and $b_k$ as $k\to \infty$ (or $s\to 0$):
\begin{align*}
\frac{2-b_k} {\sqrt{2a_k}} = \frac{r+1}{k\sqrt{c_0 s}} + O(\frac{1}{k^2 \sqrt{s}})=\frac{r+1}{t}+O(\frac{1}{k}) \ge \frac{ 3}{t} + O(\frac{1}{k}),
\end{align*}
where $\frac{3}{t}$ matches the coefficient of $\dot X_t$ in ODE (\ref{eq:low-res-ODE-NAG-C}) for NAG-C.
In terms of the time-dependent $(\delta_k,\rho_k)$ in HAG, the preceding discussion gives
\begin{align} \label{eq:HAG-delta-rho-C}
\delta_k = \sqrt{a_k/b_k} = \frac{1}{2} \sqrt{c_0 s} + O(s),\quad
1 - \rho_k = 2-b_k-a_k \ge \frac{3}{k} + O(\frac{1}{k^2}).
\end{align}
Hence $\delta_k$ is in the order of $\sqrt{s} $, similarly as $\Delta t$ in deriving the ODE (\ref{eq:low-res-ODE-NAG-C}),
and $1- \rho_k = \Omega (\frac{1}{k})$ as $k \to \infty$ (or $s\to 0$), i.e., $\rho_k$ is close to 1 (symmetric extrapolation) but with a gap being $\Omega(\frac{1}{k})$.
The conditions, (\ref{eq:HAG-delta-rho-SC}) and (\ref{eq:HAG-delta-rho-C}), are of similar forms,
with $\sqrt{\mu s}\, (=\sqrt{q})$ and $\frac{1}{k}$ replaced by each other.
Although the correspondence between $\sqrt{q}$ and $\frac{1}{k}$ can be seen in NAG-SC (\ref{eq:NAG-SC}) and NAG-C (\ref{eq:NAG-C}),
our work shows that accelerated convergence in both the strongly convex and general convex settings
can be extended to HAG algorithms
provided that the extrapolation gap $1-\rho$ is greater than $\sqrt{q}$ and $\frac{1}{k}$, up to constant factors,
in addition to other conditions discussed below.

The parameter $\phi_k=c_2\sqrt{s}$ plays a similar role of controlling the gradient correction term $(\nabla f(x_k)-\nabla f(x_{k-1}))$ in (\ref{eq:extended-NAG-SC-single-2c}), as does $\phi=(c_2+O(\sqrt{q}))\sqrt{s}$ in the strongly convex setting.
The re-weighted gradient in the momentum update (\ref{eq:update of momentum-reweighted gradients}) for $u_{k+1}$ is similar to (\ref{eq:reweighted-grad-SC}):
\begin{align*}
        (\sqrt{c_0}-c_2 )\sqrt{s}\nabla f(x_k) + c_2 \sqrt{s}\nabla f(x_{k+1}).
\end{align*}
With $ \gamma= c_0$ and $ \beta = c_2 \sqrt{c_0}-\frac{c_0}{2}$, the condition $\beta > \gamma/2$  in Theorem \ref{thm:converge-c}
reduces to $c_2^2>c_0$, which indicates that the new gradient $\nabla f(x_{k+1})$ fully dominates
with its weight, $c_2 \sqrt{s}$, greater than the total weight $\sqrt{c_0 s}$,
whereas the old gradient $\nabla f(x_k)$ has a negative weight.
This is the same as the condition, $c_2^2 \ge c_0$,
for accelerated convergence in the strongly convex setting except that the boundary case $c_2^2 = c_0$ is excluded here.
As discussed earlier, the single-variable form (\ref{eq:extended-NAG-SC-single-2c})
does not admit a meaningful interpretation for the re-weighted gradient.

\section{Outlines of Lyapunov analyses}\label{sec:outline}

We outline our Lyapunov analyses to prove the convergence results for the discrete algorithms in Sections \ref{sec:acc-sc}--\ref{sec:acc-c}.
See Supplement Section \ref{sec:tech for section ODE} for our Lyapunov analyses for the convergence of ODEs in Section \ref{sec:ODE-interpretation}.
Compared with existing ones, our Lyapunov functions are constructed to handle more general algorithms and ODEs
or to achieve simpler analysis and sometimes sharper results.
Before the outlines of our analyses, we summarize the comparison of Lyapunov analyses.

\textbf{Strongly convex setting.} We construct the discrete Lyapunov (\ref{eq:lyapunov-SC})
to establish the convergence of algorithm (\ref{eq:extended-NAG-SC}) including NAG-SC and TMM as special cases.
The auxiliary-energy term $\mu\|z_{k+1}-x^*\|^2/2$ in (\ref{eq:lyapunov-SC}) is also used in \cite{bansal2019potential} for NAG-SC and in
\cite{aspremont2021acceleration} for TMM.\ However,
the potential-energy term in our Lyapunov  is
in $\{x_k\}$ whereas the one in \cite{bansal2019potential} is in $\{y_k\}$.
The potential-energy term of the Lyapunov in \cite{aspremont2021acceleration} is in $\{x_k\}$ like ours,
but is not lower-bounded by $f(x_k)-f^*$ so that their analysis only establishes the convergence for $\{z_k\}$.
In addition, compared with the analysis of NAG-SC in \cite{shi2021understanding}, our Lyapunov function (\ref{eq:lyapunov-SC})
has fewer terms and our analysis is much simpler.

Our continuous Lyapunov function to analyze the class of low-resolution ODEs (\ref{eq:low-res-ODE-SC}) for Proposition \ref{pro:low-res-ODE-SC}
is extended from the one proposed in \cite{wilson2021lyapunov} for $(c_0,c_1)=(1,2)$ (i.e., the low-resolution ODE of NAG-SC).
Furthermore, we construct suitable Lyapunov functions to
analyze the high-resolution ODEs of NAG-SC and heavy-ball (which admit the same low-resolution ODE).
Compared with the ones used in \cite{shi2021understanding}, our Lyapunov functions lead to a sharper convergence bound
for NAG-SC and heavy-ball (Proposition \ref{pro:high-res-ODE NAG and HB}).

\textbf{General convex setting.}
We construct the discrete Lyapunov (\ref{eq:lyapunov-C}) to establish the convergence of (\ref{eq:extended-NAG-C}) or equivalently (\ref{eq:extended-NAG-C-xyz})
including those in \cite{beck2017first} and \cite{shi2021understanding} as special cases. The auxiliary-energy term $\|z_{k+1}-x^*\|^2/2$ in (\ref{eq:lyapunov-C}) is motivated by \cite{su2016differential} for NAG-C and \cite{shi2021understanding} for (\ref{eq:extended-NAG-C-shi}),
a sub-class of (\ref{eq:extended-NAG-C}) with $\gamma_k=1$, $\beta_k=\beta/\gamma$, and $\alpha_k = (k+r)/r$ for $r\ge 2$.
However, the potential-energy term in our Lyapunov (\ref{eq:lyapunov-C}) differs from the related ones in \cite{su2016differential} and \cite{shi2021understanding}.
Moreover, \cite{shi2021understanding} analyzed the three cases, $r=2$ and $\beta\le 1$, $r=2$ and $\beta>1$,
and $r>2$, separately. It seems difficult to extend their case-by-case analysis to cover the more general results in our Theorem \ref{thm:converge-c}.

The continuous Lyapunov function for the class of high-resolution ODEs (\ref{eq:high-res-ODE-C}) is the same as in \cite{shi2021understanding} when $r>2$, but involves a technical modification when $r=2$. The modification helps to establish the convergence bounds in Proposition \ref{pro:high-res-ODE-C}
for both $r>2$ and $r=2$, whereas similar bounds are provided in \cite{shi2021understanding} only explicitly for $r >2$.

\subsection{Strongly convex setting} \label{sec:outline-SC}

We provide a unified Lyapunov analysis to establish Theorems~\ref{thm:converge-sc-const}--\ref{thm:converge-sc-1}
for the class of algorithms (\ref{eq:extended-NAG-SC}). Proof details are presented in Supplement Section \ref{sec:tech for section sc}.
Although our Lyapunov function, like most existing ones, is manually designed, our analysis proceeds in several structured steps.

\textbf{Step 1. Bounding the differencing of an auxiliary energy.}
The sequence $\{z_k\}$ plays a key role in our formulation of (\ref{eq:extended-NAG-SC}) in a three-variable form.
We identify $\frac{\mu}{2} \|z_k-x^*\|^2$ as an auxiliary-energy term and bound its differencing,
$\frac{1}{1-\nu\sqrt{q}}\frac{\mu}{2}\|z_{k+1}-x^*\|-\frac{\mu}{2}\|z_k-x^*\|^2$.

\begin{lem}\label{lem:bounding diff of kinetic-SC}
    Let $f\in \mathcal{S}^1_{\mu,L} $. For any $s>0$ such that $0\leq \eta s\leq 1/L$, $\nu\geq 0$, $1-\nu\sqrt{q}>0$, $\tau>0$ and $\zeta=1+(1-\tau)\sqrt{q}\geq 0$, the iterates of (\ref{eq:extended-NAG-SC}) satisfy that for $k\geq 1$,
    \begin{equation}\label{eq:bounding diff of kinetic-SC}
        \begin{split}
            &\left( \frac{\zeta\nu}{\tau^2}(\tau+\zeta\eta\sqrt{q})+\frac{\nu\sqrt{q}}{1-\nu\sqrt{q}}\right)(f(x_k)-f^*) - \frac{\nu s}{2}\left(\frac{\nu}{1-\nu\sqrt{q}}-\frac{\zeta\eta}{\tau}\right)\|\nabla f(x_k)\|^2 \\
            &\quad + \frac{1}{1-\nu\sqrt{q}}\frac{\mu}{2}\|z_{k+1}-x^*\|^2  \\
            &\leq \frac{\zeta\nu}{\tau^2}(\tau+\zeta\eta\sqrt{q})\left(f(x_{k-1})-f^*-\frac{\eta s}{2}\|\nabla f(x_{k-1})\|^2\right) + \frac{\mu}{2}\|z_k-x^*\|^2.
        \end{split}
    \end{equation}
\end{lem}

\vspace{.05in}
\textbf{Step 2. Constructing a discrete Lyapunov function.}
We define a Lyapunov function simply from the right-hand-side of (\ref{eq:bounding diff of kinetic-SC}): for $k\geq 0$,
\begin{equation}\label{eq:lyapunov-SC}
   V_{k+1} = \frac{\zeta\nu}{\tau^2}(\tau+\zeta\eta\sqrt{q})\left(f(x_{k})-f^*-\frac{\eta s}{2}\|\nabla f(x_{k})\|^2\right) + \frac{\mu}{2}\|z_{k+1}-x^*\|^2.
\end{equation}
We refer to the first term above as a potential energy, involving both $f(x_k)-f^*$ and $\|\nabla f(x_k)\|^2$,
and the second term $\frac{\mu}{2}\|z_k-x^*\|^2$ as an auxiliary energy.

\textbf{Step 3. Identifying sufficient conditions for Lyapunov contraction.} As expected from (\ref{eq:bounding diff of kinetic-SC}),
we further bound $V_{k+1} -  (1-\nu\sqrt{q})V_k$ and identify conditions such that a contraction inequality holds for the Lyapunov function:
$V_{k+1} -  (1-\nu\sqrt{q})V_k \le 0$.

\begin{lem}\label{lem:general conditions-SC}
    Define $\one$ and $\two$ as polynomials of $\sqrt{q}$, $\eta$, $\nu$ and $\tau$ (hence functions of $\sqrt{q}$) taking the following forms: (recall that $\zeta=1+(1-\tau)\sqrt{q}$)
    \begin{equation}\label{eq:def of one and two}
    \begin{split}
        \one &= \zeta\nu(\tau+\zeta\eta\sqrt{q})-\tau^2 \\
        &=\tau(\nu-\tau) + \nu(\eta-\tau(\tau-1))\sqrt{q}-2\eta\nu(\tau-1)q+\eta\nu(\tau-1)^2q^{\frac{3}{2}}, \\
        \two &= \tau(\nu\tau-2\zeta\eta)+\zeta\eta(\nu\tau-\zeta\eta)\sqrt{q} \\
        &= \tau(\nu\tau-2\eta) + \eta(-\eta + \nu\tau + 2\tau(\tau-1) )\sqrt{q} + \eta(\tau-1)(2\eta-\nu\tau)q-\eta^2(\tau-1)^2q^{\frac{3}{2}}.
    \end{split}
    \end{equation}
    Under the condition in Lemma \ref{lem:bounding diff of kinetic-SC}, if one of the following (mutually exclusive) conditions holds:
    (i) $\one> 0$ and $\one+\sqrt{q}\two\leq 0$, or (ii) $\two> 0$ and $\frac{\mu}{L}\one + \sqrt{q}\two \leq 0$, or (iii) $\one\leq 0$ and $\two\leq 0$,
    then we have the contraction inequality $V_{k+1}\leq (1-\nu\sqrt{q})V_k$ for  $k\geq 1$.
\end{lem}

\textbf{Step 4. Verifying the contraction conditions and completing the analysis.}
For completing the analysis, the final step is to show that the sufficient conditions for Lyapunov contraction in Lemma~\ref{lem:general conditions-SC}
are satisfied under the conditions included in Theorems~\ref{thm:converge-sc-const}--\ref{thm:converge-sc-1}.
This step can be algebraically tedious, and the details are left to the Supplement.

The contraction inequality, if verified, directly leads to a convergence bound $O((1-C_1\sqrt{\mu s})^k)$ as stated in (\ref{eq:exponential bound in q}).
However, to fulfil the conditions in Lemma \ref{lem:general conditions-SC},
we find that two ranges of stepsize $s$ (or $q$) are allowed:
$0<s\lesssim \mu/L^2$ (or $0<q\lesssim \mu^2/L^2$) and $0<s\lesssim 1/L$ (or $0<q\lesssim \mu/L$).
As discussed after Theorems~\ref{thm:converge-sc-const}, the upper bound of feasible $s$
determines whether the usual or accelerated convergence is achieved, in terms of the dependency on $\mu/L$.

\subsection{General convex setting} \label{sec:outline-C}

We provide a unified Lyapunov analysis to establish Theorem~\ref{thm:converge-c}
for the class of algorithms (\ref{eq:extended-NAG-C}) or equivalently (\ref{eq:extended-NAG-C-xyz}).
Proof details are presented in Supplement Section \ref{sec:tech for section c}.
Our analysis proceeds in several structured steps, similarly as  in the strongly convex setting (Section \ref{sec:outline-SC}).

\textbf{Step 1. Bounding the differencing of an auxiliary energy.}
From our formulation of the three-variable form (\ref{eq:extended-NAG-C-xyz}),
we identify $\frac{1}{2} \|z_k-x^*\|^2$ as an auxiliary-energy term and bound its differencing,
$\frac{1}{2} \|z_{k+1}-x^*\|^2 - \frac{1}{2} \|z_k-x^*\|^2 $.

\begin{lem}\label{lem:bounding diff of kinetic-C}
    Let $f\in\mathcal{F}^1_L $. Then for each $k\geq 1$ such that $\a_k\geq 1$ and $\widetilde \a_k\geq 0$, the iterates of (\ref{eq:extended-NAG-C-xyz}) satisfy
    \begin{equation}\label{eq:bounding diff of kinetic-C}
        \begin{split}
            &\a_k\widetilde \a_k s(f(x_k)-f^*)-\frac{\widetilde \a_k^2 s^2}{2} \|\nabla f(x_k)\|^2 + \frac{1}{2}\|z_{k+1}-x^*\|^2
        \\ &\leq  \widetilde \a_k(\a_k-1) s \left((f(x_{k-1})-f^*) -\frac{(2-\b_{k-1}Ls)\b_{k-1}s}{2}\|\nabla f(x_{k-1})\|^2\right) + \frac{1}{2}\|z_{k}-x^*\|^2.
        \end{split}
    \end{equation}
\end{lem}

\vspace{.05in}
\textbf{Step 2. Constructing a discrete Lyapunov function.}
We define a Lyapunov function simply as the left-hand-side of (\ref{eq:bounding diff of kinetic-C}) up to a scalar sequence $\{\w_k\}$:
for $k\geq 0$,
\begin{equation}\label{eq:lyapunov-C}
    V_{k+1} = \w_{k+1} \left(\a_k\widetilde \a_k s(f(x_k)-f^*)-\frac{\widetilde \a_k^2 s^2}{2} \|\nabla f(x_k)\|^2 + \frac{1}{2}\|z_{k+1}-x^*\|^2  \right),
\end{equation}
The sequence $\{\w_k\}$ are introduced to later deal with the mismatching of coefficients on two sides of (\ref{eq:bounding diff of kinetic-C}).
Similarly as in Section \ref{sec:outline-SC}, we refer to the term involving $f(x_k)-f^*$ and $\|\nabla f(x_k)\|^2$ as a potential energy,
and the term $\frac{1}{2}\|z_k-x^*\|^2$ as an auxiliary energy.

\textbf{Step 3. Identifying sufficient conditions for Lyapunov contraction.} As expected from (\ref{eq:bounding diff of kinetic-C}),
we further bound $V_{k+1} - V_k$ and identify conditions such that a contraction inequality holds for the Lyapunov function:
$V_{k+1} - V_k \le 0$ or $V_{k+1} - V_k \lesssim -k^2s^2\|\nabla f(x_{k-1})\|^2$.

\begin{lem}\label{lem:general conditions-C}
    Define $\one$ and $\two$ as follows:
    \begin{equation*}
        \begin{split}
            \one &=  \w_k\a_{k-1}\widetilde\a_{k-1}-\w_{k+1}\widetilde\a_k(\a_k-1), \\
            \two  &= \w_{k+1}\widetilde\a_k(\a_k-1)\b_{k-1}(2-\b_{k-1}Ls)-\w_k\widetilde\a_{k-1}^2.
        \end{split}
    \end{equation*}
    For any $k\geq 1$ such that
    $\alpha_k\geq 1$, $\widetilde\alpha_k \geq 0$,
    $\omega_k\geq \omega_{k+1}$, $\one \geq 0$, and $ \two \geq 0$, the Lyapunov function (\ref{eq:lyapunov-C}) satisfies $V_{k+1}\leq V_{k}$. If further $ \two \ge Ck^2 $ for a constant $C>0$, then $V_{k+1}- V_{k}\le -\frac{C}{2} k^2 s^2 \|\nabla f(x_{k-1})\|^2$.
\end{lem}

\textbf{Step 4. Verifying the contraction conditions and completing the analysis.}
The final step  is to show that the sufficient conditions for Lyapunov contraction in Lemma~\ref{lem:general conditions-C}
are satisfied under the conditions included in Theorem \ref{thm:converge-c}.
The details are left to the Supplement.

The contraction inequality $V_{k+1} - V_k \le 0$ leads to the convergence bound (\ref{eq:optimal bound-c})
for the objective gap, whereas the contract inequality $V_{k+1} - V_k \lesssim -k^2s^2\|\nabla f(x_{k-1})\|^2$
leads to the inverse cubic rate (\ref{eq:inverse cubic rate}) for the squared gradient norm.
Unlike in the strongly convex setting, the stepsize $s$ can be simply set in the range $0<s\lesssim 1/L$
to fulfil the conditions in Lemma~\ref{lem:general conditions-C}.

\section{Conclusion}

Our work contributes to understanding acceleration of first-order algorithms for convex optimization.
Compared with the ODE-based approach, we directly formulate discrete algorithms as general as we can
and establish sufficient conditions for accelerated convergence using discrete Lyapunov functions.
We point out currently notable gaps between the convergence properties of the corresponding algorithms and ODEs.
We propose the Hamiltonian assisted gradient method, HAG, and demonstrate meaningful and unified interpretations of our acceleration conditions.
Future work is needed to address various open questions, including to what extent our sufficient conditions are also necessary,
further understanding the construction of discrete Lyapunov functions as well as continuous ones,
and resolving the current gaps between discrete algorithms and ODEs.

\bibliographystyle{apacite}
\bibliography{Acc_grad_c}
%%%%%%%%%%%%%%%%%%%%%%%%%%%%%%%%%%%%%%%%%%%%%%%%%%%%%%%%%%%%%%%%%%%%%%%%%%%%%%%%%%%%%%%%%%
%%%---------------------------------------------------------------------------------------
%%%%%%%%%%%%%%%%%%%%%%%%%%%%%%%%%%%%%%%%%%%%%%%%%%%%%%%%%%%%%%%%%%%%%%%%%%%%%%%%%%%%%%%%%%
\clearpage

\setcounter{page}{1}

\setcounter{section}{0}
\setcounter{equation}{0}

\setcounter{figure}{0}
\setcounter{table}{0}

\renewcommand{\theequation}{S\arabic{equation}}
\renewcommand{\thesection}{\Roman{section}}

\renewcommand\thefigure{S\arabic{figure}}
\renewcommand\thetable{S\arabic{table}}

\setcounter{lem}{0}
\renewcommand{\thelem}{S\arabic{lem}}

\begin{center}
{\Large Supplementary Material for}

{\Large ``Understanding Accelerated Gradient Methods: Lyapunov Analyses and Hamiltonian 
Assisted Interpretations"}

\vspace{.1in} {\large Penghui Fu and Zhiqiang Tan}
\end{center}

\section{Technical details in Section \ref{sec:acc-sc}}\label{sec:tech for section sc}

\subsection{Reformulation of NAG-SC and TMM}\label{subsec:reformulation of NAG-SC and TMM}
First, we show that NAG-SC (\ref{eq:NAG-SC}) with (\ref{eq:NAG-sc-z}) can be reformulated as (\ref{eq:reformulate-NAG-SC}), i.e., (\ref{eq:extended-NAG-SC}) with $\eta=\nu=\tau=1$. Apparently (\ref{subeq:NAG-SC-y}) is the same as (\ref{subeq:reformulate-NAG-SC-y}). It suffices to derive (\ref{subeq:reformulate-NAG-SC-z}) and (\ref{subeq:reformulate-NAG-SC-x}) from (\ref{eq:NAG-SC}). Because $x_0=y_0=z_0$, (\ref{eq:NAG-sc-z}) also holds for $z_0$. Therefore, for $k\geq 0$ we have
\begin{equation*}
    z_{k} = \frac{1+\sqrt{q}}{\sqrt{q}} x_{k} + \left(1-\frac{1+\sqrt{q}}{\sqrt{q}}\right)y_{k},
\end{equation*}
which implies
\begin{equation}\label{eq:supp-NAG-SC-x}
    x_k = \frac{\sqrt{q}}{1+\sqrt{q}}z_k + \frac{1}{1+\sqrt{q}}y_k,
\end{equation}
and
\begin{equation}\label{eq:supp-NAG-SC-y}
    y_k = (1+\sqrt{q})x_k - \sqrt{q}z_k.
\end{equation}
Eq.~(\ref{eq:supp-NAG-SC-x}) with $k$ replaced by $k+1$ becomes (\ref{subeq:reformulate-NAG-SC-x}). Substituting $x_{k+1}$ from (\ref{eq:supp-NAG-SC-x}) into (\ref{subeq:NAG-SC-single}), we have
$$
\frac{\sqrt{q}}{1+\sqrt{q}}z_{k+1} + \frac{1}{1+\sqrt{q}}y_{k+1}= y_{k+1} + \frac{1-\sqrt{q}}{1+\sqrt{q}}(y_{k+1}-y_k),
$$
which yields
$$
z_{k+1} = \frac{1}{\sqrt{q}}y_{k+1}-\frac{1-\sqrt{q}}{\sqrt{q}}y_k.
$$
Substituting $y_{k+1}$ from (\ref{subeq:NAG-SC-y}) into the above display, we have
$$
z_{k+1} = \frac{1}{\sqrt{q}}(x_k-s\nabla f(x_k))-\frac{1-\sqrt{q}}{\sqrt{q}}y_k,
$$
which together with (\ref{eq:supp-NAG-SC-y}) gives
\begin{align*}
    z_{k+1} &= \frac{1}{\sqrt{q}}(x_k-s\nabla f(x_k))-\frac{1-\sqrt{q}}{\sqrt{q}}\big((1+\sqrt{q})x_k - \sqrt{q}z_k\big), \\
    &= \sqrt{q}\left(x_k-\frac{1}{\mu}\nabla f(x_k)\right) + (1-\sqrt{q})z_k,
\end{align*}
which is (\ref{subeq:reformulate-NAG-SC-z}).

Second, we show that TMM (\ref{eq:TMM}), i.e., (\ref{eq:extended-NAG-SC}) with $\eta=\nu=1$ and $\tau=2$, can be reformulated as follows,
with (\ref{subeq:TMM-reformulate-z}) same as (\ref{eq:TMM-z-auxiliary}) for $z_{k+1}$:
\begin{subequations}\label{eq:TMM-reformulate}
    \begin{align}
        y_{k+1} &= x_k - s\nabla f(x_k),\label{subeq:TMM-reformulate-y}\\
        x_{k+1} &= x_k - (2-\sqrt{q})s\nabla f(x_k) + \frac{(1-\sqrt{q})^2}{1+\sqrt{q}}(y_{k+1}-y_k),\label{subeq:TMM-reformulate-x}\\
        z_{k+1} &= \frac{1+\sqrt{q}}{2\sqrt{q}} x_{k+1} + \left(1-\frac{1+\sqrt{q}}{2\sqrt{q}}\right)y_{k+1} . \label{subeq:TMM-reformulate-z}
    \end{align}
\end{subequations}
Apparently (\ref{eq:TMM-y}) is the same as (\ref{subeq:TMM-reformulate-y}). It suffices to derive (\ref{subeq:TMM-reformulate-x}) and (\ref{subeq:TMM-reformulate-z}). 
Solving (\ref{eq:TMM-x}) for $z_{k+1}$ gives (\ref{subeq:TMM-reformulate-z}). 
Because $x_0=y_0=z_0$, (\ref{subeq:TMM-reformulate-z}) also holds for $z_0$. That is, for all $k\geq 0$,
\begin{equation*}
    z_{k} = \frac{1+\sqrt{q}}{2\sqrt{q}} x_{k} + \left(1-\frac{1+\sqrt{q}}{2\sqrt{q}}\right)y_{k},
\end{equation*}
which together with (\ref{eq:TMM-z}) yields
\begin{equation*}
    \begin{split}
        z_{k+1} &= \sqrt{q}\left(x_k-\frac{1}{\mu}\nabla f(x_k)\right) + (1-\sqrt{q}) \left[\frac{1+\sqrt{q}}{2\sqrt{q}} x_{k} + \left(1-\frac{1+\sqrt{q}}{2\sqrt{q}}\right)y_{k} \right]\\
        &=\frac{1+q}{2\sqrt{q}}x_k - \frac{\sqrt{q}}{\mu}\nabla f(x_k) - \frac{(1-\sqrt{q})^2}{2\sqrt{q}}y_k.
    \end{split}
\end{equation*}
To derive (\ref{subeq:TMM-reformulate-x}), substituting the above display into (\ref{eq:TMM-x}) and using (\ref{eq:TMM-y}) yields
\begin{equation*}
    \begin{split}
        x_{k+1} &= \frac{2\sqrt{q}}{1+\sqrt{q}}\left(\frac{1+q}{2\sqrt{q}}x_k - \frac{\sqrt{q}}{\mu}\nabla f(x_k) - \frac{(1-\sqrt{q})^2}{2\sqrt{q}}y_k\right) + \frac{1-\sqrt{q}}{1+\sqrt{q}}y_{k+1} \\
        &= \frac{1+q}{1+\sqrt{q}}x_k - \frac{2s}{1+\sqrt{q}}\nabla f(x_k) + \frac{(1-\sqrt{q})^2}{1+\sqrt{q}}(y_{k+1}-y_k) + \frac{\sqrt{q}(1-\sqrt{q})}{1+\sqrt{q}}y_{k+1} \\
        &= \frac{1+q}{1+\sqrt{q}}x_k - \frac{2s}{1+\sqrt{q}}\nabla f(x_k) + \frac{(1-\sqrt{q})^2}{1+\sqrt{q}}(y_{k+1}-y_k) \\
        &\quad + \frac{\sqrt{q}(1-\sqrt{q})}{1+\sqrt{q}}(x_k-s\nabla f(x_k)) \\
        &= x_k - (2-\sqrt{q})s\nabla f(x_k) + \frac{(1-\sqrt{q})^2}{1+\sqrt{q}}(y_{k+1}-y_k),
    \end{split}
\end{equation*}
which is (\ref{subeq:TMM-reformulate-x}).

\subsection{Proofs of Lemmas~\ref{lem:bounding diff of kinetic-SC}--\ref{lem:general conditions-SC} in Section \ref{sec:outline-SC}}

\textbf{Proof of Lemma~\ref{lem:bounding diff of kinetic-SC}.}
By (\ref{subeq:extended-NAG-SC-z}), we have
\begin{align}
    &\frac{1}{1-\nu\sqrt{q}}\frac{\mu}{2}\|z_{k+1}-x^*\|^2 - \frac{\mu}{2}\|z_{k}-x^*\|^2  \notag\\
        &=\frac{1}{1-\nu\sqrt{q}}\frac{\mu}{2}\left\|\nu \sqrt{q}(x_k-\frac{1}{\mu}\nabla f(x_k)-x^*) + (1-\nu\sqrt{q}) (z_k - x^*)\right\|^2 - \frac{\mu}{2}\|z_{k}-x^*\|^2 \notag \\
        &= -\frac{\mu\nu\sqrt{q}}{2}\|z_k-x^*\|^2 + \mu\nu\sqrt{q}\langle z_k-x^*, x_k-\frac{1}{\mu}\nabla f(x_k)-x^* \rangle \notag\\
        &\quad + \frac{\nu^2 q}{1-\nu\sqrt{q}}\frac{\mu}{2}\left\|x_k-\frac{1}{\mu}\nabla f(x_k)-x^*\right\|^2 \notag \\
        &= -\frac{\mu\nu\sqrt{q}}{2}\left\|(z_k-x^*)-(x_k-\frac{1}{\mu}\nabla f(x_k)-x^*)\right\|^2 \notag\\
        &\quad + \frac{\nu \sqrt{q}}{1-\nu\sqrt{q}}\frac{\mu}{2}\left\|x_k-\frac{1}{\mu}\nabla f(x_k)-x^*\right\|^2 \notag \\
        &= -\frac{\mu\nu\sqrt{q}}{2}\left\|z_k-x_k+\frac{1}{\mu}\nabla f(x_k)\right\|^2 + \frac{\nu \sqrt{q}}{1-\nu\sqrt{q}}\frac{\mu}{2}\left\|x_k-\frac{1}{\mu}\nabla f(x_k)-x^*\right\|^2. \label{eq:diff of kinetic-SC-0}
\end{align}
Next, we bound the two terms in the above display separately.

For the second term, by the $\mu$-strong convexity of $f$,
$$
f^* - f(x_k) \geq \langle \nabla f(x_k),x^*-x_k \rangle + \frac{\mu}{2}\|x_k-x^*\|^2,
$$
from which we have
\begin{equation}\label{eq:diff of kinetic-SC-1}
    \begin{split}
        &\frac{\nu \sqrt{q}}{1-\nu\sqrt{q}}\frac{\mu}{2}\left\|x_k-\frac{1}{\mu}\nabla f(x_k)-x^*\right\|^2 \\
        &=\frac{\nu \sqrt{q}}{1-\nu\sqrt{q}}\left(\frac{\mu}{2}\|x_k-x^*\|^2-\langle x_k-x^*,\nabla f(x_k) \rangle + \frac{1}{2\mu}\|\nabla f(x_k)\|^2\right) \\
        &\leq \frac{\nu \sqrt{q}}{1-\nu\sqrt{q}}\left(f^*-f(x_k)+\frac{1}{2\mu}\|\nabla f(x_k)\|^2 \right).
    \end{split}
\end{equation}

For the first term, solving (\ref{subeq:extended-NAG-SC-x}) for $z_{k+1}$ we obtain
$$
z_{k+1} = x_{k+1} + \frac{1+(1-\tau)\sqrt{q}}{\tau\sqrt{q}}(x_{k+1}-y_{k+1}) = x_{k+1} + \frac{\zeta}{\tau\sqrt{q}}(x_{k+1}-y_{k+1}),
$$
which together with (\ref{subeq:extended-NAG-SC-y}) yields
$$
z_{k+1} = x_{k+1} + \frac{\zeta}{\tau\sqrt{q}}\left(x_{k+1}-x_{k} + \eta s \nabla f(x_{k})\right).
$$
Hence for $k\geq 1$ we have
\begin{equation}\label{eq:diff of kinetic-SC-2}
    \begin{split}
        &-\frac{\mu\nu\sqrt{q}}{2}\left\|z_k-x_k+\frac{1}{\mu}\nabla f(x_k)\right\|^2 \\
        &= -\frac{\mu\nu\sqrt{q}}{2}\left\|\frac{\zeta}{\tau\sqrt{q}}(x_{k}-x_{k-1}) + \frac{\zeta\eta s}{\tau\sqrt{q}} \nabla f(x_{k-1})+\frac{1}{\mu}\nabla f(x_k)\right\|^2 \\
        &= - \frac{\zeta^2\nu}{\tau^2\sqrt{q}}\frac{\mu}{2}\|x_k-x_{k-1}\|^2 - \frac{\zeta^2\eta^2\nu\sqrt{q}}{\tau^2}\frac{s}{2}\|\nabla f(x_{k-1})\|^2 - \frac{\nu\sqrt{q}}{2\mu}\|\nabla f(x_k)\|^2 \\
        &\quad - \frac{\zeta^2\eta\nu\sqrt{q}}{\tau^2}\langle \nabla f(x_{k-1}), x_k-x_{k-1}\rangle - \frac{\zeta\nu}{\tau}\langle \nabla f(x_{k}), x_k-x_{k-1}\rangle - \frac{\zeta\eta\nu}{\tau}s\langle\nabla f(x_{k-1}),\nabla f(x_{k})\rangle.
    \end{split}
\end{equation}
By the $L$-smoothness of $f$, we have
\begin{equation}
    \begin{dcases}
        f(x_{k-1}) \geq f(x_k) + \langle \nabla f(x_k),x_{k-1}-x_k \rangle + \frac{1}{2L}\|\nabla f(x_{k})-\nabla f(x_{k-1})\|^2, \\
        f(x_k) \leq f(x_{k-1}) + \langle\nabla f(x_{k-1}),x_k-x_{k-1} \rangle + \frac{L}{2}\|x_k-x_{k-1}\|^2.
    \end{dcases}
\end{equation}
Therefore,
\begin{equation}\label{eq:diff of kinetic-SC-3}
    \begin{split}
        &- \frac{\zeta^2\eta\nu\sqrt{q}}{\tau^2}\langle \nabla f(x_{k-1}), x_k-x_{k-1}\rangle - \frac{\zeta\nu}{\tau}\langle \nabla f(x_{k}),x_k-x_{k-1}\rangle  \\
        &\leq - \frac{\zeta^2\eta\nu\sqrt{q}}{\tau^2}\left(f(x_k)-f(x_{k-1})-\frac{L}{2}\|x_k-x_{k-1}\|^2 \right) \\
        &\quad -\frac{\zeta\nu}{\tau}\left(f(x_k)-f(x_{k-1})+\frac{1}{2L}\|\nabla f(x_k)-\nabla f(x_{k-1})\|^2 \right) \\
        &= \frac{\zeta\nu}{\tau^2}(\tau+\zeta\eta\sqrt{q})(f(x_{k-1})-f(x_k)) + \frac{\zeta^2\eta\nu\sqrt{q}}{\tau^2}\frac{L}{2}\|x_k-x_{k-1}\|^2 \\
        &\quad - \frac{\zeta\nu}{\tau}\frac{1}{2L}\|\nabla f(x_k)-\nabla f(x_{k-1})\|^2.
    \end{split}
\end{equation}
By the cosine rule,
\begin{equation}\label{eq:diff of kinetic-SC-4}
    -\frac{\zeta\eta\nu s}{\tau}\langle\nabla f(x_{k-1}),\nabla f(x_{k})\rangle = -\frac{\zeta\eta\nu s}{2\tau}\Big(\|\nabla f(x_{k-1})\|^2 + \|\nabla f(x_{k})\|^2 - \|\nabla f(x_k)-\nabla f(x_{k-1})\|^2 \Big).
\end{equation}
Combining (\ref{eq:diff of kinetic-SC-2}), (\ref{eq:diff of kinetic-SC-3}) and (\ref{eq:diff of kinetic-SC-4})
and noting $0\leq\eta s\leq 1/L$ by assumption, we obtain
\begin{equation}\label{eq:diff of kinetic-SC-5}
    \begin{split}
        &-\frac{\mu\nu\sqrt{q}}{2}\left\|z_k-x_k+\frac{1}{\mu}\nabla f(x_k)\right\|^2 \\
        &\leq -\frac{\zeta^2\mu\nu}{2\tau^2\sqrt{q}}(1-\eta Ls)\|x_k-x_{k-1}\|^2 - \frac{\zeta\nu}{2\tau L}(1-\eta Ls)\|\nabla f(x_k)-\nabla f(x_{k-1})\|^2 \\
        &\quad +  \frac{\zeta\nu}{\tau^2}(\tau+\zeta\eta\sqrt{q})(f(x_{k-1})-f(x_k)) - \frac{\zeta\eta\nu s}{2\tau^2}(\tau+\zeta\eta\sqrt{q})\|\nabla f(x_{k-1})\|^2 \\
        &\quad - \frac{\nu\sqrt{q}}{2\mu\tau}(\tau+\zeta\eta\sqrt{q})\|\nabla f(x_k)\|^2 \\
        &\leq \frac{\zeta\nu}{\tau^2}(\tau+\zeta\eta\sqrt{q})(f(x_{k-1})-f(x_k)) - \frac{\zeta\eta\nu s}{2\tau^2}(\tau+\zeta\eta\sqrt{q})\|\nabla f(x_{k-1})\|^2 \\
        &\quad- \frac{\nu\sqrt{q}}{2\mu\tau}(\tau+\zeta\eta\sqrt{q})\|\nabla f(x_k)\|^2.
    \end{split}
\end{equation}
Collecting (\ref{eq:diff of kinetic-SC-0}), (\ref{eq:diff of kinetic-SC-1}) and (\ref{eq:diff of kinetic-SC-5}) completes the proof.
\hfill $\Box$ \vspace{.1in}

\noindent\textbf{Proof of Lemma~\ref{lem:general conditions-SC}.}
The three sets of conditions are mutually exclusive because if $\one >0$ and $\two >0$, then $\one + \sqrt{q}\two >0$ and $\frac{\mu}{L}\one + \sqrt{q}\two> 0$. Applying (\ref{eq:bounding diff of kinetic-SC}) to the Lyapunov function (\ref{eq:lyapunov-SC}) and using the fact that $\frac{1}{2L}\|\nabla f(x_k)\|^2\leq f(x_k)-f^*\leq \frac{1}{2\mu}\|\nabla f(x_k)\|^2$, we have
\begin{equation*}
    \begin{split}
        &V_{k+1} - (1-\nu\sqrt{q})V_k \\
        &\leq \frac{\nu}{2\tau^2}\Big(2\sqrt{q}(\one) (f(x_k)-f^*) + s(\two) \|\nabla f(x_k)\|^2\Big)  \\
        &\leq \begin{dcases}
            \frac{\nu}{2\tau^2}\frac{\sqrt{q}}{\mu}\Big(\one + \sqrt{q}\two\Big)\|\nabla f(x_k)\|^2, & \text{if } \one\geq 0;\\
            \frac{\nu}{2\tau^2}\frac{\sqrt{q}}{\mu}\Big( \frac{\mu}{L}\one + \sqrt{q}\two \Big)\|\nabla f(x_k)\|^2, & \text{if } \one\leq 0;\\
            \frac{\nu}{\tau^2}\frac{\sqrt{q}L}{\mu}\Big( \frac{\mu}{L}\one + \sqrt{q}\two \Big)(f(x_k)-f^*), & \text{if } \two\geq 0; \\
            \frac{\nu}{\tau^2}\sqrt{q}\Big(\one + \sqrt{q}\two\Big) (f(x_k)-f^*), & \text{if } \two\leq 0.
        \end{dcases}
    \end{split}
\end{equation*}
The rest is straightforward.
\hfill $\Box$ \vspace{.1in}

\subsection{Proofs of Theorems~\ref{thm:converge-sc-const}--\ref{thm:converge-sc-1}} \label{sec:prf-thm-SC}

To prepare for the proofs of Theorems~\ref{thm:converge-sc-const}--\ref{thm:converge-sc-1}, we show that under Assumption \ref{ass:regular func of sqrt-q},
the leading coefficients in the
Taylor expansions of $\one$ and $\two$ in $\sqrt{q}$ can be used to verify the conditions in Lemma \ref{lem:general conditions-SC}.
Throughout, the range $0< q\lesssim \mu^2/L^2$ or  $0< q\lesssim \mu/L$ is interpreted as, respectively,
$0<q\leq C_0 \mu^2/L^2$ or $0<q\leq C_0 \mu/L$ for a constant $C_0 >0$.

\begin{lem}\label{lem:general conditions-coef-SC}
    Under Assumption \ref{ass:regular func of sqrt-q}, denote the Taylor expansions of $\one$ and $\two$  as
    $$
    \one = \sum_{n=0}^\infty a_n (\sqrt{q})^n,\quad \two = \sum_{m=0}^\infty b_m (\sqrt{q})^m,
    $$
    where $\{a_n\}_{n\geq 0}$ and $\{b_m\}_{m\geq 0}$ are real sequences. When $\one$ is not constant $0$, define $N$ as the minimal of $n$ such that $a_n\neq 0$. Define $M$ in a similar manner for $\two$. Then    \\
    \indent (ia) If $a_N >0$, $b_M<0$ and $M\leq N-2$, then $\one> 0$ and $\one+\sqrt{q}\two\leq 0$ for $0< q\lesssim \mu/L$.  \\
    \indent (ib) If $a_N >0$, $b_M<0$, $M=N-1$, and the first nonzero element of $\{a_n+b_{n-1}\}_{n \geq N}$ (i.e., the first nonzero coefficient in the expansion of $\one+\sqrt{q}\two$) is negative or the entire sequence is $0$ (i.e., $\one+\sqrt{q}\two\equiv 0$), then $\one> 0$ and $\one+\sqrt{q}\two\leq 0$ for $0< q\lesssim \mu/L$. \\
    \indent (iia) If $a_N <0$, $b_M >0$ and $M=N$, then $\two> 0$ and $\frac{\mu}{L}\one + \sqrt{q}\two \leq 0$ for $0< q\lesssim \mu^2/L^2$. \\
    \indent (iib) If $a_N <0$, $b_M >0$ and $M\geq N+1$, then $\two> 0$ and $\frac{\mu}{L}\one + \sqrt{q}\two \leq 0$ for $0< q\lesssim \mu/L$. \\
    \indent (iii) If $a_N<0$ (or $N$ does not exist) and $b_M<0$ (or $M$ does not exist), then $\one\leq 0$ and $\two\leq 0$ for $0< q\lesssim \mu/L$.
\end{lem}

\begin{prf}
    With $0< \mu/L \le 1$,  $q$ can be made sufficiently small by picking $C_0$
    in the range $0<q\leq C_0 \mu^2/L^2$ or $0<q\leq C_0 \mu/L$. Hence it suffices to study the leading terms of $\one$ and $\two$.

    For (ia) and (ib), $a_N>0$ ensures that $\one \sim a_N(\sqrt{q})^N >0$. For (ia),
    with $M+1<N$ and $b_M<0$, we have $\one+\sqrt{q}\two \sim b_M(\sqrt{q})^{M+1}< 0$.
    For (ib), with $M+1=N$, we have $\one + \sqrt{q}\two = \sum_{n=N} (a_n + b_{n-1})(\sqrt{q})^n\leq 0$ if the first nonzero element of $\{a_n+b_{n-1}\}_{n \geq N}$ is negative or the entire sequence is $0$.

    For (iia) and (iib), $b_M>0$ ensures that $\two\sim b_M(\sqrt{q})^M >0$. With $a_N <0$ and $b_M >0$, we have that
    for sufficiently small $q>0$, $\one < \frac{a_N}{2}(\sqrt{q})^{N}$ and $\two < 2b_M(\sqrt{q})^M$. Hence, $\frac{\mu}{L}\one + \sqrt{q}\two< \frac{\mu a_N}{2L}(\sqrt{q})^{N} + 2b_M(\sqrt{q})^{M+1} = (\sqrt{q})^N (\frac{\mu a_N}{2L}+ 2b_M(\sqrt{q})^{M+1-N})$. For (iia), $M=N$ and $0< q\leq C_0^2\mu^2/L^2$ imply that $\frac{\mu}{L}\one + \sqrt{q}\two<(\sqrt{q})^N\frac{\mu}{2L}(a_N + 4C_0 b_M)<0$ by
    picking sufficiently small $C_0$ with $a_N<0$. For (iib), $M\geq N+1$ and $0< q\leq C_0 \mu/L$ imply $\frac{\mu}{L}\one + \sqrt{q}\two<(\sqrt{q})^N (\frac{\mu a_N}{2L}+ 2b_M(\sqrt{q})^{2})\leq (\sqrt{q})^N\frac{\mu}{2L}(a_N + 4C_0 b_M)<0$
    again by picking sufficiently small $C_0$ with $a_N<0$.

    The case (iii) is straightforward to verify. If $N$ (or $M$) does not exist, then $\one\equiv 0$ (or $\two\equiv 0$).

    We notice that the constant $C_0$ in the range of $q$ is picked, depending only on $\{a_n\}$ and $\{b_m\}$, which are determined by the algorithm parameters $\tilde \eta$, $\tilde \nu$ and $\tilde \tau$.
\end{prf}

Next, we show Theorems~\ref{thm:converge-sc-const}$^\ast$, \ref{thm:converge-sc-0}$^\ast$, and \ref{thm:converge-sc-1},
where Theorems~\ref{thm:converge-sc-const}$^\ast$ and \ref{thm:converge-sc-0}$^\ast$ are the same as
Theorems~\ref{thm:converge-sc-const} and \ref{thm:converge-sc-0} except with conditions (ia) and (iia) replaced by (ia$^\ast$) and (iia$^\ast$) as follows:
        \begin{itemize} \addtolength{\itemsep}{-.05in}
        \item[(ia$^\ast$)] $0 < \nu_0 < \tau_0$, and $0 \le \eta_0 < \nu_0\tau_0/2$;
        \item[(iia$^\ast$)] $0< \nu_0 < \tau_0$, and $\eta_0\geq\nu_0\tau_0/2$.
\end{itemize}
Conditions (ia) and (iia) are the symmetrized (hence weaker) versions of (ia$^\ast$) and (iia$^\ast$),
by allowing either $0 < \nu_0 < \tau_0$ or $0 < \tau_0 < \nu_0$.
To show Theorems~\ref{thm:converge-sc-const}$^\ast$, \ref{thm:converge-sc-0}$^\ast$, and \ref{thm:converge-sc-1},
it suffices to verify that the conditions in Lemma \ref{lem:general conditions-SC} are satisfied.
In fact, the contraction inequality in Lemma \ref{lem:general conditions-SC}
directly implies that that for $k\geq 1$,
\begin{align*}
 & \quad C (f(x_k)-f^*) \leq\frac{\zeta\nu}{\tau^2}(\tau+\zeta\eta\sqrt{q})(1-\eta Ls)(f(x_k)-f^*) \\
 & \leq V_{k+1} \leq (1-\nu\sqrt{q})^k V_1\leq (1-\frac{\nu_0}{2}\sqrt{q})^k V_1,
\end{align*}
for some constant $C>0$, where the first inequality holds by noting $\nu_0, \tau_0>0$ in each condition of Theorems~\ref{thm:converge-sc-const}$^\ast$, \ref{thm:converge-sc-0}$^\ast$, and \ref{thm:converge-sc-1}
and picking sufficiently small $C_0$ in $0 < q \le C_0 \mu/L$ such that, for example, $\nu \ge \nu_0/2$, $\tau_0/2\le \tau \le 2 \tau_0$, $\zeta \ge 1/2$,
$0\leq \eta s\leq 1/(2L)$ (i.e., $ 0 \le \eta q \le \mu/(2L)$) with $C=\nu_0/(16\tau_0)$.
Moreover, by the definition of $V_1$,
\begin{equation}\label{eq:bounding V1-SC}
    \begin{split}
        V_1 &\leq \frac{\zeta\nu}{\tau^2}(\tau+\zeta\eta\sqrt{q})(f(x_0)-f^*) + \frac{\mu}{2}\|z_1-x^*\|^2 \\
    &= \frac{\zeta\nu}{\tau^2}(\tau+\zeta\eta\sqrt{q})(f(x_0)-f^*) + \frac{\mu}{2}\|x_0 -x^*- \frac{\nu\sqrt{q}}{\mu}\nabla f(x_0)\|^2 \\
    &\lesssim L \|x_0-x^*\|^2 + \mu\|x_0-x^*\|^2 + s\|\nabla f(x_0)\|^2 \\
    &\leq (L+\mu+sL^2)\|x_0-x^*\|^2 \\
    &\lesssim L\|x_0-x^*\|^2.
    \end{split}
\end{equation}
Then $f(x_k)-f^*\lesssim L(1-\frac{\nu_0}{2}\sqrt{q})^k\|x_0-x^*\|^2$, which is (\ref{eq:exponential bound in q}).
The conditions in Lemma \ref{lem:bounding diff of kinetic-SC}, which are required in Lemma \ref{lem:general conditions-SC},
can be easily verified by noting $\nu_0, \tau_0>0$ and picking sufficiently small $C_0$.
Therefore, it remains to verify the conditions involving $\one$ and $\two$ in Lemma \ref{lem:general conditions-SC}.\vspace{.1in}

%%%%%%%%%%%%%%%%%%%%%%%%
\noindent\textbf{Proof of Theorems~\ref{thm:converge-sc-0}$^\ast$ and \ref{thm:converge-sc-1}.}
We apply Lemma \ref{lem:general conditions-coef-SC} to verify the conditions involving $\one$ and $\two$ in Lemma \ref{lem:general conditions-SC}.
The Taylor expansions of $\one$ and $\two$ up to $\sqrt{q}$-terms are
\begin{equation}\label{eq:expansion of one and two}
    \begin{split}
        \one &= \underbrace{\tau_0(\nu_0-\tau_0)}_{a_0}  + \big[\tau_0(\nu_1-\tau_1)+\tau_1(\nu_0-\tau_0)+\nu_0(\eta_0-\tau_0(\tau_0-1)) \big]\sqrt{q} + O(q), \\
        \two &= \underbrace{\tau_0(\nu_0\tau_0-2\eta_0)}_{b_0}+ \big[\tau_0(\nu_1\tau_1-2\eta_1)+\tau_1(\nu_0\tau_0-2\eta_0) \\
        &\quad +\eta_0\big(2\tau_0(\tau_0-1)+\nu_0\tau_0-\eta_0\big) \big]\sqrt{q} + O(q).
    \end{split}
\end{equation}
Consider the following scenarios.

\textbf{Scenario 1:} $\bm{\eta_0=0,\, 0<\nu_0,\tau_0}$. Then $b_M=b_0=\nu_0\tau_0^2>0$. Because $M=0$, only case (iia) in Lemma \ref{lem:general conditions-coef-SC} is feasible, which holds when $N=0$ and $a_N=a_0=\tau_0(\nu_0-\tau_0)<0$, i.e., $\nu_0<\tau_0$. To conclude, if $\eta_0=0$, $0<\nu_0<\tau_0$, then Lemma \ref{lem:general conditions-SC} holds for $0<q\lesssim \mu^2/L^2$.

Assume $\eta_0$, $\nu_0$ and $\tau_0$ are all positive. We notice that $0<\nu_0\leq \tau_0$ is necessary. Otherwise, $a_0=\tau_0(\nu_0-\tau_0)>0$ and $N=0$. Then only cases (ia) and (ib) in Lemma \ref{lem:general conditions-coef-SC} are feasible, which require $M\leq N-1$, contradicting $N=0$. To proceed, we further split $0<\nu_0\leq \tau_0$ into Scenario 2 ($0<\nu_0<\tau_0$) and Scenario 3 ($0<\nu_0=\tau_0$) as below.

\textbf{Scenario 2:} $\bm{\eta_0>0,\, 0<\nu_0<\tau_0}$. Then $a_0=\tau_0(\nu_0-\tau_0)<0$ and $N=0$. We notice that either $M$ does not exist (i.e., $\two\equiv 0$), or $M$ exists and has $b_M>0$ or $b_M<0$ for some $M\geq N=0$. Therefore, one of case (iia), case (iib), and case (iii) in Lemma \ref{lem:general conditions-coef-SC} is valid. Case (iia) holds if and only if $M=N=0$ and $b_M=b_0=\tau_0(\nu_0\tau_0-2\eta_0)>0$, i.e., $\eta_0< \frac{\nu_0\tau_0}{2}$. Therefore, if $\eta_0\geq \frac{\nu_0\tau_0}{2}$, case (iib) or case (iii) holds. To conclude, if $0<\eta_0<\frac{\nu_0\tau_0}{2}$ and $0<\nu_0<\tau_0$, then Lemma \ref{lem:general conditions-SC} holds for $0<q\lesssim \mu^2/L^2$. If $\eta_0\geq \frac{\nu_0\tau_0}{2}$ and $0<\nu_0<\tau_0$, then the range of $q$ is relaxed to $0<q\lesssim \mu/L$.

\textbf{Scenario 3:} $\bm{\eta_0>0,\, 0<\nu_0=\tau_0}$. In this scenario, (\ref{eq:expansion of one and two}) reduces to
\begin{equation*}
    \begin{split}
        \one &= \underbrace{\tau_0\big[\nu_1-\tau_1+\eta_0-\tau_0(\tau_0-1) \big]}_{a_1}\sqrt{q} + O(q), \\
        \two &= \underbrace{\tau_0(\tau^2_0-2\eta_0)}_{b_0}+ \underbrace{\big[\tau_0(\nu_1\tau_1-2\eta_1)+\tau_1(\tau^2_0-2\eta_0) +\eta_0\big(\tau_0(3\tau_0-2)-\eta_0\big) \big]}_{b_1}\sqrt{q} + O(q).
    \end{split}
\end{equation*}
Then $a_0=0$ and hence $N\geq 1$. We point out that $\eta_0\geq \frac{\tau_0^2}{2}$ is necessary. Otherwise, $b_0=\tau_0(\tau^2_0-2\eta_0)>0$ and hence $M=0$. Then only case (iia) or (iib) in Lemma \ref{lem:general conditions-coef-SC} is possible, which requires $M\geq N$. But this contradicts the fact that $N\geq 1$ and $M=0$. To proceed, we split $\eta_0\geq \frac{\tau_0^2}{2}$ into Scenario 3.1 ($\eta_0 > \frac{\tau_0^2}{2}$) and 3.2 ($\eta_0=\frac{\tau_0^2}{2}$) based on $\eta_0$.

\textbf{Scenario 3.1:} $\bm{\eta_0>\frac{\tau_0^2}{2},\, 0<\nu_0=\tau_0}$. Then $b_0=\tau_0(\tau^2_0-2\eta_0)<0$ and $M=0$. 
Case (ia), case (ib) and case (iii) in Lemma \ref{lem:general conditions-coef-SC} are each possible. We consider several special cases involving only $\eta_1$, $\nu_1$ and $\tau_1$.
\begin{itemize}
\item For case (ia) to hold, even higher-order coefficients are needed, and we skip this case.

\item For case (ib) to hold, let $N=M+1=1$ and $a_N=a_1=\tau_0\big[\nu_1-\tau_1+\eta_0-\tau_0(\tau_0-1) \big]>0$, which is equivalent to $\nu_1-\tau_1+\eta_0-\tau_0(\tau_0-1)>0$. Moreover, let $a_N+b_{N-1}=a_1+b_0=\tau_0\big[\nu_1-\tau_1+\eta_0-(\tau_0-1)\tau_0\big]+\tau_0(\tau^2_0-2\eta_0)=\tau_0[\nu_1-\tau_1+\tau_0-\eta_0]<0$, which is equivalent to $\nu_1-\tau_1<\eta_0-\tau_0$.
\item For case (iii) to hold, let $a_1=\tau_0\big[\nu_1-\tau_1+\eta_0-\tau_0(\tau_0-1) \big]<0$, which is equivalent to $\nu_1-\tau_1<\tau_0(\tau_0-1)-\eta_0$.
\end{itemize}
To conclude, if $\eta_0>\frac{\tau_0^2}{2}$, $0<\nu_0=\tau_0$, and
either $\tau_0(\tau_0-1)-\eta_0 < \nu_1-\tau_1<\eta_0-\tau_0$ or $\nu_1-\tau_1<\tau_0(\tau_0-1)-\eta_0\}$,
then Lemma \ref{lem:general conditions-SC} holds for $0<q\lesssim \mu/L$.

\textbf{Scenario 3.2:} $\bm{\eta_0=\frac{\tau_0^2}{2},\, 0<\nu_0=\tau_0}$. Then (\ref{eq:expansion of one and two}) further reduces to
\begin{equation*}
    \one = \underbrace{\tau_0\big(\nu_1-\tau_1+\tau_0(1-\frac{\tau_0}{2}) \big)}_{a_1}\sqrt{q} + O(q), \quad
    \two = \underbrace{\tau_0\big[\nu_1\tau_1-2\eta_1 +\frac{\tau_0^2}{2}(\frac{5}{2}\tau_0-2) \big]}_{b_1}\sqrt{q} + O(q).
\end{equation*}
Then $a_0=b_0=0$. We consider several special cases where $\eta_1$, $\nu_1$, and $\tau_1$ are enough to determine the convergence. Let $a_1=\tau_0\left(\nu_1-\tau_1+\tau_0(1-\frac{\tau_0}{2}) \right)<0$, which is equivalent to $\nu_1-\tau_1+\tau_0(1-\frac{\tau_0}{2})<0$. Then $N=1$.
We distinguish three cases by the sign of $b_1=\tau_0\big[\nu_1\tau_1-2\eta_1 +\frac{\tau_0^2}{2}(\frac{5}{2}\tau_0-2) \big]$.
\begin{itemize}
    \item If $b_1<0$, then $M=1$ and $b_M<0$, and hence  case (iii) holds.
    \item If $b_1>0$, then $M=N=1$ and $b_M>0$, and hence case (iia) holds. Note that the range for $q$ is $0<q\lesssim \mu^2/L^2$.
    \item If $b_1=0$, then $M\geq 2=N+1$ or $M$ does not exist. If $M$ does not exist, case (iii) is valid. If $M$ exists, then either $b_M >0$ or $b_M<0$. The former satisfies case (iib) and the latter satisfies case (iii).
\end{itemize}

Collecting the results of Scenario 1, 2 and 3 concludes the proof for Theorem \ref{thm:converge-sc-0}$^\ast$ and \ref{thm:converge-sc-1}.
\hfill $\Box$ \vspace{.1in}

\noindent\textbf{Proof of Theorems~\ref{thm:converge-sc-const}$^\ast$.}
When $\eta$, $\nu$, and $\tau$ are constants, the above analysis still holds, but we can unfold Scenario 3 without imposing strong conditions on $\eta_1$, $\nu_1$ and $\tau_1$.

We continue with Scenario 3, that is, assume $\eta\geq \frac{\tau^2}{2}$ and $0<\nu=\tau$. Then we have the Taylor expansions in finite terms:
\begin{equation}\label{eq:expansion of one and two-const}
    \begin{split}
        \one &= {\tau(\eta-(\tau-1)\tau)}\sqrt{q}-2\eta\tau(\tau-1) q + \eta\tau(\tau-1)^2 q^{\frac{3}{2}}, \\
        \two &= \tau(\tau^2-2\eta)+\eta[\tau(3\tau-2)-\eta]\sqrt{q}+\eta(\tau-1)(2\eta-\tau^2)q-\eta^2(\tau-1)^2q^{\frac{3}{2}}, \\
        \one + \sqrt{q}\two &= \tau(\tau-\eta)\sqrt{q} + \eta(\tau^2-\eta)q + \eta(\tau-1)(2\eta-\tau)q^{\frac{3}{2}}-\eta^2(\tau-1)^2q^2.
    \end{split}
\end{equation}
To proceed, we split Scenario 3 into Scenario 3.1$^\ast$ and 3.2$^\ast$ based on $\eta$.

\textbf{Scenario 3.1$^\ast$:} $\bm{\eta>\frac{\tau^2}{2},\, 0<\nu=\tau}.$ Then $b_0=\tau(\tau^2-2\eta)<0$, and $M=0$.
\begin{itemize}
    \item $\tau>2$ (hence $\tau<\frac{\tau^2}{2}<\tau(\tau-1)$).
    \begin{itemize}
        \item If $\eta>\tau(\tau-1)$, then $a_1=\tau(\eta-\tau(\tau-1))>0$ and $N=M+1=1$. Moreover, $a_1+b_{0}=\tau(\tau-\eta)<0$, and hence case (ib) in Lemma \ref{lem:general conditions-coef-SC} holds.
        \item If $\eta=\tau(\tau-1)$, then $a_1=\tau(\eta-\tau(\tau-1))=0$, and $a_2=-2\eta\tau(\tau-1)=-2\tau^2(\tau-1)^2<0$. Hence case (iii) in Lemma \ref{lem:general conditions-coef-SC} holds.
        \item If $\frac{\tau^2}{2}<\eta<\tau(\tau-1)$, then $a_1=\tau(\eta-\tau(\tau-1))<0$. Hence case (iii) in Lemma \ref{lem:general conditions-coef-SC} holds.
    \end{itemize}

    \item $\tau=2$ (hence $\tau=\frac{\tau^2}{2}=\tau(\tau-1)=2$). Then $\one=2(\eta-2)\sqrt{q}+O(q)$ and $\one + \sqrt{q}\two = 2(2-\eta)\sqrt{q} + O(q)$. By $\eta>\frac{\tau^2}{2}=2$, we have $a_1=2(\eta-2)>0$ and $a_1+b_0=2(2-\eta)<0$. Case (ib) holds.

    \item $0<\tau<2$ (hence $\tau>\frac{\tau^2}{2}>\tau(\tau-1)$). Then $\eta>\frac{\tau^2}{2}$ implies that $\eta>\tau(\tau-1)$. Hence $a_1=\tau(\eta-\tau(\tau-1))>0$ and $N=M+1=1$. Only case (ib) in Lemma \ref{lem:general conditions-coef-SC} is possible, which requires the first non-zero coefficient of $\one+\sqrt{q}\two$ in (\ref{eq:expansion of one and two-const}) to be negative, or $\one+\sqrt{q}\two\equiv 0$.
    \begin{itemize}
        \item If $\eta>\tau$, then $\one+\sqrt{q}\two=\tau(\tau-\eta)\sqrt{q}+O(q)<0$.
        \item If $\eta = \tau$, then $\one+\sqrt{q}\two=\tau^2(\tau-1)q+\tau^2(\tau-1) q^{\frac{3}{2}} -\tau^2(\tau-1)^2 q^2$. Hence for $0<\tau\leq 1$, $\one+\sqrt{q}\two\leq 0$.
        \item If $\frac{\tau^2}{2}<\eta<\tau$, then $\one+\sqrt{q}\two = \tau(\tau-\eta)\sqrt{q}+O(q)>0$, which violates case (ib).
    \end{itemize}
\end{itemize}

\textbf{Scenario 3.2$^\ast$:} $\bm{\eta=\frac{\tau^2}{2},\, 0<\nu=\tau}.$ Then (\ref{eq:expansion of one and two-const}) only involves $\tau$:
$\one=\tau^2(1-\frac{\tau}{2})\sqrt{q}-\tau^3(\tau-1) q + \frac{1}{2}\tau^3(\tau-1)^2 q^{\frac{3}{2}}$, $\two = \frac{1}{2}\tau^3(\frac{5}{2}\tau-2)\sqrt{q}-\frac{1}{4}\tau^4(\tau-1)^2 q^{\frac{3}{2}}$, and $\one+\sqrt{q}\two = \tau^2(1-\frac{\tau}{2})\sqrt{q}+O(q)$.
\begin{itemize}
    \item If $\tau>2$, then $a_N=a_1=\tau^2(1-\frac{\tau}{2})<0$, and $b_M=b_1=\frac{1}{2}\tau^3(\frac{5}{2}\tau-2)>0$. Hence case (iia) in Lemma \ref{lem:general conditions-coef-SC} holds. Note that the range of $q$ is $0<q\lesssim \mu^2/L^2$.
    \item If $\tau=2$, then $\one=-8q+4q^{\frac{3}{2}}$ and $\two=12\sqrt{q}-4q^{\frac{3}{2}}$. Hence $a_N=a_2=-8<0$ and $b_M=b_1=12>0$, with $1=M<N=2$, which contradicts case (iia) and (iib) requiring $M\ge N$. This case is invalid.
    \item If $0<\tau<2$, then $a_N=a_1=\tau^2(1-\frac{\tau}{2})>0$. Only case (ia) or (ib) in Lemma \ref{lem:general conditions-coef-SC} is possible. But then        $\one+\sqrt{q}\two = \tau^2(1-\frac{\tau}{2})\sqrt{q}+O(q)>0$, contradicting the conclusion in (ia) and (ib). This case is invalid.
\end{itemize}

Collecting the results in Scenario 3.1$^\ast$ and 3.2$^\ast$ concludes the proof for Theorem \ref{thm:converge-sc-const}$^\ast$.
\hfill $\Box$ \vspace{.1in}
%%%%%%%%%%%

Finally, we derive Theorems~\ref{thm:converge-sc-const}  and \ref{thm:converge-sc-0} from
Theorems~\ref{thm:converge-sc-const}$^\ast$  and \ref{thm:converge-sc-0}$^\ast$
by exploiting symmetrization in the case of $\nu_0\not=\tau_0$ due to Lemma \ref{lem:matching of two forms of extended-NAG-SC}
in Section \ref{sec:prf-cor-SC}.
\vspace{.1in}

\noindent\textbf{Proof of Theorems~\ref{thm:converge-sc-const}  and \ref{thm:converge-sc-0}.}
It suffices to only deal with conditions (ia) and (iia) in Theorem~\ref{thm:converge-sc-0},
which directly implies the conclusions from conditions (ia) and (iia) in Theorem~\ref{thm:converge-sc-const}.

For $\eta\sim\eta_0\geq 0$, $\nu\sim\nu_0>0$, $\tau\sim\tau_0>0$, and $\nu_0\neq \tau_0$, by Lemma \ref{lem:matching of two forms of extended-NAG-SC}, algorithm (\ref{eq:extended-NAG-SC}) can be first put into (\ref{eq:extended-NAG-SC-single}), with $R_1$, $R_2$, $R_3$ and $h_1=\frac{\zeta\eta+\nu\tau}{1+\sqrt{q}}$. Next, we keep $\eta_0$ and the remainder terms, but exchange the role of $\nu_0$ and $\tau_0$ by setting $\bar \nu_0=\tau_0$ and $\bar \tau_0=\nu_0$ and then translate (\ref{eq:extended-NAG-SC-single}) back to (\ref{eq:extended-NAG-SC}) with new parameters $\bar \eta\sim\eta_0$, $\bar \nu\sim\ \bar\nu_0= \tau_0$ and $\bar \tau\sim\bar\tau_0=\nu_0$ ($\zeta$ is also
translated to the new $\bar\zeta = 1+(1-\bar\tau)\sqrt{q}$), and a possibly nonzero $h_2= \frac{\bar\zeta\bar\eta+\bar\nu\bar\tau-(\zeta\eta+\nu\tau)}{\bar\tau(1-\bar\nu\sqrt{q})}$ in $z_0$. In other words, the original algorithm (\ref{eq:extended-NAG-SC})
can be reformulated such that the leading constants in $\nu$ and $\tau$ are exchanged and the algorithm now starts with $x_0$ and possibly $z_0\neq x_0$.

For the reformulated algorithm, the proof of Theorem \ref{thm:converge-sc-0}$^\ast$ remains valid except for the bound of $V_1$ in (\ref{eq:bounding V1-SC}) with the new $z_0$. Nevertheless, an inspection of (\ref{eq:bounding V1-SC}) reveals that $V_1\lesssim L\|x_0-x^*\|^2$ still holds because $h_2$ can be easily shown to be bounded. Hence, the desired result follows by symmetrizing the conditions (ia$^\ast$) and (iia$^\ast$) in Theorem \ref{thm:converge-sc-0}$^\ast$.
\hfill $\Box$ \vspace{.1in}

\subsection{Proof of Corollary \ref{cor:symmetrized conditions for extended-NAG-SC-single}} \label{sec:prf-cor-SC}

To facilitate interpretation of (\ref{eq:extended-NAG-SC}) and prepare for the proof of Corollary \ref{cor:symmetrized conditions for extended-NAG-SC-single},
we study the single-variable form of (\ref{eq:extended-NAG-SC}). The following lemma shows that the two forms can be transformed into each other, provided that the leading constants in $\nu$ and $\tau$ differ from each other.
The initial points need to be aligned because (\ref{eq:extended-NAG-SC}) starts from $x_0$ and $z_0$ while (\ref{eq:extended-NAG-SC-single}) starts from $x_0$ and $x_1$.

\begin{lem}\label{lem:matching of two forms of extended-NAG-SC}
    Let $\zeta = 1+(1-\tau)\sqrt{q}$ as in Lemma \ref{lem:bounding diff of kinetic-SC}.
    (i) Algorithm (\ref{eq:extended-NAG-SC}) with tuning parameters $\eta$, $\nu$ and $\tau$ under Assumption \ref{ass:regular func of sqrt-q}
    admits the single-variable form (\ref{eq:supp-extended-NAG-SC-single}), which can be expressed as
    \begin{equation}\label{eq:extended-NAG-SC-single}
        \begin{split}
            x_{k+1} &= x_k - (\nu_0\tau_0+R_1) s\nabla f(x_k)  + (1-(\nu_0+\tau_0)\sqrt{q} + R_2)(x_k-x_{k-1}) \\
            &\quad -(\eta_0+R_3)s(\nabla f(x_k)-\nabla f(x_{k-1})),
        \end{split}
    \end{equation}
with $x_0$, $x_1=x_0-h_1 s\nabla f(x_0)$, where $R_1=O(\sqrt{q})$, $R_2=O(q)$, $R_3=O(\sqrt{q})$ and $h_1=\frac{\zeta\eta+\nu\tau}{1+\sqrt{q}}$ are analytic functions of $\sqrt{q}$ around $0$.
(ii) Conversely, given any analytic functions of $\sqrt{q}$: $R_1=O(\sqrt{q})$, $R_2=O(q)$, $R_3=O(\sqrt{q})$, $h_1$, and three scalars $\eta_0\geq 0$, $\nu_0, \tau_0>0$, $\nu_0\neq \tau_0$, there exist $\eta\sim\eta_0$, $\nu\sim\nu_0$ and $\tau\sim\tau_0$ satisfying Assumption \ref{ass:regular func of sqrt-q} such that (\ref{eq:extended-NAG-SC-single}) starting from $x_0$ and $x_1=x_0-h_1s\nabla f(x_0)$ is equivalent to (\ref{eq:extended-NAG-SC}) starting from $x_0$ and $z_0=x_0+h_2\frac{\sqrt{q}}{\mu}\nabla f(x_0)$ where
    $$
    h_2 = \frac{\zeta\eta+\nu\tau-(1+\sqrt{q})h_1}{\tau(1-\nu\sqrt{q})}.
    $$
\end{lem}

\begin{prf}
First, we show that (\ref{eq:extended-NAG-SC}) admits the single-variable form (\ref{eq:supp-extended-NAG-SC-single}), i.e.,
\begin{equation*}
    \begin{split}
        x_{k+1} &= x_k - \frac{\nu(\tau+\zeta\eta\sqrt{q})}{1+\sqrt{q}} s\nabla f(x_k) + \frac{\zeta(1-\nu\sqrt{q})}{1+\sqrt{q}}(x_k-x_{k-1}) \\
        &\quad - \frac{\zeta\eta(1-\nu\sqrt{q})}{1+\sqrt{q}} s(\nabla f(x_k)-\nabla f(x_{k-1})),
    \end{split}
\end{equation*}
for $k\geq 1$ starting from $x_0$ and $x_1=x_0-\frac{\zeta\eta+\nu\tau}{1+\sqrt{q}}s\nabla f(x_0)$. The calculation for $x_1$ is straightforward and hence omitted. To show (\ref{eq:supp-extended-NAG-SC-single}), from (\ref{subeq:extended-NAG-SC-x}) and (\ref{subeq:extended-NAG-SC-y}) we have for $k\geq 0$,
\begin{equation*}
    z_{k+1} = \frac{1+\sqrt{q}}{\tau\sqrt{q}}x_{k+1} - \frac{\zeta}{\tau\sqrt{q}}y_{k+1} = \frac{1+\sqrt{q}}{\tau\sqrt{q}}x_{k+1} - \frac{\zeta}{\tau\sqrt{q}}(x_k-\eta s\nabla f(x_k)).
\end{equation*}
Substituting the above display with subscript $k+1$ and $k$ for $k\geq 1$ into (\ref{subeq:extended-NAG-SC-z}), we have
\begin{equation*}
    \begin{split}
        \frac{1+\sqrt{q}}{\tau\sqrt{q}}x_{k+1} - \frac{\zeta}{\tau\sqrt{q}}(x_k-\eta s\nabla f(x_k)) &= \nu\sqrt{q}(x_k-\frac{1}{\mu}\nabla f(x_k)) \\
        &\quad + (1-\nu\sqrt{q})\left(\frac{1+\sqrt{q}}{\tau\sqrt{q}}x_{k} - \frac{\zeta}{\tau\sqrt{q}}(x_{k-1}-\eta s\nabla f(x_{k-1}))\right).
    \end{split}
\end{equation*}
After rearrangement we obtain
\begin{equation*}
    \begin{split}
        \frac{1+\sqrt{q}}{\tau\sqrt{q}}x_{k+1} &= \left(\frac{1+\sqrt{q}}{\tau\sqrt{q}} + \frac{\zeta(1-\nu\sqrt{q})}{\tau\sqrt{q}}\right)x_k - \frac{\zeta(1-\nu\sqrt{q})}{\tau\sqrt{q}}x_{k-1} \\
        &\quad - \frac{\zeta\eta+\nu\tau}{\tau\sqrt{q}}s\nabla f(x_k) + \frac{\zeta\eta(1-\nu\sqrt{q})}{\tau\sqrt{q}}s\nabla f(x_{k-1}).
    \end{split}
\end{equation*}
Solving for $x_{k+1}$ yields (\ref{eq:supp-extended-NAG-SC-single}). Expanding the coefficients in (\ref{eq:supp-extended-NAG-SC-single}) yields (\ref{eq:extended-NAG-SC-single}).

Second, we prove the reverse statement. Given $R_1=O(\sqrt{q})$, $R_2=O(q)$, $R_3=O(\sqrt{q})$ and the leading constants $\eta_0$, $\nu_0$ and $\tau_0$, we determine $\eta$, $\nu$ and $\tau$ by solving the following equations
\begin{equation*}
    \begin{split}
        &\nu(\tau+\zeta\eta\sqrt{q}) = (\nu_0\tau_0+R_1)(1+\sqrt{q}) = g_1 \sim \nu_0\tau_0 + O(\sqrt{q}) , \\
        &\zeta(1-\nu\sqrt{q}) = (1-(\nu_0+\tau_0)\sqrt{q} + R_2)(1+\sqrt{q}) = g_2 \sim 1- (\nu_0+\tau_0-1)\sqrt{q} + O(q), \\
        &\zeta\eta(1-\nu\sqrt{q}) = (\eta_0+ R_3)(1+\sqrt{q}).
    \end{split}
\end{equation*}
Note that $g_1$ and $g_2$ above are known. Solving the equations we obtain
$$
\eta = \frac{\eta_0 + R_3}{1-(\nu_0+\tau_0)\sqrt{q} + R_2}\sim \eta_0,
$$
and $\nu$ depending on $\tau$ as
$$
\nu = \frac{g_1}{\tau + \zeta\eta\sqrt{q}} \sim \frac{\nu_0\tau_0}{\tau_0}= \nu_0,
$$
and $\tau$ as a root for the quadratic equation
$$
\alpha_2\cdot \tau^2 + \alpha_1\cdot \tau + \alpha_0=0,
$$
where
\begin{equation*}
    \begin{split}
        &\alpha_2 = -\sqrt{q}(1-\eta q), \\
        &\alpha_1 = 1 - g_2 + \sqrt{q} + (g_1-2\eta+\eta g_2)q - 2\eta q^{\frac{3}{2}} = \sqrt{q}(\nu_0+\tau_0+O(\sqrt{q})), \\
        &\alpha_0 = (\eta-g_1-\eta g_2)\sqrt{q} + (2\eta-g_1-\eta g_2)q + \eta q^{\frac{3}{2}} = \sqrt{q}(-\nu_0\tau_0 + O(\sqrt{q})).
    \end{split}
\end{equation*}
The discriminant is $\Delta = \alpha_1^2-4\alpha_2\cdot \alpha_0=q[(\nu_0-\tau_0)^2+O(\sqrt{q})]$.
For $\nu_0\neq \tau_0>0$, the root
\begin{equation*}
    \begin{split}
        \tau &= \frac{-\alpha_1 \pm \sqrt{\Delta}}{2\cdot \alpha_2}=\frac{\alpha_1/\sqrt{q} \pm \sqrt{\Delta/q}}{-2\cdot \alpha_2/\sqrt{q}}\\
        &=\frac{\nu_0+\tau_0+O(\sqrt{q})\pm \sqrt{(\nu_0-\tau_0)^2 +O(\sqrt{q})}}{2(1-\eta q)}\sim \frac{\nu_0+\tau_0\pm|\nu_0-\tau_0|}{2}
    \end{split}
\end{equation*}
is well-defined for small $q$. The sign is determined by the sign of $\nu_0-\tau_0$ to make $\tau\sim \tau_0$.
An inspection of the expressions suggests that $\eta$, $\nu$ and $\tau$ are all analytic functions of $\sqrt{q}$.

From the above calculation, the updating formula for $x_k$ is matched between (\ref{eq:extended-NAG-SC}) and (\ref{eq:extended-NAG-SC-single}),
starting from $x_2$. The initial points can be aligned by picking suitable $z_0$ such that $x_1$ from (\ref{eq:extended-NAG-SC})
is exactly $x_0-h_1 s\nabla f(x_0)$. The calculation is straightforward and hence omitted.
\end{prf}

Interestingly, the leading coefficients in the single-variable form (\ref{eq:extended-NAG-SC-single}) for (\ref{eq:extended-NAG-SC}) are symmetric with regard to $\nu_0$ and $\tau_0$. This suggests that the convergence properties of (\ref{eq:extended-NAG-SC}) should also be symmetric in $\nu_0$ and $\tau_0$.
Indeed, such symmetrization in the case of $\nu_0 \not=\tau_0$  is exploited in our derivation of Theorems~\ref{thm:converge-sc-const}--\ref{thm:converge-sc-0}
from  Theorems~\ref{thm:converge-sc-const}$^\ast$--\ref{thm:converge-sc-0}$^\ast$ in Section \ref{sec:prf-thm-SC}.

Corollary \ref{cor:symmetrized conditions for extended-NAG-SC-single} is derived by translating
the conclusions from Theorem~\ref{thm:converge-sc-0} in the case of $\nu_0 \not= \tau_0$
in terms of the leading coefficients in a single-variable form. \vspace{.1in}

\noindent\textbf{Proof of Corollary \ref{cor:symmetrized conditions for extended-NAG-SC-single}.}
   Let $c_1 = \nu_0 + \tau_0$, $c_0=\nu_0\tau_0$, $c_2\sqrt{c_0}-c_0/2=\eta_0$ and solve for $\eta_0$, $\nu_0$ and $\tau_0$. The conditions $c_1^2 > 4c_0$ and $c_0,c_1>0$ ensure that $\nu_0$ and $\tau_0$ exist, satisfying $\nu_0,\tau_0>0$ and $\nu_0\neq \tau_0$.
   The conditions on $c_2$ can be directly translated into the conditions on $\eta_0$. The desired result
   follows by applying Lemma \ref{lem:matching of two forms of extended-NAG-SC} and Theorem \ref{thm:converge-sc-0}.
\hfill $\Box$ \vspace{.1in}

\section{Technical details in Section \ref{sec:acc-c}}\label{sec:tech for section c}

\subsection{Proof of three-variable form (\ref{eq:extended-NAG-C-xyz})}

By comparing (\ref{eq:extended-NAG-C-xyz}) with (\ref{eq:extended-NAG-C}), the update of $y_{k+1}$ is the same. Rearranging $z_{k+1}=\a_{k+1}x_{k+1}+(1-\a_{k+1})y_{k+1}$, we obtain the update of $x_{k+1}$ in (\ref{subeq:extended-NAG-C-x}). It suffices to prove the update of $z_{k+1}$ in (\ref{subeq:extended-NAG-C-z}). With (\ref{eq:extended-NAG-C}) in place, we have
\begin{align*}
    z_{k+1} &= \a_{k+1}x_{k+1}+(1-\a_{k+1})y_{k+1} \\
    &= \a_{k+1}[x_k - \gamma_k s \nabla f(x_k) + \sigma_{k+1}(y_{k+1}-y_k)] + (1-\a_{k+1})y_{k+1} \\
    &= \a_{k+1}(x_k - \gamma_k s \nabla f(x_k)) + (\a_k-1)(y_{k+1}-y_k) + (1-\a_{k+1})y_{k+1} \\
    &= \a_{k+1}(x_k - \gamma_k s \nabla f(x_k)) + (\a_k-\a_{k+1})y_{k+1} + (1-\a_k)y_k \\
    &= \a_{k+1}(x_k - \gamma_k s \nabla f(x_k)) + (\a_k-\a_{k+1})(x_k - \beta_k s \nabla f(x_k)) + (1-\a_k)y_k \\
    &= \a_k x_k + (1-\a_k)y_k - (\gamma_k\a_{k+1} + \beta_k(\a_k-\a_{k+1}) )s\nabla f(x_k) \\
    &= z_{k} - \widetilde \a_{k} s\nabla f(x_k).
\end{align*}
The proof is completed.

\subsection{Proofs of Lemmas~\ref{lem:bounding diff of kinetic-C} and \ref{lem:general conditions-C} in Section \ref{sec:outline-C}}

\textbf{Proof of Lemma~\ref{lem:bounding diff of kinetic-C}.}
Using (\ref{subeq:extended-NAG-C-z}), we have
\begin{equation*}
    \begin{split}
        \frac{1}{2}\|z_{k+1}-x^*\|^2 - \frac{1}{2}\|z_{k}-x^*\|^2 &= \frac{1}{2}\|z_k-x^*-\widetilde \a_k s\nabla f(x_k) \|^2 - \frac{1}{2}\|z_{k}-x^*\|^2 \\
        &= -\widetilde \a_k s\langle z_k-x^*,\nabla f(x_k) \rangle + \frac{\widetilde \a_k^2s^2}{2}\|\nabla f(x_k)\|^2.
    \end{split}
\end{equation*}
Substituting $z_k = x_k + (\a_k-1)(x_k-y_k)$ into the above display, we have
\begin{equation}\label{eq:bounding diff-1}
    \begin{split}
        &\frac{1}{2}\|z_{k+1}-x^*\|^2 - \frac{1}{2}\|z_{k}-x^*\|^2\\
        &= -\widetilde\a_k s\langle x_k-x^*,\nabla f(x_k) \rangle - \widetilde\a_k (\a_k-1) s\langle x_k-y_k,\nabla f(x_k) \rangle + \frac{\widetilde \a_k^2s^2}{2}\|\nabla f(x_k)\|^2 \\
        &\leq -\widetilde\a_k s(f(x_k)-f^*) - \widetilde\a_k(\a_k-1)s(f(x_k)-f(y_k)) +  \frac{\widetilde \a_k^2s^2}{2}\|\nabla f(x_k)\|^2 \\
        &=-\a_k\widetilde\a_ks(f(x_k)-f^*) + \widetilde\a_k(\a_k-1)s(f(y_k)-f^*) + \frac{\widetilde \a_k^2s^2}{2}\|\nabla f(x_k)\|^2,
    \end{split}
\end{equation}
where the inequality holds because $\langle x_k-x^*,\nabla f(x_k) \rangle\geq f(x_k)-f^*$ and $\langle x_k-y_k,\nabla f(x_k) \rangle\geq f(x_k)-f(y_k)$ by the convexity of $f$ and the assumption that $\widetilde\a_k\geq 0$ and $\a_k\geq 1$. When $k\geq 1$, by the $L$-smoothness of $f$ and (\ref{subeq:extended-NAG-C-y}), we have
\begin{equation}\label{eq:bounding diff-2}
    \begin{split}
        f(y_k)-f^* &\leq f(x_{k-1}) -f^* + \langle \nabla f(x_{k-1}), y_k-x_{k-1} \rangle + \frac{L}{2}\|y_k-x_{k-1}\|^2 \\
    &= f(x_{k-1}) -f^* - (2-\b_{k-1}Ls)\frac{\b_{k-1}s}{2}\|\nabla f(x_{k-1})\|^2.
    \end{split}
\end{equation}
Combining (\ref{eq:bounding diff-1}) and (\ref{eq:bounding diff-2}) yields (\ref{eq:bounding diff of kinetic-C}), which completes the proof.
\hfill $\Box$ \vspace{.1in}

\noindent \textbf{Proof of Lemma~\ref{lem:general conditions-C}.}
    By Lemma \ref{lem:bounding diff of kinetic-C}, for any $k\geq 1$ such that $\alpha_k\geq 1$, $\widetilde\alpha_k\geq 0$, we have
    \begin{equation}\label{eq:bounding diff of Lyapunov}
    \begin{split}
        V_{k+1}-V_k \leq -\frac{1}{2}(\w_k-\w_{k+1})\|z_k-x^*\|^2 - \underbrace{(\w_k\a_{k-1}\widetilde\a_{k-1}-\w_{k+1}\widetilde\a_k(\a_k-1))}_{\one}s(f(x_{k-1})-f^*) \\- \underbrace{\left(\w_{k+1}\widetilde\a_k(\a_k-1)\b_{k-1}(2-\b_{k-1}Ls)-\w_k\widetilde\a_{k-1}^2\right)}_{\two}\frac{s^2}{2}\|\nabla f(x_{k-1})\|^2.
    \end{split}
\end{equation}
If $\w_k\geq\w_{k+1}$, then $-(\w_k-\w_{k+1})\|z_k-x^*\|^2/2\le 0$. The rest is straightforward.
\hfill $\Box$ \vspace{.1in}

\subsection{Proof of Theorem \ref{thm:converge-c}}

To prepare for the proof for Theorem \ref{thm:converge-c}, we provide a simple bound which will be used in the last step of the proof.
\begin{lem}\label{lem:simple bound-C}
    Let $f\in\mathcal{F}^1_L $. When $0<s\leq C_0/L$ for some constant $C_0>0$, the iterates of (\ref{eq:extended-NAG-C-xyz}) satisfy that
     for any fixed $K$, there exists a constant $C$ such that $s(f(x_K)-f^*)$, $s^2\|\nabla f(x_K)\|^2$ and $\|z_{K+1}-x^*\|^2$ are upper-bounded by $C\|x_0-x^*\|^2$. The constant $C$ depends only on $C_0$, $K$ and the algorithm parameters $\{\a_k\}$, $\{\b_k\}$ and $\{\gamma_k\}$.
\end{lem}

\begin{prf}
For notational simplicity, we assume $C_0=1$ so that $0< Ls\leq 1$. First, by convexity of $f$, we have
$$
s(f(x_K)-f^*) \leq \frac{Ls}{2}\|x_K-x^*\|^2 \leq \frac{1}{2}\|x_K-x^*\|^2.
$$
Second, by the $L$-smoothness, we have
$$
s^2\|\nabla f(x_K)\|^2=s^2\|\nabla f(x_K)-\nabla f(x^*)\|^2\leq (Ls)^2 \|x_K-x^*\|^2\leq \|x_K-x^*\|^2.
$$
Third, by (\ref{subeq:extended-NAG-C-z}) and the Cauchy--Schwartz inequality, we have
\begin{equation*}
    \begin{split}
        \frac{1}{2}\|z_{K+1}-x^*\|^2 = \frac{1}{2} \|z_K-x^*-\widetilde\a_K s\nabla f(x_K)\|^2 \leq \|z_K-x^*\|^2 + \widetilde\a_K^2 s^2\|\nabla f(x_K)\|^2,
    \end{split}
\end{equation*}
which together with the preceding bound on $s^2\|\nabla f(x_K)\|^2$ yields
\begin{equation*}
    \frac{1}{2}\|z_{K+1}-x^*\|^2 \leq \|z_K-x^*\|^2 + \widetilde\a_K^2 \|x_K-x^*\|^2.
\end{equation*}
Because $z_0=x_0$, it suffices to bound $\|x_k-x^*\|^2$ by $\|x_0-x^*\|^2$ up to a constant for general $k\ge 1$.
The fact that (\ref{eq:extended-NAG-C-xyz}) is a first-order method implies that for each $k \ge 1$, there exist scalars $\{c_{k,i}\}$ (depending only on the algorithm parameters) for $i=0,\ldots,{k-1}$ such that
$$
x_k = x_{k-1} + \sum_{i=0}^{k-1} c_{k,i}s\nabla f(x_i).
$$
Then by the Cauchy--Schwartz inequality,
\begin{equation*}
    \begin{split}
        \|x_k-x^*\|^2 &= \|x_{k-1} - x^* + \sum_{i=0}^{k-1} c_{k,i}s\nabla f(x_i)\|^2 \leq k \left(\|x_{k-1} - x^*\|^2 + \sum_{i=0}^{k-1}c_{k,i}^2 s^2\|\nabla f(x_i)\|^2\right) \\
        &\leq k \left(\|x_{k-1} - x^*\|^2 + \sum_{i=0}^{k-1}c_{k,i}^2 \|x_i-x^*\|^2\right).
    \end{split}
\end{equation*}
Because $K$ is fixed, the proof is completed by applying the preceding bound for $1\le k \le K$.
\end{prf}

With Lemma \ref{lem:general conditions-C} and Lemma \ref{lem:simple bound-C} in place, we are ready to prove Theorem \ref{thm:converge-c}. \vspace{.1in}

\noindent\textbf{Proof of Theorem \ref{thm:converge-c}.} By the definition of $\widetilde\a_k$ and condition (i)-(ii), we have
\begin{equation*}
    \lim_{k\to \infty} \frac{\widetilde\a_k}{\a_k} =
    \lim_{k\to \infty} \b_k + (\g_k-\b_k)\frac{\a_{k+1}}{\a_k} = \g >0 .
\end{equation*}
Then $\a_k=\Omega(k)$ implies that $\widetilde\a_k=\Omega(k)$. By Lemma \ref{lem:general conditions-C}, the key is to bound $\one$ and $\two$ from below for all $k$ sufficiently large. The beginning $V_k$ can be dealt with by the simple bound in Lemma \ref{lem:simple bound-C} so they do not affect the convergence. To proceed, we consider two choices of $\{\w_k\}$ depending on the monotonicity of $\{\widetilde\a_k/\a_k\}$.
Note that $\{\w_k\}$ needs to be non-increasing, to ensure that the first term $-(\w_k-\w_{k+1})\|z_k-x^*\|^2/2$ on the right-hand-side of (\ref{eq:bounding diff of Lyapunov}) is always non-positive.

\textbf{Choice 1:} $\bm{\w_k\equiv 1.}$ We pick $w_k$ simply as constant $1$ when $\{\widetilde\a_k/\a_k\}$ is non-increasing.
As for $\one$, by condition (ii) we have
\begin{equation*}
    \begin{split}
        \one &= \a_{k-1}\widetilde\a_{k-1}-\widetilde\a_k(\a_k-1)= \a_{k-1}^2\left( \frac{\widetilde\a_{k-1}}{\a_{k-1}} - \frac{\a_k(\a_k-1)}{\a_{k-1}^2}\cdot\frac{\widetilde\a_k}{\a_k} \right) \\
        &\geq \a_{k-1}^2 \left( \frac{\widetilde\a_{k-1}}{\a_{k-1}} - \frac{\widetilde\a_k}{\a_k} \right) \geq 0.
    \end{split}
\end{equation*}
As for $\two$, we have
\begin{equation*}
    \begin{split}
        \two &= \widetilde\a_k(\a_k-1)\b_{k-1}(2-\b_{k-1}Ls)-\widetilde\a_{k-1}^2 \\
        &= \a_k(\a_k-1)\left(\frac{\widetilde\a_k}{\a_k}\b_{k-1}(2-\b_{k-1}Ls)-\frac{\widetilde\a_{k-1}^2}{\a_k(\a_{k}-1)} \right).
    \end{split}
\end{equation*}
For $0<s\leq (2-\g/\b)/(2\b L)$ with $\b>\gamma/2>0$,
\begin{equation*}
    \begin{split}
        \two \geq \a_k(\a_k-1)\left[\frac{\widetilde\a_k}{\a_k}\b_{k-1}\left(2-\frac{\b_{k-1}}{2\b}(2-\frac{\g}{\b})\right)-\frac{\widetilde\a_{k-1}^2}{\a_k(\a_{k}-1)} \right],
    \end{split}
\end{equation*}
where the limit of the term in the square brackets is
$$
\lim_{k\to \infty} \frac{\widetilde\a_k}{\a_k}\b_{k-1}\left(2-\frac{\b_{k-1}}{2\b}(2-\frac{\g}{\b})\right)-\frac{\widetilde\a_{k-1}^2}{\a_k(\a_{k}-1)}=\g(\b - \frac{\g}{2})>0.
$$
Combining the preceding two displays and $\alpha_k = \Omega(k)$ shows that
there exist constants $K$ and $C>0$ depending only on algorithm parameters such that for $0<s \le (2-\g/\b)/(2\b L)$ and $k\ge K$,
$$
\two \ge \frac{\g}{2}(\b - \frac{\g}{2})\a_k(\a_k-1)\geq Ck^2.
$$

\textbf{Choice 2:} $\bm{\w_{k+1} = \a_k/\widetilde\a_k.}$ When $\{\widetilde\a_k/\a_k\}$ is non-decreasing, $\{\w_k\}$ is non-increasing. In this case we have $\one = \a_{k-1}^2-\a_k(\a_k-1)\geq 0$. For $0<s \le (2-\g/\b)/(2\b L)$,
\begin{equation*}
    \begin{split}
        \two &= \a_k(\a_k-1)\b_{k-1}(2-\b_{k-1}Ls)-\a_{k-1}\widetilde\a_{k-1} \\
        &= \a_k(\a_k-1)\left(\b_{k-1}(2-\beta_{k-1}Ls) - \frac{\a_{k-1}^2}{\a_k(\a_k-1)}\cdot\frac{\widetilde\a_{k-1}}{\a_{k-1}}\right) \\
        & \geq \a_k(\a_k-1)\left[\b_{k-1}\left(2-\frac{\beta_{k-1}}{2\b}(2-\frac{\g}{\b})\right) - \frac{\a_{k-1}^2}{\a_k(\a_k-1)}\cdot\frac{\widetilde\a_{k-1}}{\a_{k-1}}\right],
    \end{split}
\end{equation*}
where the limit of the term in the square brackets is
$$
\lim_{k\to\infty}\b_{k-1}\left(2-\frac{\beta_{k-1}}{2\b}(2-\frac{\g}{\b})\right) - \frac{\a_{k-1}^2}{\a_k(\a_k-1)}\cdot\frac{\widetilde\a_{k-1}}{\a_{k-1}} = \b - \frac{\g}{2}>0.
$$
Similarly as in the first choice, there exist constants $K$ and $C>0$ depending only on algorithm parameters such that for $0<s \le (2-\g/\b)/(2\b L)$ and $k\ge K$,
$$
\two \ge (\b - \frac{\g}{2})\a_k(\a_k-1)\ge  Ck^2 .
$$

Combining the two choices above, we see that when $\{\widetilde\a_k/\a_k\}$ is either non-increasing or non-decreasing in $k$, there exist constants $K$ and $C>0$ such that for $k\geq K$, we have $\one\geq 0$ and
$\inf_{0<s\leq C_0/L}\two \ge C k^2$,
where $C_0=(2-\g/\b)/(2\b)$.
By Lemma \ref{lem:general conditions-C},
$$V_{k+1}-V_k\le - \frac{C}{2} k^2 s^2\|\nabla f(x_{k-1})\|^2 \le 0.
$$
To complete the proof of (\ref{eq:optimal bound-c}) for the objective gap, it suffices to bound $V_k$ from below.

By convexity of $f$, we have $\|\nabla f(x_k)\|^2\leq (2L) (f(x_k)-f^*)$ which together with
the definition of $V_{k+1}$ in (\ref{eq:lyapunov-C}) yields that for $k\geq K$ and $0<s\leq C_0/L$,
\begin{equation*}
    \begin{split}
        V_K \geq V_{k+1} \geq  \w_{k+1}(\a_k\widetilde\a_k-\widetilde\a_k^2Ls)s(f(x_k)-f^*)= \w_{k+1}\a_k\widetilde\a_k\left(1-\frac{\widetilde\a_k}{\a_k}Ls \right)s(f(x_k)-f^*).
    \end{split}
\end{equation*}
Because $\lim_k\w_k$ is $1$ or $1/\g>0$ in the Choice 1 or 2 above, and $\lim_k \widetilde\a_k/\a_k=\gamma>0$, we reset $K$ large enough
such that the above is further bounded from below as
$$
V_K \geq V_{k+1} \geq  \frac{1}{2}(\frac{1}{\g} \vee 1) \a_k\widetilde\a_k(1-2\g Ls)s(f(x_k)-f^*).
$$
By resetting $C_0=\frac{2-\g/\b}{2\b} \wedge \frac{1}{4\g}$, we have that for $k\geq K$ and $0<s\leq C_0/L$,
$$
V_K\geq V_{k+1}\geq  \frac{1}{4}(\frac{1}{\g} \vee 1)  \a_k\widetilde\a_k s(f(x_k)-f^*).
$$
Using Lemma \ref{lem:simple bound-C}, we have $f(x_k)-f^*=O(V_K/(\a_k\widetilde\a_k s))=O(V_K/sk^2)=O(\|x_0-x^*\|^2/(sk^2))$,
which is the optimal rate (\ref{eq:optimal bound-c}).

Finally, we complete the proof of (\ref{eq:inverse cubic rate}) for the squared gradient norm.
For $k\geq K$ and $0<s\leq C_0/L$, using $V_{k+1}-V_k\le -\frac{C}{2}
k^2s^2\|\nabla f(x_{k-1})\|^2$, we have
\begin{equation*}
     \frac{C}{2}\left(\sum_{i=K}^k i^2\right)s^2\min_{0\leq i \leq k-1} \|\nabla f(x_i)\|^2 \leq \frac{C}{2}\sum_{i=K}^k i^2s^2\|\nabla f(x_{i-1})\|^2 \le \sum_{i=K}^k (V_i-V_{i+1})=V_K-V_{k+1} \leq V_K,
\end{equation*}
which together with Lemma \ref{lem:simple bound-C} implies that for $k\geq 2K$,
\begin{equation*}
    \min_{0\leq i \leq k-1}\|\nabla f(x_i)\|^2 \leq \frac{2V_K}{Cs^2\sum_{i=K}^k i^2} \le \frac{4V_K}{Cs^2\sum_{i=1}^k i^2}=\frac{24V_K}{Cs^2k(k+1)(2k+1)} = O\left(\frac{\|x_0-x^*\|^2}{s^2k^3} \right),
\end{equation*}
which is the desired bound (\ref{eq:inverse cubic rate}).
\hfill $\Box$ \vspace{.1in}

\subsection{Proof of Lemma \ref{lem:counterexample for alphak}}
First, we show $\a_k=\Omega(k)$. When $k$ is odd, $\a_k=(1+\sqrt{1+4\a^2_{k-1}})/2>(1+2\a_{k-1})/2=(1+2(k+r-1)/r)/2=(k+r-1)/r + 1/2$. So for any $r>0$, $\a_k=\Omega(k)$.

Second, we show (\ref{eq:recursive condition}), $\a_k(\a_k-1) \le \a_{k-1}^2$ for $k\ge 1$. When $k$ is odd, by construction (\ref{eq:recursive condition}) holds. When $k\geq 2$ and $k$ is even, $k-1$ is odd and $k-2$ is even.
Then $\a_k=(k+r)/r$ and $\a_{k-1}=(1+\sqrt{1+4\a^2_{k-2}})/2=(1+\sqrt{r^2+4(k+r-2)^2 }/r)/2$. Simple algebra yields
\begin{equation}\label{eq:diff of recursive}
    \a_k(\a_k-1)-\a_{k-1}^2= \frac{4-r}{r^2}k - \frac{\sqrt{(k+r-2)^2+r^2/4}}{r} - \frac{1}{2} - \left(\frac{r-2}{r}\right)^2.
\end{equation}
When $r\geq 4$, (\ref{eq:diff of recursive}) is negative for all $k$. For $0<r<4$, we further rearrange (\ref{eq:diff of recursive}) as
\begin{align*}\label{eq:diff of recursive-2}
   & \quad   \frac{(4-r)^2k^2-r^2\left((k+r-2)^2 + \frac{r^2}{4} \right)}{r^2\left((4-r)k + r\sqrt{(k+r-2)^2+ \frac{r^2}{4}}  \right)} - \frac{1}{2} - \left(\frac{r-2}{r}\right)^2  \\
    & = -\frac{8(r-2)k^2 + 2r^2(r-2)k + r^2\left((r-2)^2+\frac{r^2}{4}\right)}{r^2\left((4-r)k + r\sqrt{(k+r-2)^2+ \frac{r^2}{4}}  \right)} - \frac{1}{2} - \left(\frac{r-2}{r}\right)^2,
\end{align*}
which remains negative for $2\leq r<4$ but blows up to $\infty$ as $k\rightarrow \infty$ if $0<r<2$.

Third, we calculate $\lim_k k(1-\sigma_{k+1})$ along $\{k'\}=\{2k\}$ and $\{k''\}=\{2k+1\}$. By the definition of $\sigma_{k+1}$ we have
\begin{equation}\label{eq:calculate-limit}
    \lim_{k\to\infty} k(1-\sigma_{k+1}) = \lim_{k\to\infty} k\left(1-\frac{\a_k-1}{\a_{k+1}} \right)=\lim_{k\to\infty} \frac{k}{\a_{k+1}}\left(\a_{k+1}-\a_k+1\right).
\end{equation}
Along $\{k'\}=\{2k\}$,
\begin{equation*}
    \lim_{k\to\infty} \frac{2k}{\a_{2k+1}} = \lim_{k\to\infty} \frac{4k}{1+\sqrt{1+4\a^2_{2k}}} = \lim_{k\to\infty} \frac{4k}{2\a_{2k}}= \lim_{k\to\infty} \frac{4k}{2(2k+r)/r} = r,
\end{equation*}
and
\begin{equation*}
    \begin{split}
        \lim_{k\to\infty} \a_{2k+1}-\a_{2k}&= \lim_{k\to\infty} \frac{1+\sqrt{1+4\a_{2k}^2}}{2} - \a_{2k} \\
        &= \lim_{k\to\infty} \frac{1}{2}\left(1 + 2\a_{2k}\sqrt{1+\frac{1}{4\a^2_{2k}}} \right)-\a_{2k} \\
        &= \lim_{k\to\infty} \frac{1}{2} + \a_{2k}\left(1+O\left(\frac{1}{\a^2_{2k}}\right)\right) - \a_{2k} \\
        &= \frac{1}{2}.
    \end{split}
\end{equation*}
Using (\ref{eq:calculate-limit}) we have 
\begin{align}
\lim_{k'\to\infty} k'(1-\sigma_{k'+1})=r(1+1/2)=3r/2. \label{eq:sigma-limit-a}
\end{align}
Along $\{k''\}=\{2k+1\}$,
$$
\lim_{k\to\infty} \frac{2k+1}{\a_{2k+2}} = \lim_{k\to\infty} \frac{2k+1}{(2k+2+r)/r} = r,
$$
and
\begin{equation*}
    \begin{split}
        \lim_{k\to\infty} \a_{2k+2} - \a_{2k+1} &= \lim_{k\to\infty} \a_{2k+2} - \frac{1+\sqrt{1+4\a^2_{2k}}}{2} \\
        &=\lim_{k\to\infty} \a_{2k+2} - \frac{1}{2}\left(1+2\a_{2k}\sqrt{1+\frac{1}{4\a^2_{2k}}} \right) \\
        &=\lim_{k\to\infty} \a_{2k+2}-\frac{1}{2}-\a_{2k}\left(1+O\left(\frac{1}{\a^2_{2k}}\right) \right) \\
        &= \lim_{k\to\infty} \frac{2k+2+r}{r}- \frac{2k+r}{r}- \frac{1}{2}\\
        &= \frac{2}{r} - \frac{1}{2}.
    \end{split}
\end{equation*}
Using (\ref{eq:calculate-limit}) we have
\begin{align}
\lim_{k''\to\infty} k''(1-\sigma_{k''+1}) = r\left(\frac{2}{r} - \frac{1}{2} + 1 \right)=2 + \frac{r}{2}. \label{eq:sigma-limit-b}
\end{align}

Lastly, combining the limits (\ref{eq:sigma-limit-a}) and (\ref{eq:sigma-limit-b}) gives $\lim_k \sigma_{k+1}=1$, which implies that $\lim_{k} \alpha_{k+1}/\alpha_k=1$. Therefore, condition (ii) holds. The proof of Lemma \ref{lem:counterexample for alphak} is completed.

\section{Technical details in Section \ref{sec:ODE-interpretation}} \label{sec:tech for section ODE}

\subsection{Convergence of ODEs (\ref{eq:gradient flow}) and (\ref{eq:low-res-ODE-NAG-SC})}

In the general convex setting where $f\in \mathcal{F}^1 $, the convergence rates of gradient flow (\ref{eq:gradient flow}) and the ODE (\ref{eq:low-res-ODE-NAG-C}) for NAG-C are proved in \cite{su2016differential}. For completeness, we present the proofs for the convergence rates stated in the strongly convex setting where $f\in\mathcal{S}^1_{\mu} $.

For gradient flow (\ref{eq:gradient flow}), consider the Lyapunov function $V_t=f(X_t)-f^*$, i.e., the potential gap itself. Then by (\ref{eq:gradient flow}), we have $\dot V_t=\langle \nabla f(X_t),\dot X_t \rangle = -\|\nabla f(X_t)\|^2$. By strong convexity, $f(X_t)-f^*\leq \frac{1}{2\mu}\|\nabla f(X_t)\|^2$. Therefore, $\dot V_t\leq -2\mu (f(X_t)-f^*)=-2\mu V_t$, which implies that $f(X_t)-f^*=V_t \leq V_0 \me^{-2\mu t}=\me^{-2\mu t}(f(x_0)-f^*)$.

For the ODE (\ref{eq:low-res-ODE-NAG-SC}) corresponding to NAG-SC (\ref{eq:NAG-SC-single}), consider the Lyapunov function $V_t=f(X_t)-f^*+\frac{1}{2}\|\dot X_t+\sqrt{\mu}(X_t-x^*)\|^2$. Then $V_0=f(x_0)-f^* + \frac{\mu}{2}\|x_0-x^*\|^2\leq 2(f(x_0)-f^*)$ by strong convexity. Using the ODE (\ref{eq:low-res-ODE-NAG-SC}), we have
\begin{equation*}
    \begin{split}
        \frac{\dif V_t}{\dif t}&= \langle \nabla f(X_t),\dot X_t \rangle + \langle \dot X_t + \sqrt{\mu}(X_t-x^*), \ddot X_t + \sqrt{\mu}\dot X_t \rangle \\
        &=\langle \nabla f(X_t),\dot X_t \rangle - \langle \dot X_t + \sqrt{\mu}(X_t-x^*), \sqrt{\mu}\dot X_t + \nabla f(X_t) \rangle \\
        &=-\sqrt{\mu}\|\dot X_t\|^2 - \sqrt{\mu}\langle \dot X_t, \sqrt{\mu}(X_t-x^*)\rangle - \sqrt{\mu}\langle X_t-x^*,\nabla f(X_t) \rangle,
    \end{split}
\end{equation*}
which together with the inequality $\langle X_t-x^*,\nabla f(X_t) \rangle \geq f(X_t)-f^*+\frac{\mu}{2}\|X_t-x^*\|^2$ by $\mu$-strong convexity suggests that
\begin{equation*}
    \frac{\dif V_t}{\dif t} \leq - \frac{\sqrt{\mu}}{2}\|\dot X_t\|^2-\sqrt{\mu}\left(f(X_t)-f^*+\frac{1}{2}\|\dot X_t + \sqrt{\mu}(X_t-x^*)\|^2\right)\leq -\sqrt{\mu} V_t.
\end{equation*}
Hence, $f(X_t)-f^*\leq V_t\leq V_0 \me^{-\sqrt{\mu}t}\leq 2\me^{-\sqrt{\mu}t}(f(x_0)-f^*)$.

\subsection{Proof of Proposition \ref{pro:low-res-ODE-SC}}

Consider the Lyapunov function $V_t=c_0(f(X_t)-f^*)+\frac{1}{2}\|\dot X_t + \lambda \sqrt{\mu}(X_t-x^*)\|^2$, where $\lambda$ is to be chosen later. Using the ODE (\ref{eq:low-res-ODE-SC}), we have
\begin{equation*}
    \begin{split}
        \frac{\dif V_t}{\dif t} &= c_0 \langle \nabla f(X_t), \dot X_t \rangle + \langle \dot X_t + \lambda \sqrt{\mu}(X_t-x^*), \ddot X_t + \lambda\sqrt{\mu}\dot X_t \rangle \\
        &= c_0 \langle \nabla f(X_t), \dot X_t \rangle + \langle \dot X_t + \lambda \sqrt{\mu}(X_t-x^*), -(c_1-\lambda)\sqrt{\mu}\dot X_t - c_0\nabla f(X_t) \rangle \\
        &= -(c_1-\lambda)\sqrt{\mu}\|\dot X_t\|^2 - (c_1-\lambda)\sqrt{\mu}\langle \dot X_t,\lambda\sqrt{\mu}(X_t-x^*) \rangle - c_0\lambda\sqrt{\mu}\langle X_t-x^*,\nabla f(X_t) \rangle \\
        &= -\frac{c_1-\lambda}{2}\sqrt{\mu}\|\dot X_t\|^2 - \frac{c_1-\lambda}{2}\sqrt{\mu}\|\dot X_t + \lambda \sqrt{\mu}(X_t-x^*)\|^2 \\
        &\quad + \frac{c_1-\lambda}{2}\lambda^2\mu^{\frac{3}{2}}\|X_t-x^*\|^2 - c_0\lambda\sqrt{\mu}\langle X_t-x^*,\nabla f(X_t) \rangle \\
        &= -2C\sqrt{\mu}V_t - \frac{c_1-\lambda}{2}\sqrt{\mu}\|\dot X_t\|^2 - \left(\frac{c_1-\lambda}{2}-C\right)\sqrt{\mu}\|\dot X_t + \lambda \sqrt{\mu}(X_t-x^*)\|^2 \\
        &\quad + 2Cc_0\sqrt{\mu}(f(X_t)-f^*) + \frac{c_1-\lambda}{2}\lambda^2\mu^{\frac{3}{2}}\|X_t-x^*\|^2 - c_0\lambda\sqrt{\mu}\langle X_t-x^*,\nabla f(X_t) \rangle,
    \end{split}
\end{equation*}
which together with the inequality $\langle X_t-x^*,\nabla f(X_t) \rangle \geq f(X_t)-f^*+\frac{\mu}{2}\|X_t-x^*\|^2$ by $\mu$-strong convexity suggests that
\begin{equation*}
    \begin{split}
        \frac{\dif V_t}{\dif t} &\leq -2C\sqrt{\mu}V_t - \frac{c_1-\lambda}{2}\sqrt{\mu}\|\dot X_t\|^2 - \left(\frac{c_1-\lambda}{2}-C\right)\sqrt{\mu}\|\dot X_t + \lambda \sqrt{\mu}(X_t-x^*)\|^2 \\
        &\quad - (\lambda-2C)c_0\sqrt{\mu}(f(X_t)-f^*) - \frac{\lambda\mu^{\frac{3}{2}}}{2}(c_0-(c_1-\lambda)\lambda)\|X_t-x^*\|^2.
    \end{split}
\end{equation*}
Because $c_0>0$ and $c_1>0$, we pick $\lambda$ such that $0<\lambda < c_1$ and $(c_1-\lambda)\lambda\leq c_0$. Moreover, we pick $C$ such that $0<C\leq \frac{c_1-\lambda}{2} \wedge \frac{\lambda}{2}$. Then
$$
\frac{\dif V_t}{\dif t}\leq -2C\sqrt{\mu}V_t \implies V_t\leq V_0\cdot\me^{-2C\sqrt{\mu}t}.
$$
In particular, if $c_1^2\leq 4c_0$, we pick $\lambda = c_1/2$ and $C=c_1/4$. If $c_1^2>4c_0$, we pick $\lambda = \frac{c_1 + \sqrt{c_1^2-4c_0}}{2}$ and $C=\frac{c_1 - \sqrt{c_1^2-4c_0}}{4}$. The conclusion then follows by bounding $V_0=c_0(f(x_0)-f^*)+\frac{\lambda^2\mu}{2}\|x_0-x^*\|^2\leq (c_0+\lambda^2)(f(x_0)-f^*)$ because $\frac{\mu}{2}\|x_0-x^*\|^2\leq f(x_0)-f^*$.

When $c_0\leq 0$ or $c_1\leq 0$, $f(X_t)-f^*$ is not guaranteed to converge. A simple counterexample is the harmonic oscillator, with $f(x)=\mu x^2/2$ for $x\in \bbR$,
\begin{equation}\label{eq:harmonic oscillator}
    \ddot X_t + c_1\sqrt{\mu}\dot X_t + c_0\mu X_t=0,
\end{equation}
starting from the initial position $X(0) =x_0$ with velocity $\dot{X}(0)=0$.
Since (\ref{eq:harmonic oscillator}) is a second-order linear ODE with constant coefficients, its general solutions admit closed forms. Consider the characteristic equation $w^2 + c_1\sqrt{\mu}w+c_0\mu =0$ with discriminant $\Delta = (c_1^2-4c_0)\mu$.

If $c_1^2=4c_0$, then there are two identical real roots $w=w_{1}=w_2=-c_1\sqrt{\mu}/2$, and $X_t=(\alpha_0+\alpha_1t)\me^{wt}$ for some real numbers $\alpha_0$ and $\alpha_1$. By the initial condition $\dot X(0)=0$, we find $\alpha_1=-\alpha_0w$. Then $X_t=\alpha_0(1-wt)\me^{wt}$. To achieve $X_t\to 0$, we need $w<0$, which means $c_1>0$ and $c_0=c_1^2/4>0$ too.

If $c_1^2>4c_0$, then there are two distinct real roots $w_1=\frac{-c_1+ \sqrt{c_1^2-4c_0}}{2}\sqrt{\mu}$ and $w_2=\frac{-c_1- \sqrt{c_1^2-4c_0}}{2}\sqrt{\mu}$,
with $w_1 > w_2$, and $X_t= \alpha_0\me^{w_1t}+\alpha_1\me^{w_2t}$. By the initial condition $\dot X(0)=0$, we find $\alpha_0w_1+\alpha_1w_2=0$.
If $w_2=0$, then $w_1>0$ and hence $\alpha_0 =0$ and $X_t \equiv \alpha_1=x_0$, which contradicts $X_t \to 0$ for any $x_0\not= 0$.
If $w_2\ne 0$, then $\alpha_1=-\alpha_0w_1/w_2$, and $X_t=\alpha_0(\me^{w_1 t}-\frac{w_1}{w_2}\me^{w_2 t})$ with $w_1\not=w_2$.  To achieve $X_t\to 0$, we also need $w_1<0$. Then $-c_1+\sqrt{c_1^2-4c_0}<0$, which implies that $c_0>0$ and $c_1>0$.

If $c_1^2<4c_0$, then there are two complex roots $w_{1,2}=\frac{-c_1\pm \sqrt{4c_0-c_1^2}i}{2}\sqrt{\mu}$, and $X_t$ is a linear combination of $\me^{-c_1\sqrt{\mu}t/2}\sin(\sqrt{\mu(4c_0-c_1^2)}t/2)$ and $\me^{-c_1\sqrt{\mu}t/2}\cos(\sqrt{\mu(4c_0-c_1^2)}t/2)$. We need $c_1>0$ to make $X_t \to 0$. Then $c_0>c_1^2/4>0$ too.

\subsubsection{Proof of Proposition \ref{pro:high-res-ODE NAG and HB}}

For (\ref{eq:high-res-ODE-NAG-SC}), we consider the Lyapunov function
$$
V_t = (1+\sqrt{\mu s})(f(X_t)-f^*) + \frac{1}{2}\|\dot X_t + \sqrt{\mu}(X_t-x^*)+\sqrt{s}\nabla f(X_t)\|^2.
$$
Using the ODE (\ref{eq:high-res-ODE-NAG-SC}), we have
\begin{equation*}
    \begin{split}
        \frac{\dif V_t}{\dif t} &= (1+\sqrt{\mu s}) \langle \nabla f(X_t), \dot X_t \rangle + \langle \dot X_t + \sqrt{\mu}(X_t-x^*)+\sqrt{s}\nabla f(X_t), \ddot X_t + \sqrt{\mu}\dot X_t + \sqrt{s}\nabla^2 f(X_t)\dot X_t \rangle \\
        &= (1+\sqrt{\mu s}) \langle \nabla f(X_t), \dot X_t \rangle + \langle \dot X_t + \sqrt{\mu}(X_t-x^*)+\sqrt{s}\nabla f(X_t), -\sqrt{\mu}\dot X_t - (1+\sqrt{\mu s})\nabla f(X_t) \rangle \\
        &= -\sqrt{\mu}\|\dot X_t\|^2 - \sqrt{\mu}\langle \dot X_t, \sqrt{\mu}(X_t-x^*)+\sqrt{s}\nabla f(X_t)\rangle \\
        &\quad - (1+\sqrt{\mu s})\langle \nabla f(X_t),  \sqrt{\mu}(X_t-x^*)+\sqrt{s}\nabla f(X_t)\rangle \\
        &= -\frac{\sqrt{\mu}}{2}\|\dot X_t\|^2 - \frac{\sqrt{\mu}}{2}\|\dot X_t + \sqrt{\mu}(X_t-x^*)+\sqrt{s}\nabla f(X_t)\|^2 \\
        &\quad + \frac{\sqrt{\mu}}{2}\|\sqrt{\mu}(X_t-x^*)+\sqrt{s}\nabla f(X_t)\|^2 - (1+\sqrt{\mu s})\sqrt{\mu}\langle X_t-x^*,\nabla f(X_t) \rangle \\
        &\quad- (1+\sqrt{\mu s})\sqrt{s} \|\nabla f(X_t)\|^2 \\
        &= -\sqrt{\mu}V_t + \sqrt{\mu}(1+\sqrt{\mu s})(f(X_t)-f^*) - \frac{\sqrt{\mu}}{2}\|\dot X_t\|^2  \\
        &\quad + \frac{\mu^{\frac{3}{2}}}{2}\|X_t-x^*\|^2 - \sqrt{\mu}\langle X_t-x^*,\nabla f(X_t) \rangle - \sqrt{s}(1+\frac{\sqrt{\mu s}}{2})\|\nabla f(X_t)\|^2.
     \end{split}
\end{equation*}
By the strong convexity, $\langle X_t-x^*,\nabla f(X_t) \rangle \geq f(X_t)-f^*+\frac{\mu}{2}\|X_t-x^*\|^2$. Then the above display yields
\begin{equation*}
    \frac{\dif V_t}{\dif t}\leq -\sqrt{\mu}V_t + \sqrt{\mu}\sqrt{\mu s}(f(X_t)-f^*)- \sqrt{s}(1+\frac{\sqrt{\mu s}}{2})\|\nabla f(X_t)\|^2.
\end{equation*}
By the strong convexity again, we have $f(X_t)-f^*\leq \frac{1}{2\mu}\|\nabla f(X_t)\|^2$. Hence
$$
\frac{\dif V_t}{\dif t}\leq -\sqrt{\mu}V_t - \frac{\sqrt{s}}{2}(1+\sqrt{\mu s})\|\nabla f(X_t)\|^2\leq -\sqrt{\mu}V_t.
$$
The conclusion follows by integrating with $t$.

For (\ref{eq:high-res-ODE-HB}), we consider the Lyapunov function
$$
V_t = (1+\sqrt{\mu s})(f(X_t)-f^*) + \frac{1}{2}\|\dot X_t + \sqrt{\mu}(X_t-x^*)\|^2.
$$
Using the ODE (\ref{eq:high-res-ODE-HB}), we have
\begin{equation*}
    \begin{split}
        \frac{\dif V_t}{\dif t}&= (1+\sqrt{\mu s})\langle \nabla f(X_t), \dot X_t \rangle + \langle \dot X_t + \sqrt{\mu}(X_t-x^*), \ddot X_t + \sqrt{\mu}\dot X_t  \rangle \\
        &= (1+\sqrt{\mu s})\langle \nabla f(X_t), \dot X_t \rangle - \langle \dot X_t + \sqrt{\mu}(X_t-x^*), \sqrt{\mu}\dot X_t + (1+\sqrt{\mu s})\nabla f(X_t) \rangle \\
        &= -\sqrt{\mu}\|\dot X_t\|^2 - \sqrt{\mu}\langle \dot X_t,\sqrt{\mu}(X_t-x^*) \rangle - \sqrt{\mu}(1+\sqrt{\mu s})\langle X_t-x^*, \nabla f(X_t) \rangle  \\
        &= -\frac{\sqrt{\mu}}{2}\|\dot X_t\|^2  - \frac{\sqrt{\mu}}{2}\|\dot X_t + \sqrt{\mu}(X_t-x^*)\|^2 + \frac{\mu^{\frac{3}{2}}}{2}\|X_t-x^*\|^2 - \sqrt{\mu}(1+\sqrt{\mu s})\langle X_t-x^*, \nabla f(X_t) \rangle\\
        &\leq -\sqrt{\mu} V_t + \sqrt{\mu}(1+\sqrt{\mu s})(f(X_t)-f^*) + \frac{\mu^{\frac{3}{2}}}{2}\|X_t-x^*\|^2 - \sqrt{\mu}(1+\sqrt{\mu s})\langle X_t-x^*, \nabla f(X_t) \rangle.
    \end{split}
\end{equation*}
By the strong convexity, $\langle X_t-x^*,\nabla f(X_t) \rangle \geq f(X_t)-f^*+\frac{\mu}{2}\|X_t-x^*\|^2$. Then the above display yields
$$
\frac{\dif V_t}{\dif t} \leq -\sqrt{\mu}V_t - \frac{\mu^{2}\sqrt{s}}{2}\|X_t-x^*\|^2\leq -\sqrt{\mu}V_t.
$$
The conclusion follows by integrating with $t$.

\subsection{Proof of Proposition \ref{pro:high-res-ODE-C}}

The construction of our continuous-time Lyapunov function is motivated by \cite{su2016differential} and \cite{shi2021understanding}.
We consider the auxiliary-energy term defined as
\begin{equation*}
    V^A_t = \frac{1}{2}\left\|r(X_t-x^*)+t\left(\dot X_t + \frac{\beta\sqrt{s}}{\gamma}\nabla f(X_t)\right)\right\|^2.
\end{equation*}
Using the ODE (\ref{eq:high-res-ODE-C}), for $t\ge t_0$ we have
\begin{equation*}
    \begin{split}
        \frac{\dif V_t^A}{\dif t}
        &= \big\langle r(X_t-x^*)+t\big(\dot X_t + \frac{\beta\sqrt{s}}{\gamma}\nabla f(X_t)\big),  r\dot X_t + \dot X_t +  \frac{\beta\sqrt{s}}{\gamma}\nabla f(X_t) + t\big(\ddot X_t + \frac{\beta\sqrt{s}}{\gamma}\nabla^2 f(X_t)\dot X_t\big)   \big\rangle \\
        & = \big\langle r(X_t-x^*)+t\big(\dot X_t + \frac{\beta\sqrt{s}}{\gamma}\nabla f(X_t)\big), - \big(t+ (\frac{r+1}{2}-\frac{\beta}{\gamma})\sqrt{s} \big)\nabla f(X_t) \big\rangle \\
        &= -r\big(t+ (\frac{r+1}{2}-\frac{\beta}{\gamma})\sqrt{s} \big) \langle X_t-x^*,\nabla f(X_t)    \rangle - t\big(t+ (\frac{r+1}{2}-\frac{\beta}{\gamma})\sqrt{s} \big) \langle \dot X_t, \nabla f(X_t)\rangle \\
        &\quad - \frac{\beta\sqrt{s}}{\gamma}t\big(t+ (\frac{r+1}{2}-\frac{\beta}{\gamma})\sqrt{s} \big) \|\nabla  f(X_t)\|^2.
    \end{split}
\end{equation*}
Let $C\geq 0$ be a constant to be chosen later, and introduce a factor of $\frac{t+C\sqrt{s}}{t}$ for technical adjustment. Because $\frac{t+C\sqrt{s}}{t}=1+\frac{C\sqrt{s}}{t}$ is decreasing in $t$, we have
\begin{equation}\label{eq:supp-df-kinetic}
    \begin{split}
        \frac{\dif}{\dif t}\left(\frac{t+C\sqrt{s}}{t} V_t^A\right) &=\frac{\dif}{\dif t}\left(\frac{t+C\sqrt{s}}{t}\right)\cdot V_t^A + \frac{t+C\sqrt{s}}{t}\cdot\frac{\dif V_t^A}{\dif t}  \\
        &\leq \frac{t+C\sqrt{s}}{t}\cdot\frac{\dif V_t^A}{\dif t} \\
        &= -r \frac{t+C\sqrt{s}}{t}\big(t+ (\frac{r+1}{2}-\frac{\beta}{\gamma})\sqrt{s} \big) \langle X_t-x^*,\nabla f(X_t)\rangle \\
        &\quad - (t+C\sqrt{s})\big(t+ (\frac{r+1}{2}-\frac{\beta}{\gamma})\sqrt{s} \big) \langle \dot X_t, \nabla f(X_t)\rangle \\
        &\quad -\frac{\beta\sqrt{s}}{\gamma}(t+C\sqrt{s})\big(t+ (\frac{r+1}{2}-\frac{\beta}{\gamma})\sqrt{s} \big) \|\nabla  f(X_t)\|^2.
    \end{split}
\end{equation}
To eliminate the term of $\langle \dot X_t, \nabla f(X_t) \rangle$, we define the potential-energy term as
\begin{equation*}
    V_t^P = (t+C\sqrt{s})\big(t+ (\frac{r+1}{2}-\frac{\beta}{\gamma})\sqrt{s} \big)(f(X_t)-f^*).
\end{equation*}
Then
\begin{equation}\label{eq:supp-df-potential}
    \begin{split}
        \frac{\dif V_t^P}{\dif t} &= \big(2t+(\frac{r+1}{2}-\frac{\beta}{\gamma}+C)\sqrt{s}\big)(f(X_t)-f^*)\\
        &\quad + (t+C\sqrt{s})\big(t+ (\frac{r+1}{2}-\frac{\beta}{\gamma})\sqrt{s} \big)\langle \nabla f(X_t), \dot X_t \rangle.
    \end{split}
\end{equation}
Define the continuous Lyapunov function as $V_t = V_t^P + \frac{t+C\sqrt{s}}{t}V_t^A$. Set $t_0^* = t_0\vee \{2|C|\sqrt{s}\} \vee \{2|\frac{r+1}{2}-\frac{\beta}{\gamma}|\sqrt{s}\}>0$. Then for $t\ge t_0^*$, we have $V_t\ge 0$ and
\begin{align} \label{eq:t0-star}
(t+C\sqrt{s})\big(t+ (\frac{r+1}{2}-\frac{\beta}{\gamma})\sqrt{s} \big) \ge \frac{t}{2}\cdot \frac{t}{2}= \frac{t^2}{4}.
\end{align}
Combining (\ref{eq:supp-df-kinetic}) and (\ref{eq:supp-df-potential}), we obtain
\begin{equation*}
    \begin{split}
        \frac{\dif V_t}{\dif t} &\leq  -r \frac{t+C\sqrt{s}}{t}\big(t+ (\frac{r+1}{2}-\frac{\beta}{\gamma})\sqrt{s} \big) \langle X_t-x^*,\nabla f(X_t)\rangle \\
        &\quad + \big(2t+(\frac{r+1}{2}-\frac{\beta}{\gamma}+C)\sqrt{s}\big)(f(X_t)-f^*) \\
        &\quad -\frac{\beta\sqrt{s}}{\gamma}(t+C\sqrt{s})\big(t+ (\frac{r+1}{2}-\frac{\beta}{\gamma})\sqrt{s} \big) \|\nabla  f(X_t)\|^2,
    \end{split}
\end{equation*}
which together with the convexity inequality $\langle X_t-x^*,\nabla f(X_t)\rangle \geq f(X_t)-f^*$ suggests that for $r\geq 0$ and $t\geq t_0^*$,
\begin{equation*}
    \begin{split}
        \frac{\dif V_t}{\dif t} &\leq -\underbrace{\left((r-2)t+(r-1)C\sqrt{s}+(r-1+\frac{rC\sqrt{s}}{t})(\frac{r+1}{2}-\frac{\beta}{\gamma})\sqrt{s}\right)}_{\one_t}\cdot (f(X_t)-f^*) \\
        &\quad - \underbrace{\frac{\beta}{\gamma}(t+C\sqrt{s})\big(t+ (\frac{r+1}{2}-\frac{\beta}{\gamma})\sqrt{s} \big)}_{\two_t}\cdot\sqrt{s} \|\nabla  f(X_t)\|^2.
    \end{split}
\end{equation*}
From (\ref{eq:t0-star}), if $\frac{\beta}{\gamma}\geq 0$, then $\two_t\ge \frac{\beta}{4\gamma}t^2\ge 0$ for $t\ge t_0^*$. It remains to deal with $\one_t$.
% This condition is a prerequisite.

When $r>2$, we set $C=0$. Then $\one_t = (r-2)t+(r-1)(\frac{r+1}{2}-\frac{\beta}{\gamma})\sqrt{s} \geq 0$ for $t \geq -\frac{r-1}{r-2}(\frac{r+1}{2}-\frac{\beta}{\gamma})\sqrt{s}$. In this case, we pick $t_1=t_0^* \vee \{-\frac{r-1}{r-2}(\frac{r+1}{2}-\frac{\beta}{\gamma})\sqrt{s} \}$.

When $r=2$, $\one_t=C\sqrt{s} + (1+\frac{2C\sqrt{s}}{t})(\frac{3}{2}-\frac{\beta}{\gamma})\sqrt{s}$.
We set the constant $C>\frac{\beta}{\gamma}-\frac{3}{2}$ (e.g., $\frac{\beta}{\gamma}$). Then when $t$ is large (e.g., $t\geq \frac{4\beta}{3\gamma}(\frac{\beta}{\gamma}-\frac{3}{2})\sqrt{s}$), $\one_t\geq 0$ always holds. In this case, we pick $t_1 = t_0^* \vee \{\frac{4\beta}{3\gamma}(\frac{\beta}{\gamma}-\frac{3}{2})\sqrt{s} \}$.

To summarize, when $r\geq 2$ and $\frac{\beta}{\gamma}\geq 0$, there exist $C\ge 0$ and $t_1$ such that when $t\geq t_1$, $\frac{\dif V_t}{\dif t}\leq 0$.
If further $\frac{\beta}{\gamma}> 0$, then $-\frac{\dif V_t}{\dif t}\ge \frac{\beta\sqrt{s}}{4\gamma} t^2 \|\nabla f(X_t)\|^2$.
The remaining proof is similar to the proof of Theorem \ref{thm:converge-c}, but in a continuous way. We  sketch the reasoning below.

If $\frac{\dif V_t}{\dif t}\leq 0$, then for $t\ge t_1$,
$$
V_{t_1}\geq V_t \geq V_t^P = (t+C\sqrt{s})\big(t+ (\frac{r+1}{2}-\frac{\beta}{\gamma})\sqrt{s} \big)(f(X_t)-f^*)\ge \frac{t^2}{4}(f(X_t)-f^*),
$$
which implies that $f(X_t)-f^* = O\left(\frac{V_{t_1}}{t^2}\right)$.

If further $-\frac{\dif V_t}{\dif t}\ge \frac{\beta\sqrt{s}}{4\gamma} t^2 \|\nabla f(X_t)\|^2$ with $\frac{\beta}{\gamma}>0$, then for $t>t_1$, we have $V_t\ge 0$, and
\begin{equation*}
    \begin{split}
        V_{t_1}\geq V_{t_1}-V_t = (-V_t)-(-V_{t_1}) = \int_{t_1}^t \frac{\dif}{\dif u}(-V_u)\dif u \ge \frac{\beta\sqrt{s}}{4\gamma} \int_{t_1}^t u^2\|\nabla f(X_u)\|^2\dif u \\
        \geq \frac{\beta\sqrt{s}}{4\gamma} \inf_{t_1\leq u\leq t}\|\nabla f(X_u)\|^2 \int_{t_1}^t u^2 \dif u =\frac{\beta\sqrt{s}(t^3-t_1^3)}{12\gamma}\inf_{t_1\leq u\leq t}\|\nabla f(X_u)\|^2,
    \end{split}
\end{equation*}
which implies that
$$
\inf_{t_1\leq u\leq t}\|\nabla f(X_u)\|^2\le\frac{12\gamma V_{t_1}}{\beta\sqrt{s}(t^3-t_1^3)}=O\left(\frac{\gamma V_{t_1}}{\beta\sqrt{s}(t^3-t_1^3)}\right).
$$

\section{Technical details in Section \ref{sec:HAG-interpretation}}\label{sec:tech for section HAG}

\subsection{Single-variable form for HAG (\ref{eq:HAG})}
We derive the single-variable form (\ref{eq:HAG-single}) for HAG. From (\ref{eq:HAG}) solving the first equation for $u_k$, and plugging the expressions of $u_k$ and $u_{k+1}$ into the second equation, we obtain
\begin{equation*}
    \begin{split}
        \frac{x_{k+2}-x_{k+1}+a_{k+1}\nabla f(x_{k+1})}{\sqrt{a_{k+1}b_{k+1}}} &= (b_k-1)\frac{x_{k+1}-x_{k}+a_{k}\nabla f(x_{k})}{\sqrt{a_{k}b_{k}}}\\
        &\quad- \sqrt{a_kb_k}\nabla f(x_k) - \phi_k(\nabla f(x_{k+1})-\nabla f(x_k)).
    \end{split}
\end{equation*}
Rearrange the above display to conclude.

\subsection{Verification of monotonicity condition for (\ref{eq:HAG-C})}
We verify that the monotonicity condition in Theorem \ref{thm:converge-c} is satisfied by (\ref{eq:HAG-C}), i.e., HAG under the configuration (\ref{eq:config-HAG-C}). First, we notice that (\ref{eq:HAG-C}) can be put into (\ref{eq:extended-NAG-C-single}) with
\begin{equation*}
    \beta_{k-1} = \frac{1}{\sigma_{k+1}}\left(c_2\sqrt{c_0}-\frac{c_0}{b_{k-1}}\right), \quad \gamma_k = c_0\left(\frac{1}{b_{k-1}}+\frac{1}{b_k}\right)-\sigma_{k+1}(\beta_k-\beta_{k-1}).
\end{equation*}
When $\alpha_k=\frac{k+r}{r}$, we have $\sigma_{k+1}=\frac{\alpha_k-1}{\alpha_{k+1}}=\frac{k}{k+r+1}$, $b_{k-1}=1+\sigma_{k+1}=\frac{2k+r+1}{k+r+1}$, and $\frac{\alpha_{k+1}}{\alpha_k}=\frac{k+r+1}{k+r}$. Without loss of generality, take $c_0=1$. 
Substituting $\sigma_{k+1}$, $b_{k-1}$ and $b_k$ into the display above and after some algebra, we obtain
\begin{equation*}
    \begin{split}
        \frac{\widetilde \alpha_k}{\alpha_k} &= \gamma_k\frac{\alpha_{k+1}}{\alpha_k} - \beta_k\left(\frac{\alpha_{k+1}}{\alpha_k}-1\right) \\
        &= \gamma_k \frac{k+r+1}{k+r}-\beta_k\frac{1}{k+r} \\
        &= \frac{2k^2+(7-2c_2+4r)k + 6-3c_2+7r-c_2r+2r^2}{(k+r)(2k+r+3)},
    \end{split}
\end{equation*}
which is monotone in $k$ when $k\geq K$ for some $K$ depending only on $c_2$ and $r$.

\section{Numerical experiments} \label{sec:experiment}

We apply algorithms (\ref{eq:extended-NAG-SC-single-2}) and (\ref{eq:extended-NAG-C-shi}) to various objective functions
(strongly or general convex) under different parameters and step sizes.
It is particularly of interest to compare the performances from the parameter choices falling inside (or on the boundaries) versus outside
the sufficient conditions for achieving accelerated convergence in our theoretical results.

The numerical results below are overall consistent with our theoretical results in both strongly and general convex settings.
Better performances are observed from algorithms inside (or on the boundaries) our sufficient conditions for accelerated convergence,
that is, with more friction (over- or critical-damping, $c_1^2 \ge 4c_0$, vs under-damping, $c_1^2 < 4c_0$, in (\ref{eq:extended-NAG-SC-single-2})
and $r \ge 2$ vs $r <1$ in (\ref{eq:extended-NAG-C-shi}))
and larger gradient correction ($c_2^2 \ge c_0$ vs $c_2^2 < c_0$ in (\ref{eq:extended-NAG-SC-single-2})
and $\beta/\gamma \ge 1/2$ vs $\beta/\gamma < 1/2$ in (\ref{eq:extended-NAG-C-shi})).
The effect of gradient correction becomes more significant as the step size increases.

\begin{figure}[htbp]
    \centering
    \begin{subfigure}{0.49\textwidth}
       \includegraphics[width=\textwidth]{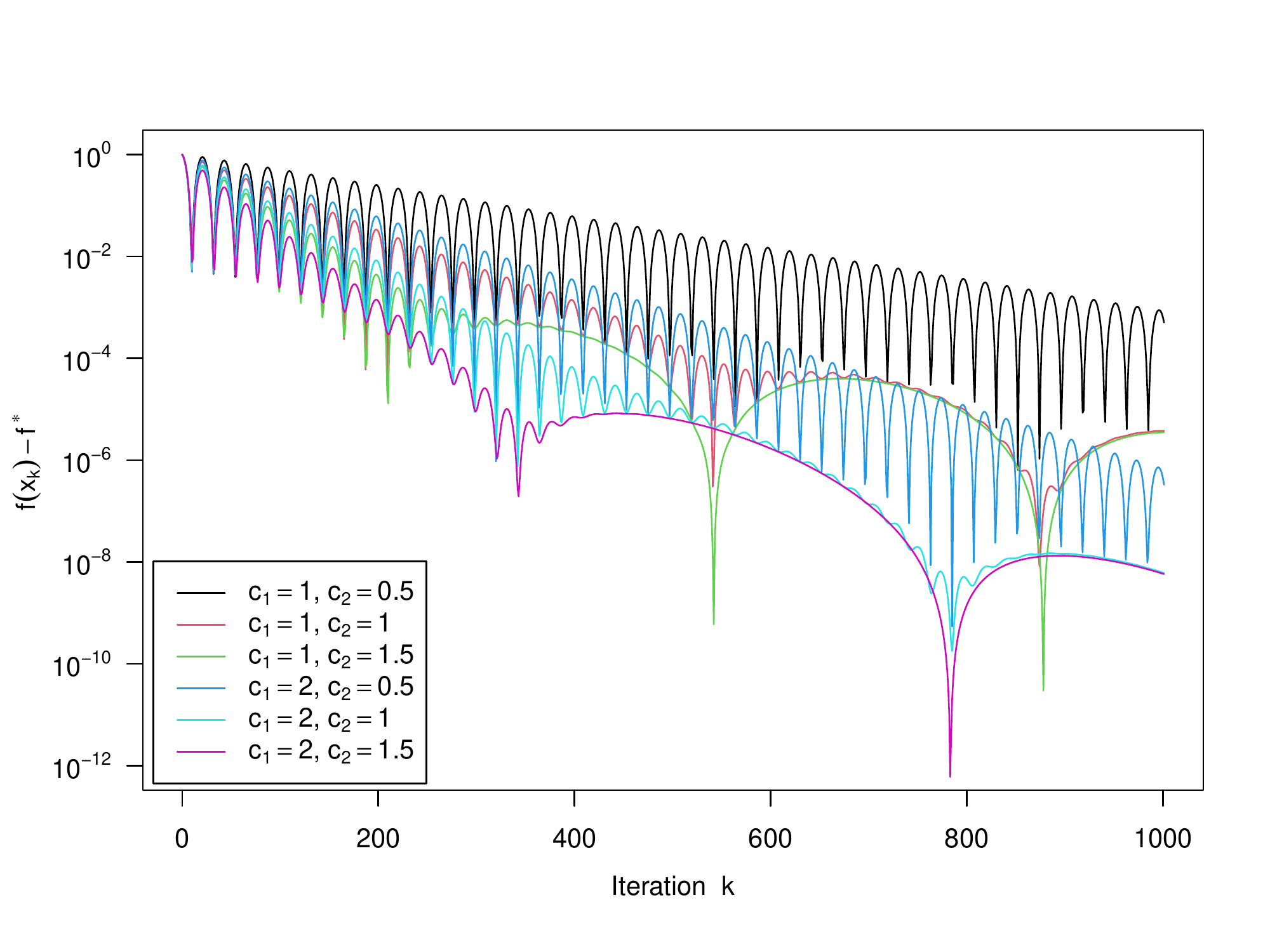}
       \vspace{-0.45in}
       \caption{Errors $f(x_k)-f^*$ ($s=0.01$)}
    \end{subfigure}
    \hfill
    \begin{subfigure}{0.49\textwidth}
       \includegraphics[width=\textwidth]{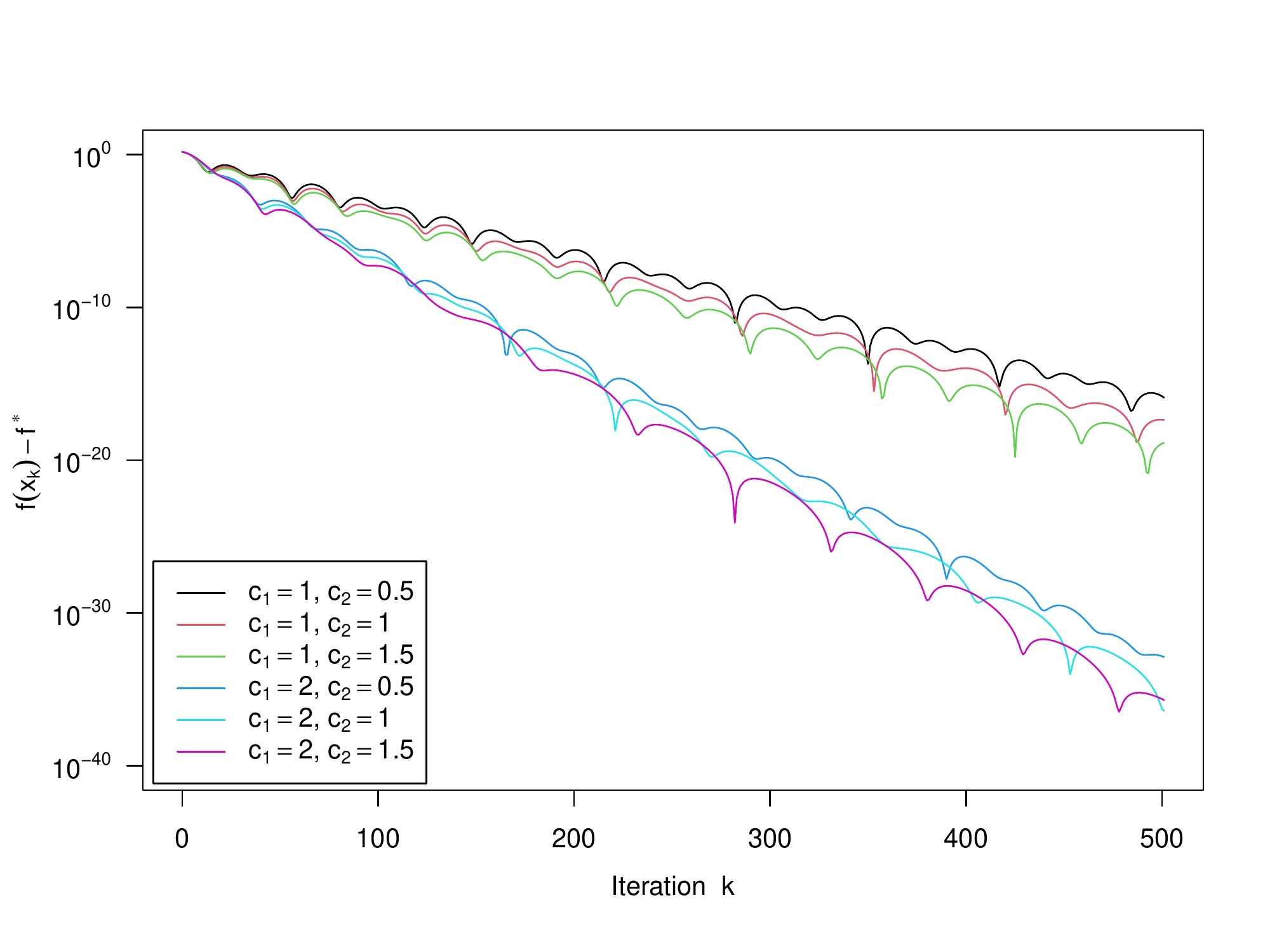}
       \vspace{-0.45in}
       \caption{Errors $f(x_k)-f^*$ ($s=0.01$)}
    \end{subfigure}
    \vspace{-0.2in}

    \begin{subfigure}{0.49\textwidth}
       \includegraphics[width=\textwidth]{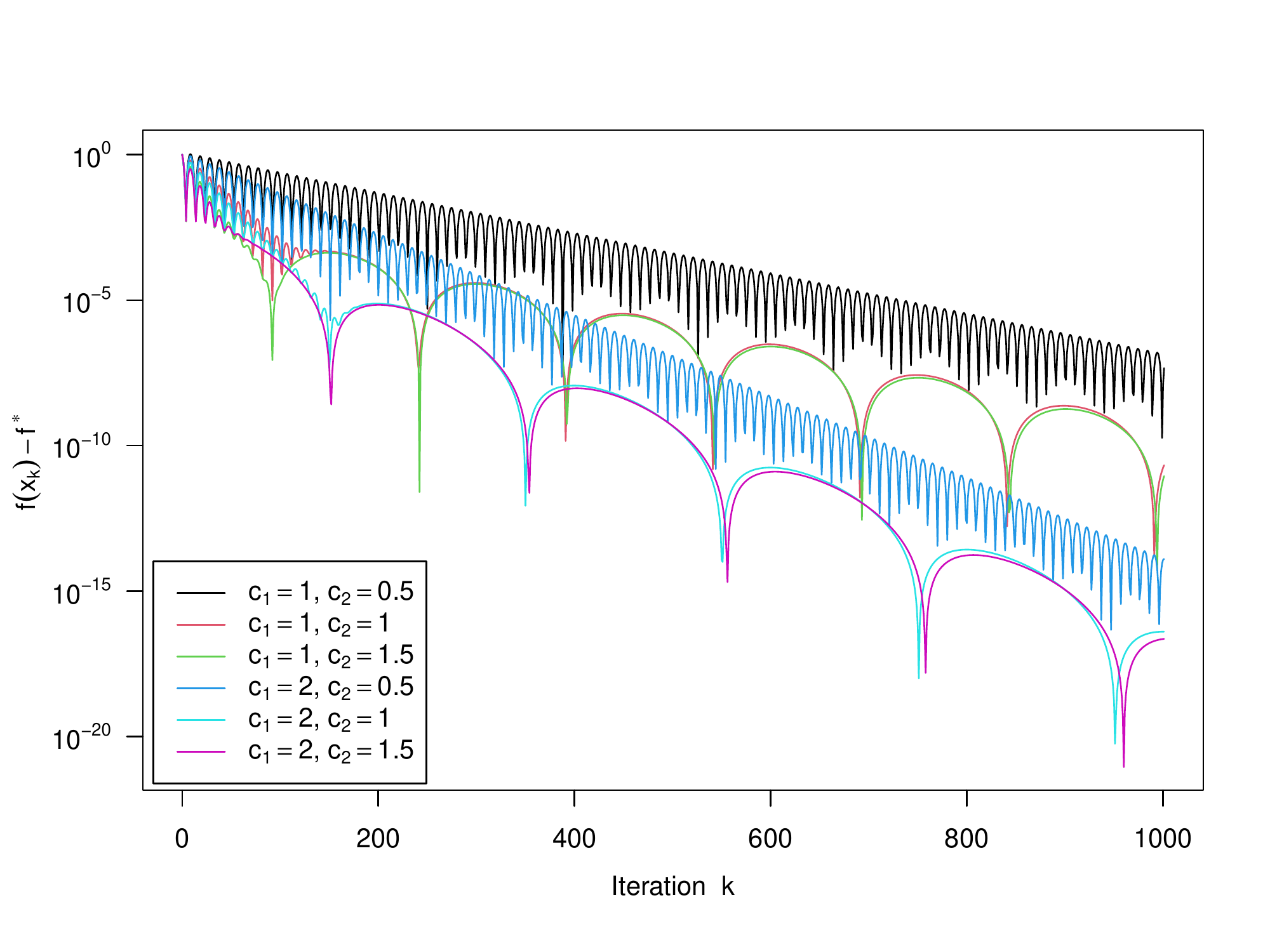}
       \vspace{-0.45in}
       \caption{Errors $f(x_k)-f^*$ ($s=0.05$)}
    \end{subfigure}
    \hfill
    \begin{subfigure}{0.49\textwidth}
       \includegraphics[width=\textwidth]{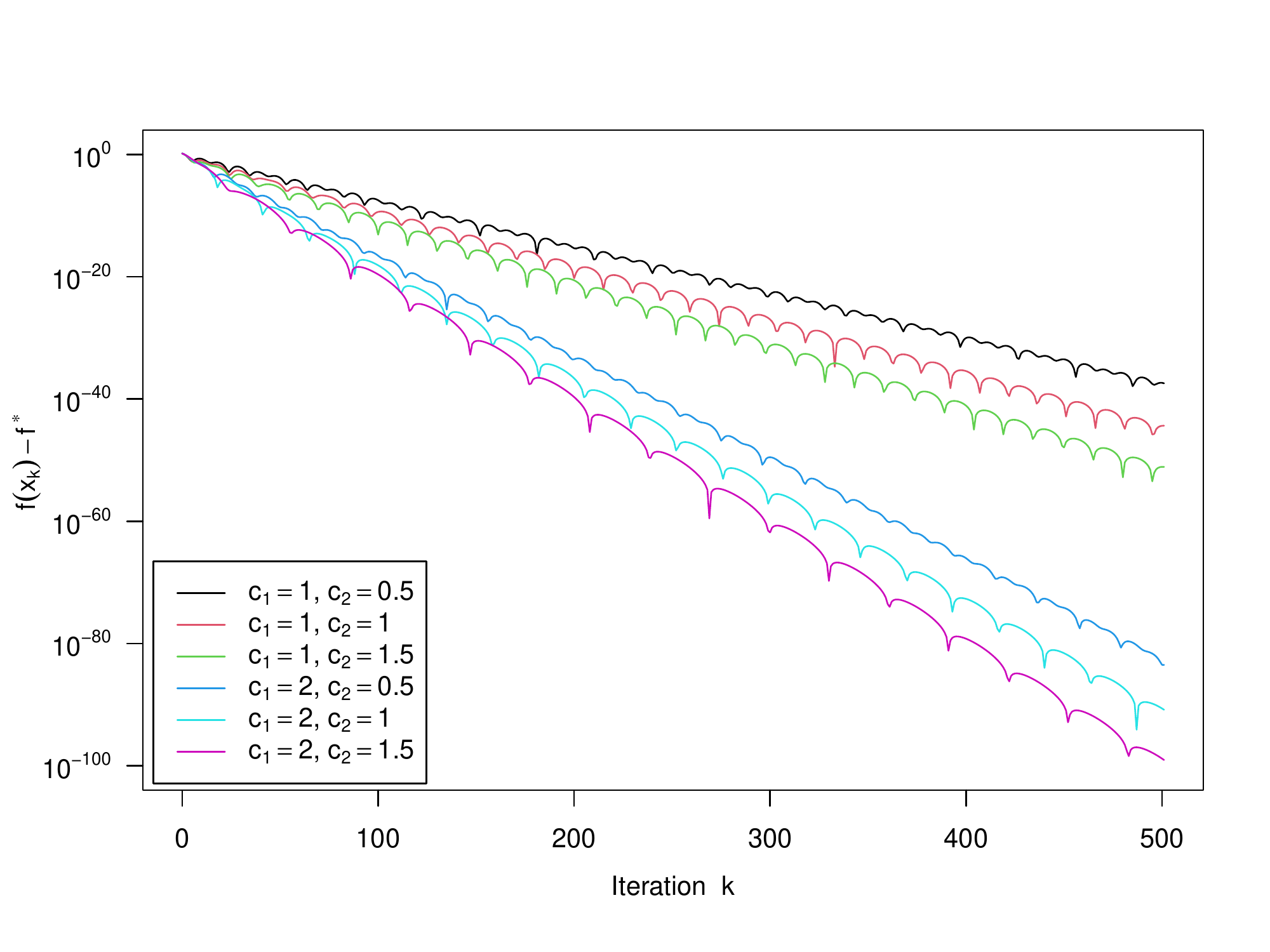}
       \vspace{-0.45in}
       \caption{Errors $f(x_k)-f^*$ ($s=0.05$)}
    \end{subfigure}
    \vspace{-0.2in}

    \begin{subfigure}{0.49\textwidth}
       \includegraphics[width=\textwidth]{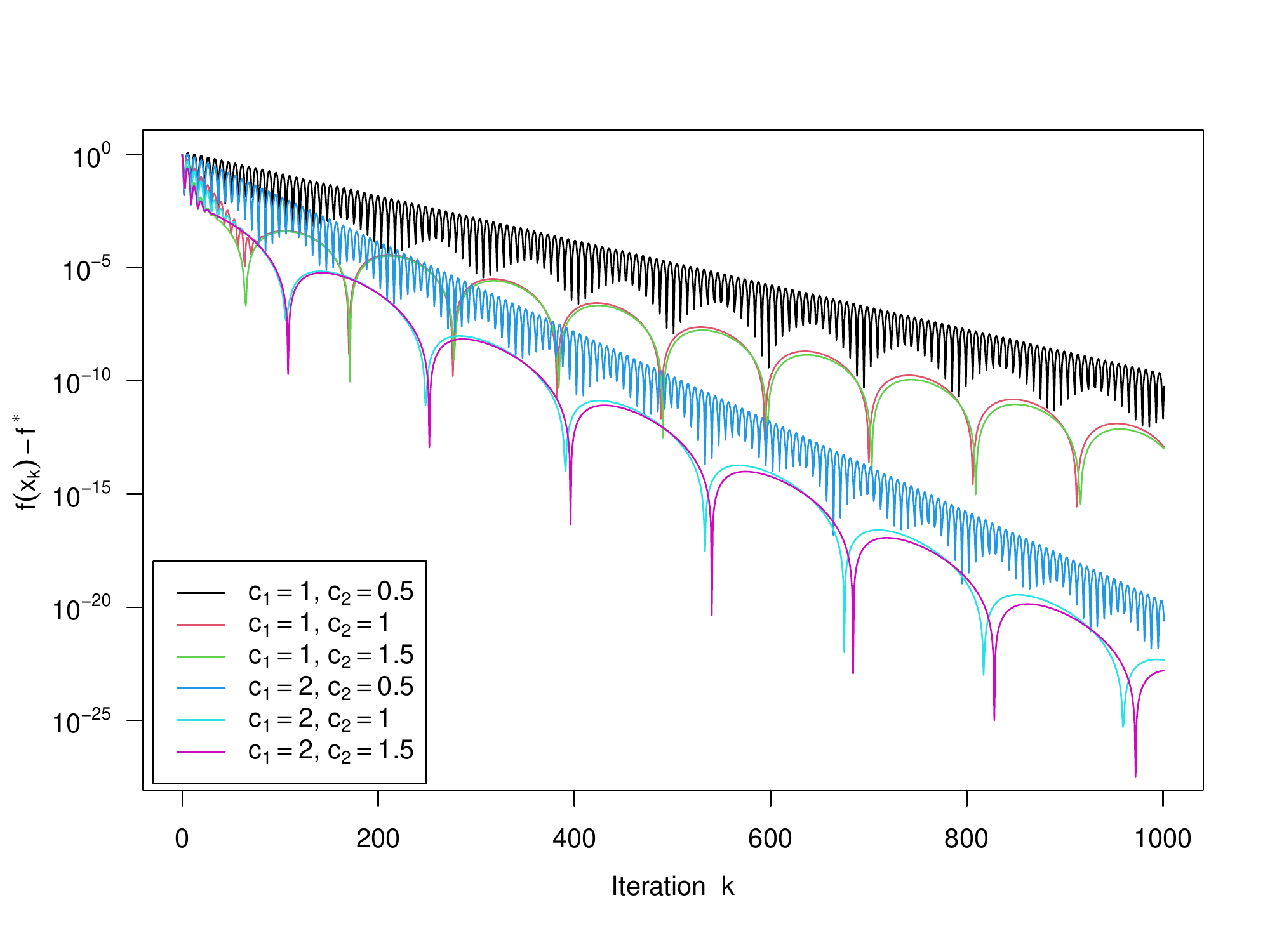}
       \vspace{-0.45in}
       \caption{Errors $f(x_k)-f^*$ ($s=0.1$)}
    \end{subfigure}
    \hfill
    \begin{subfigure}{0.49\textwidth}
       \includegraphics[width=\textwidth]{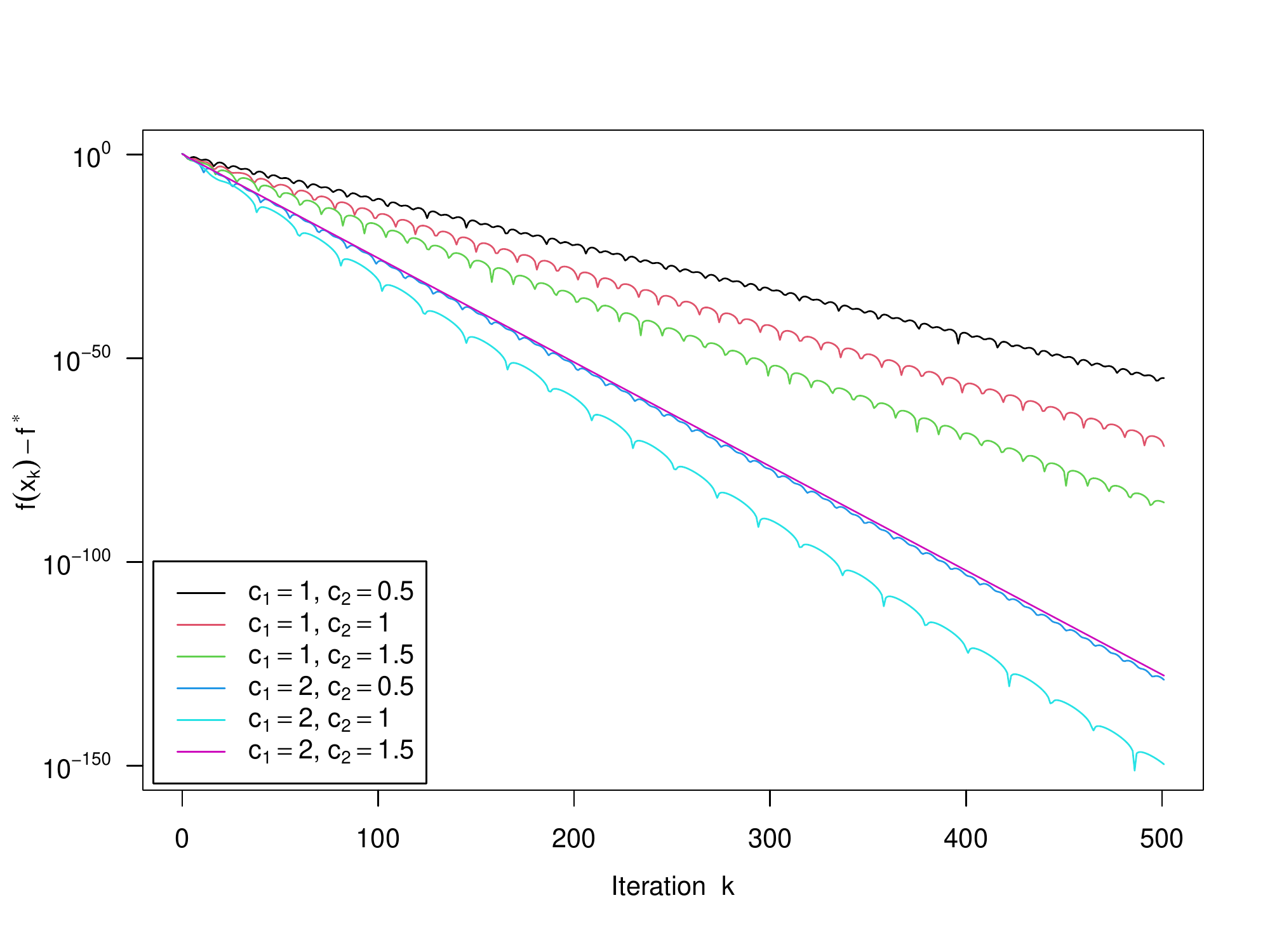}
       \vspace{-0.45in}
       \caption{Errors $f(x_k)-f^*$ ($s=0.1$)}
    \end{subfigure}

    \caption{Minimizing strongly convex functions by (\ref{eq:extended-NAG-SC-single-2}). The left column is for ill-conditioned $f(x_1,x_2)=5\times 10^{-3}x_1^2 + x_2^2$, and the right column is for well-conditioned $f(x_1,x_2) = 5\times 10^{-1}x_1^2 + x_2^2$. Fix $c_0=1$. The initial iterates are $x_0=(1,1)$ and $x_1=x_0-2s\nabla f(x_0)/(1+\sqrt{\mu s})$.}
    \label{fig: strongly convex}
\end{figure}

\textbf{Strongly convex setting.}
Figure \ref{fig: strongly convex} presents several numerical results of (\ref{eq:extended-NAG-SC-single-2}) applied to two quadratic functions on $\bbR^2$
similarly as in \cite{su2016differential}, with $c_0$ fixed at $1$, $c_1=1$, $2$ (boundary) and $c_2= 1/2, 1$ (boundary), $3/2$. All remainder terms $R_1$, $R_2$, and $R_3$ are taken to be $0$. In particular, $c_1=1, 2$ corresponds to under-damping ($c_1^2<4c_0$) and critical-damping ($c_1^2=4c_0$) respectively, and $c_2=1/2, 1, 3/2$ corresponds to the leading constant of the gradient correction coefficient in (\ref{eq:extended-NAG-SC-single-2}) being $0, 1/2, 1$ respectively. For NAG-SC (\ref{eq:NAG-SC-single}), $c_1=2$ and $c_2=3/2$, and for heavy-ball (\ref{eq:HB}), $c_1=2$ and $c_2=1/2$,
but both with nonzero remainder terms $R_2$ and $R_3$.
Hence the algorithms studied are not exactly NAG-SC or heavy-ball.

From Figure \ref{fig: strongly convex} we observe that critical-damping ($c_1=2$) results in a faster convergence than under-damping ($c_1=1$). Increasing $c_2$ (gradient correction) also tends to improve the performance, but to a small extent in the ill-conditioned case once $c_2\geq 1$. The error plots when $c_2=1/2$ in the ill-conditioned case appear to be a sum of two oscillations.

\begin{figure}[htbp]
    \centering
    \begin{subfigure}{0.49\textwidth}
       \includegraphics[width=\textwidth]{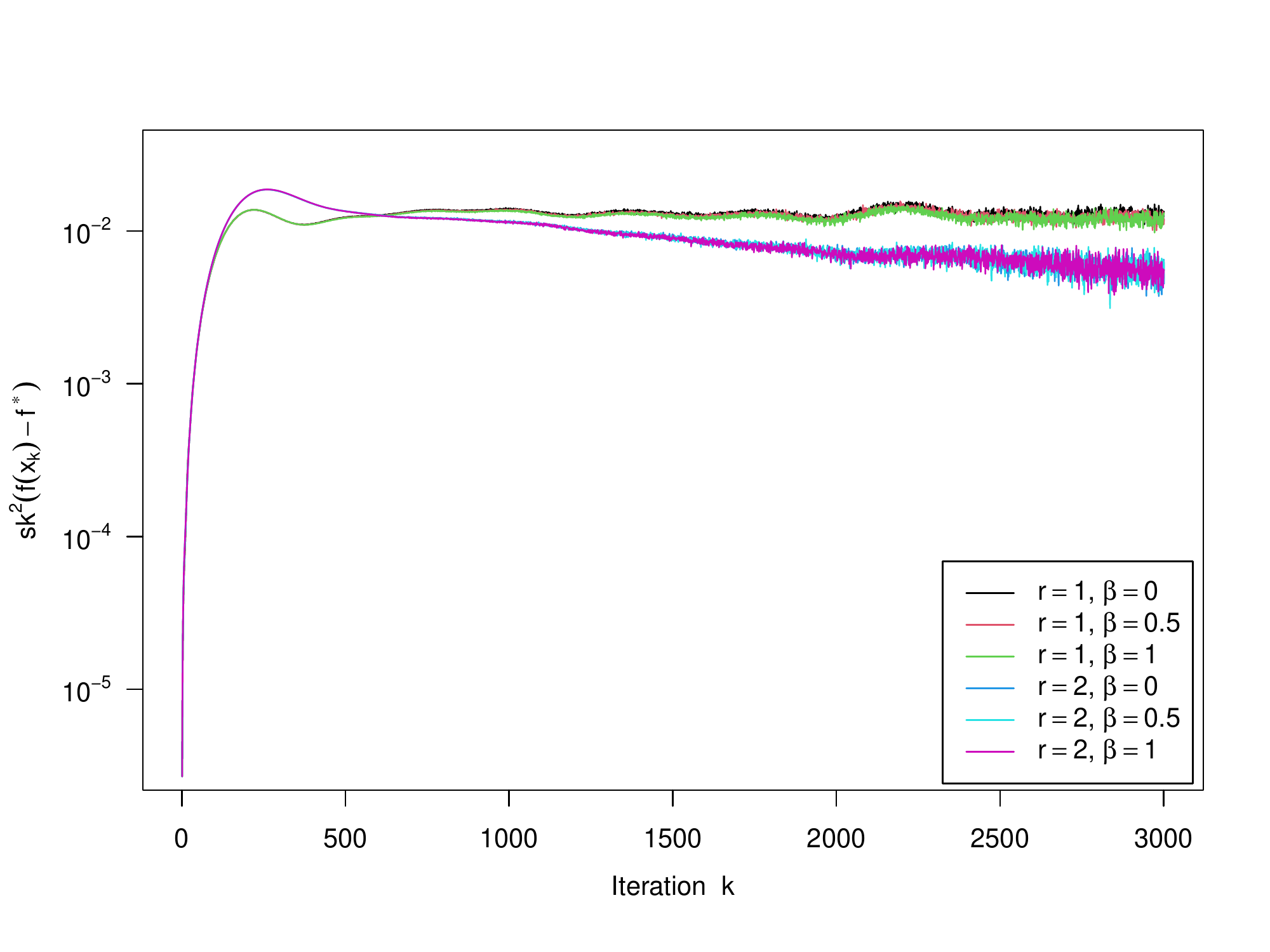}
       \vspace{-0.45in}
       \caption{$s=0.05/\|A\|$}
    \end{subfigure}
    \hfill
    \begin{subfigure}{0.49\textwidth}
       \includegraphics[width=\textwidth]{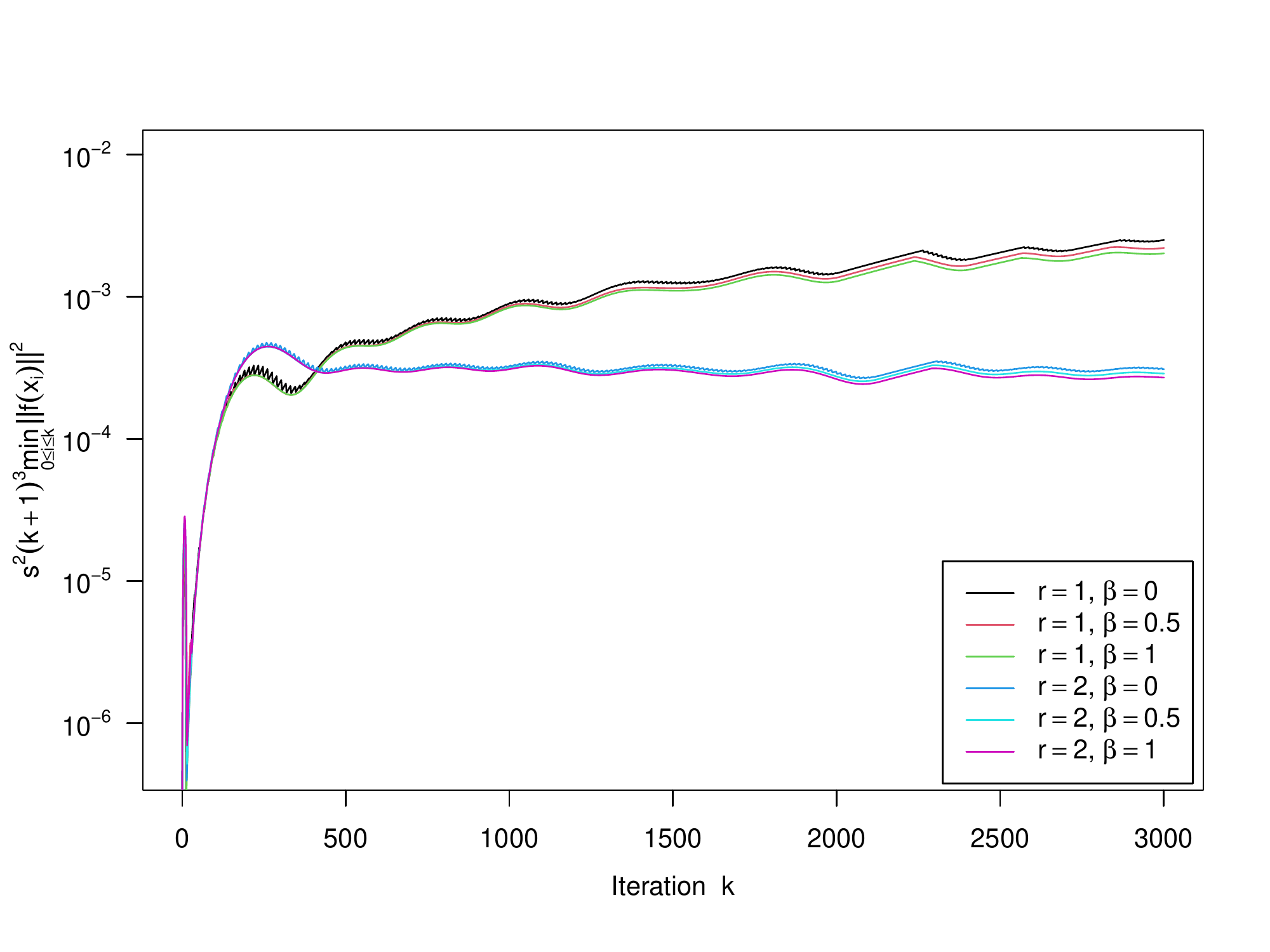}
       \vspace{-0.45in}
       \caption{$s=0.05/\|A\|$}
    \end{subfigure}
    \vspace{-0.2in}

    \begin{subfigure}{0.49\textwidth}
       \includegraphics[width=\textwidth]{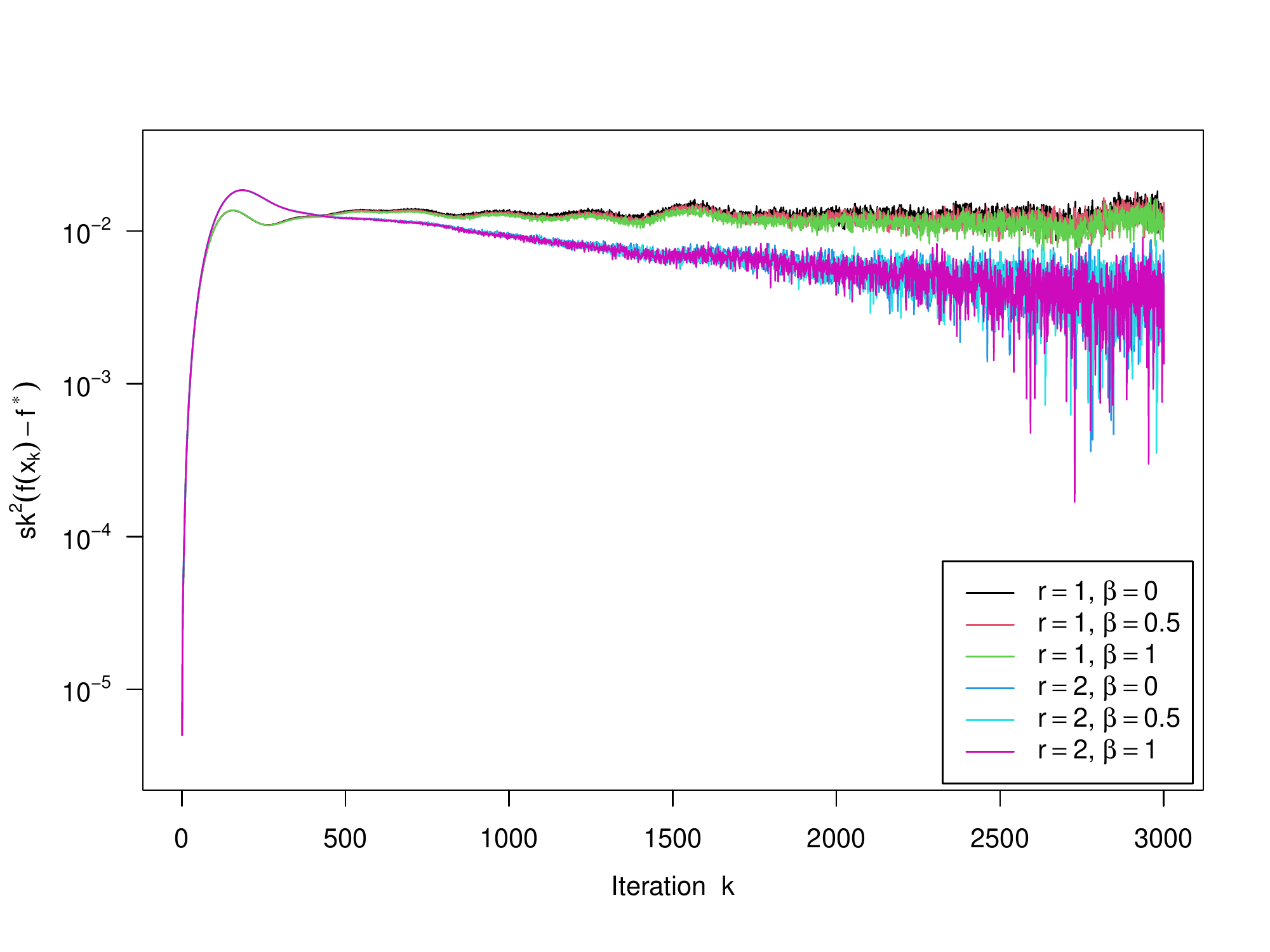}
       \vspace{-0.45in}
       \caption{$s=0.1/\|A\|$}
    \end{subfigure}
    \hfill
    \begin{subfigure}{0.49\textwidth}
       \includegraphics[width=\textwidth]{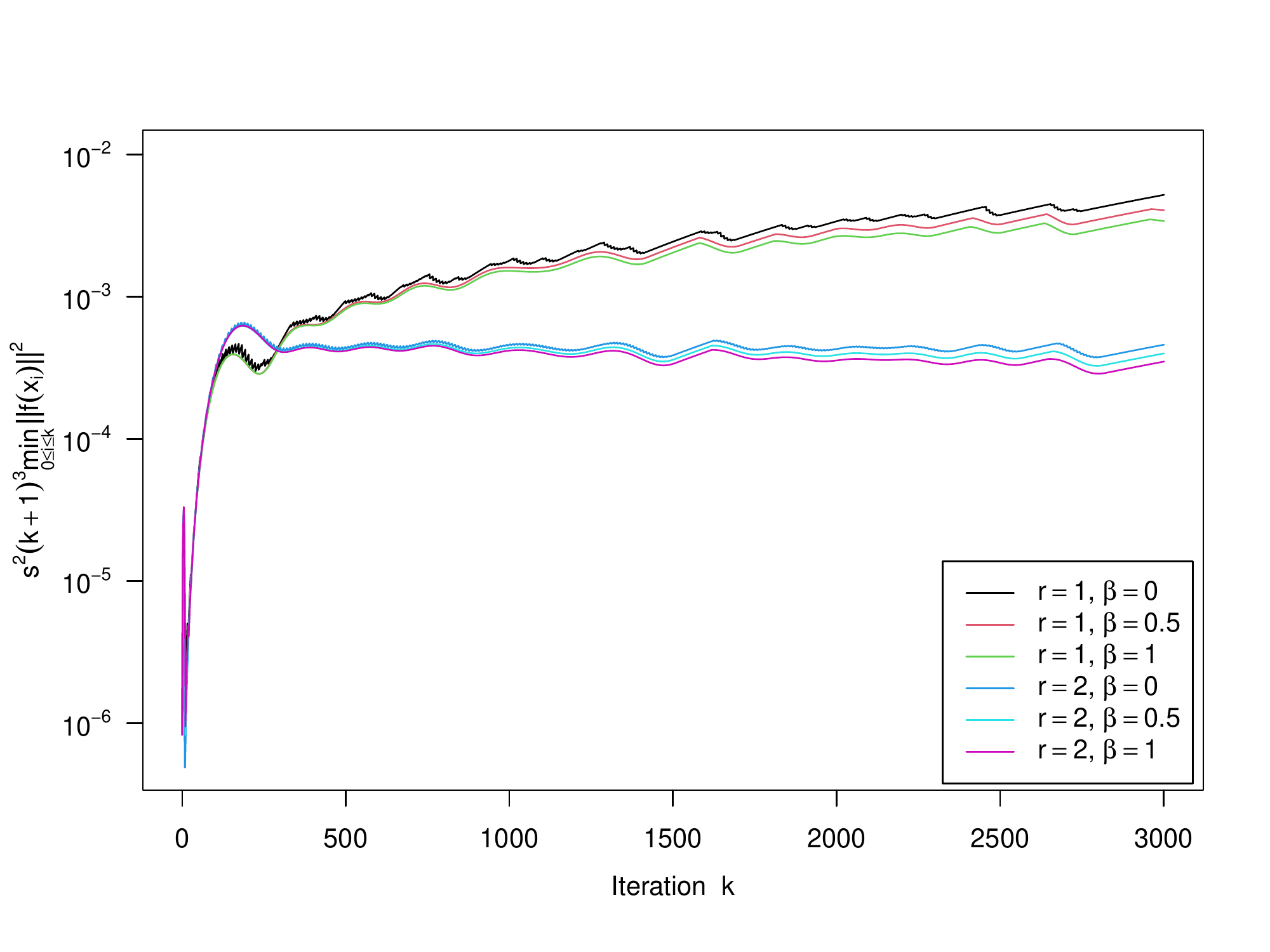}
       \vspace{-0.45in}
       \caption{$s=0.1/\|A\|$}
    \end{subfigure}
    \vspace{-0.2in}

    \begin{subfigure}{0.49\textwidth}
       \includegraphics[width=\textwidth]{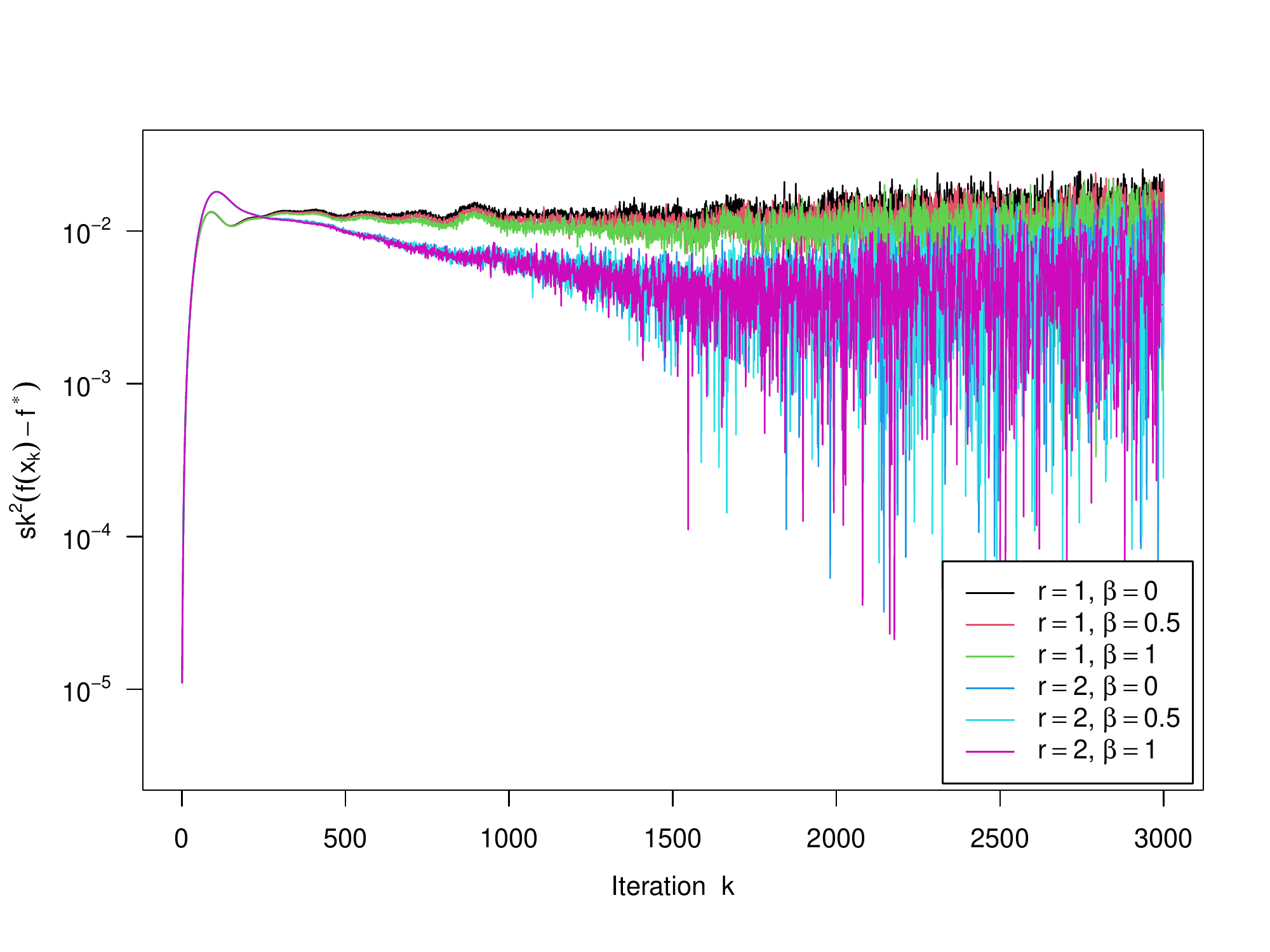}
       \vspace{-0.45in}
       \caption{$s=0.3/\|A\|$}
    \end{subfigure}
    \hfill
    \begin{subfigure}{0.49\textwidth}
       \includegraphics[width=\textwidth]{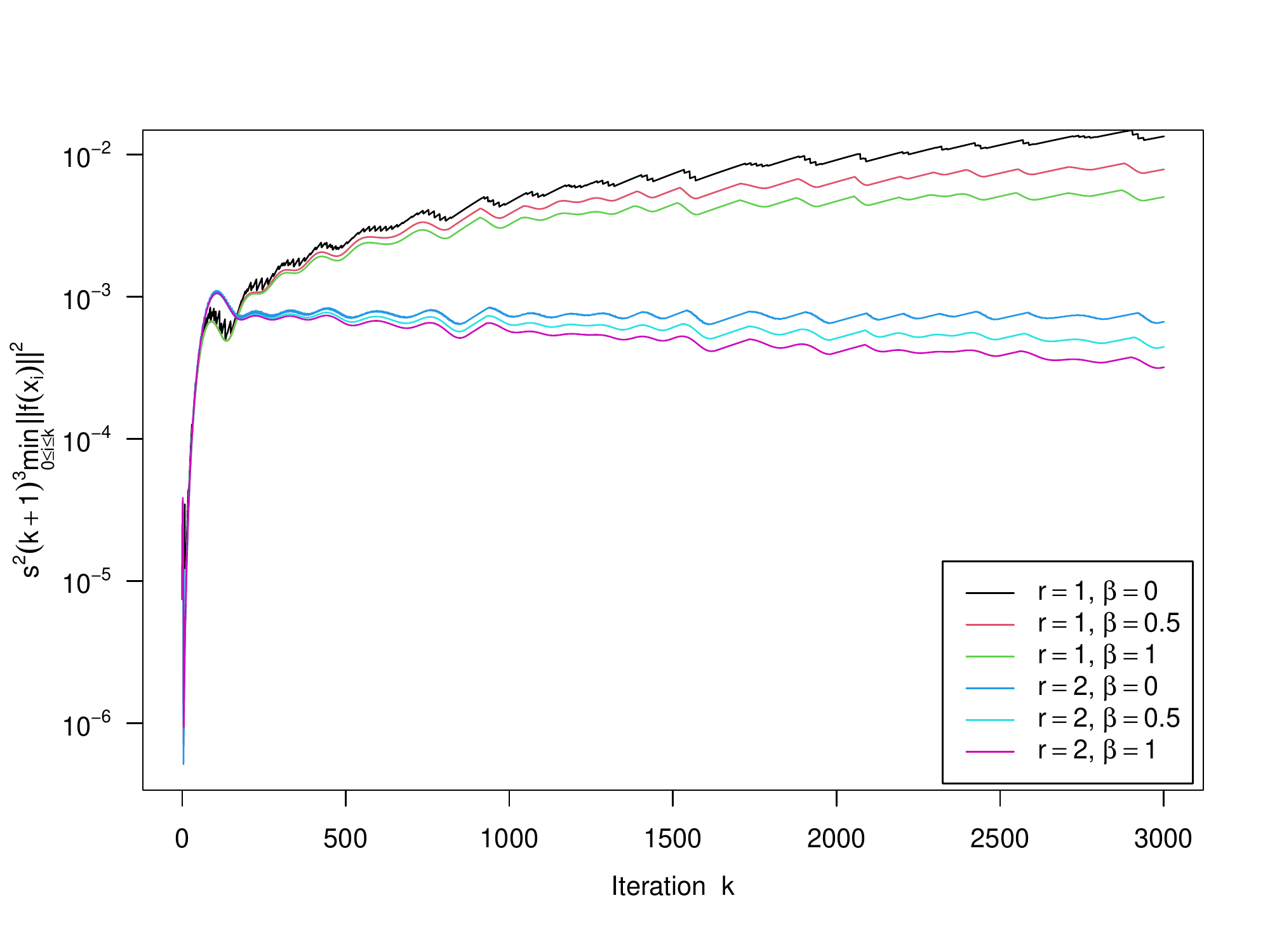}
       \vspace{-0.45in}
       \caption{$s=0.3/\|A\|$}
    \end{subfigure}

    \caption{Scaled errors and squared gradient norms in minimizing $f(x)=\frac{1}{2}x^\T Ax + b^\T x$ by (\ref{eq:extended-NAG-C-shi}) under different step sizes, where $A=B^\T B$ for $B\in\bbR^{500\times 500}$, $b\in\bbR^{500}$. All entries in $B$ and $b$ are {i.i.d. draws} from $U(0,1)$, and $\|A\|$ is the spectral norm of $A$.  We take $\gamma = 1$, $\sigma_{k+1} = \frac{k}{k+r+1}$ for $r=1,2$ and $\beta=0, 0.5, 1$ in (\ref{eq:extended-NAG-C-shi}).}
    \label{fig: convex quadratic}
\end{figure}

\begin{figure}
    \centering
    \begin{subfigure}{0.49\textwidth}
       \includegraphics[width=\textwidth]{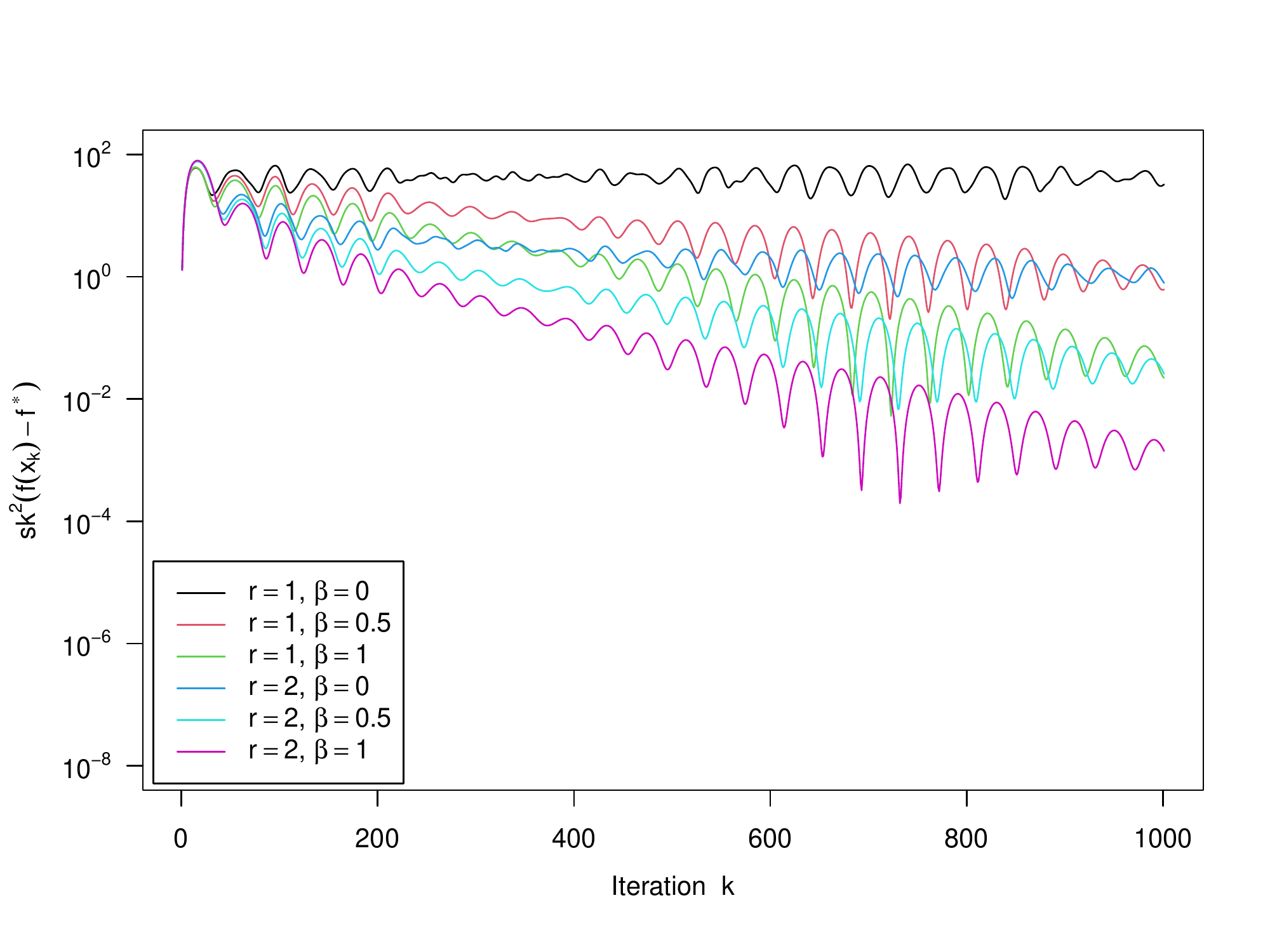}
       \vspace{-0.45in}
       \caption{$s=0.5$}
    \end{subfigure}
    \hfill
    \begin{subfigure}{0.49\textwidth}
       \includegraphics[width=\textwidth]{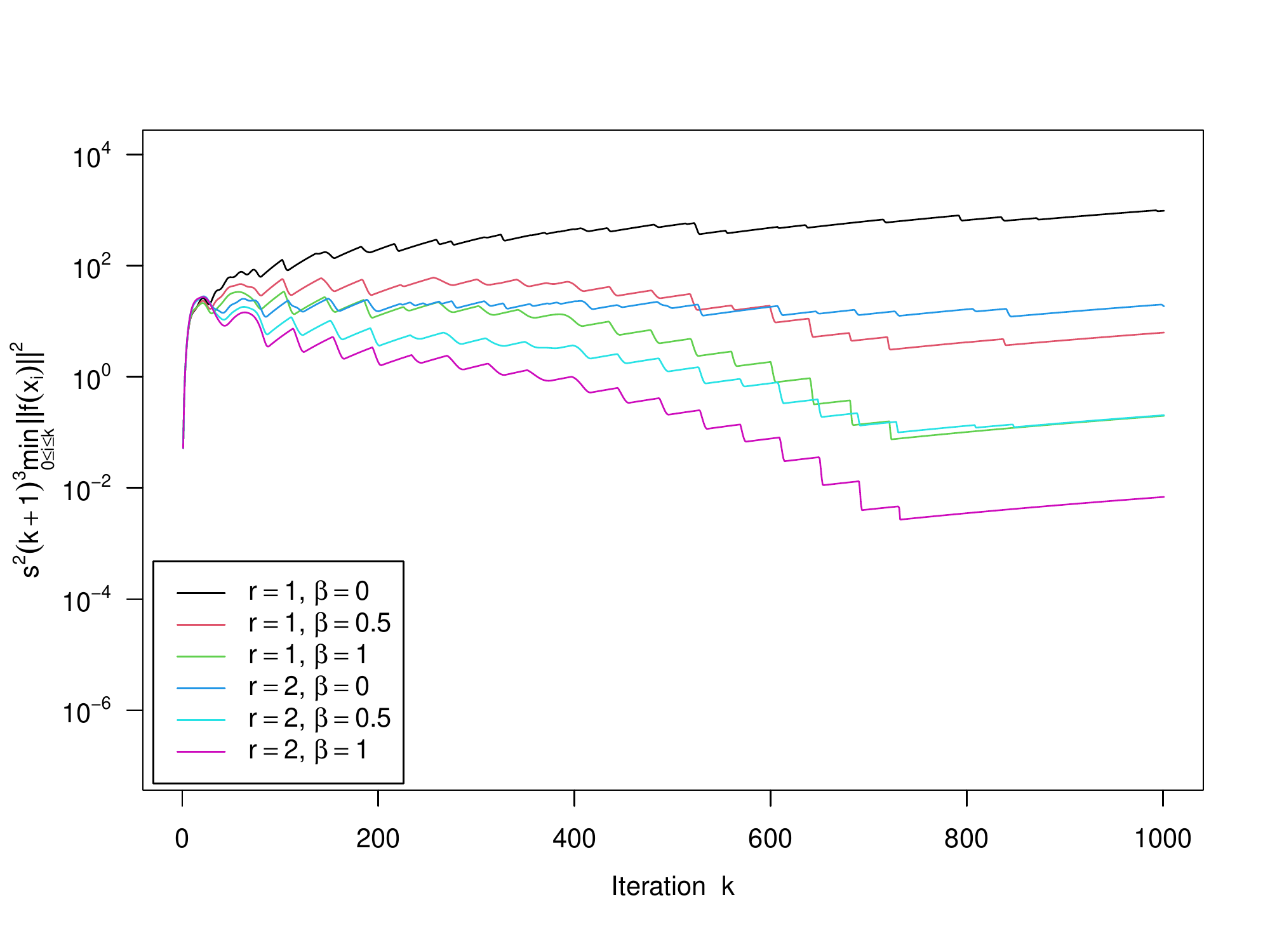}
       \vspace{-0.45in}
       \caption{$s=0.5$}
    \end{subfigure}
    \vspace{-0.2in}

    \begin{subfigure}{0.49\textwidth}
       \includegraphics[width=\textwidth]{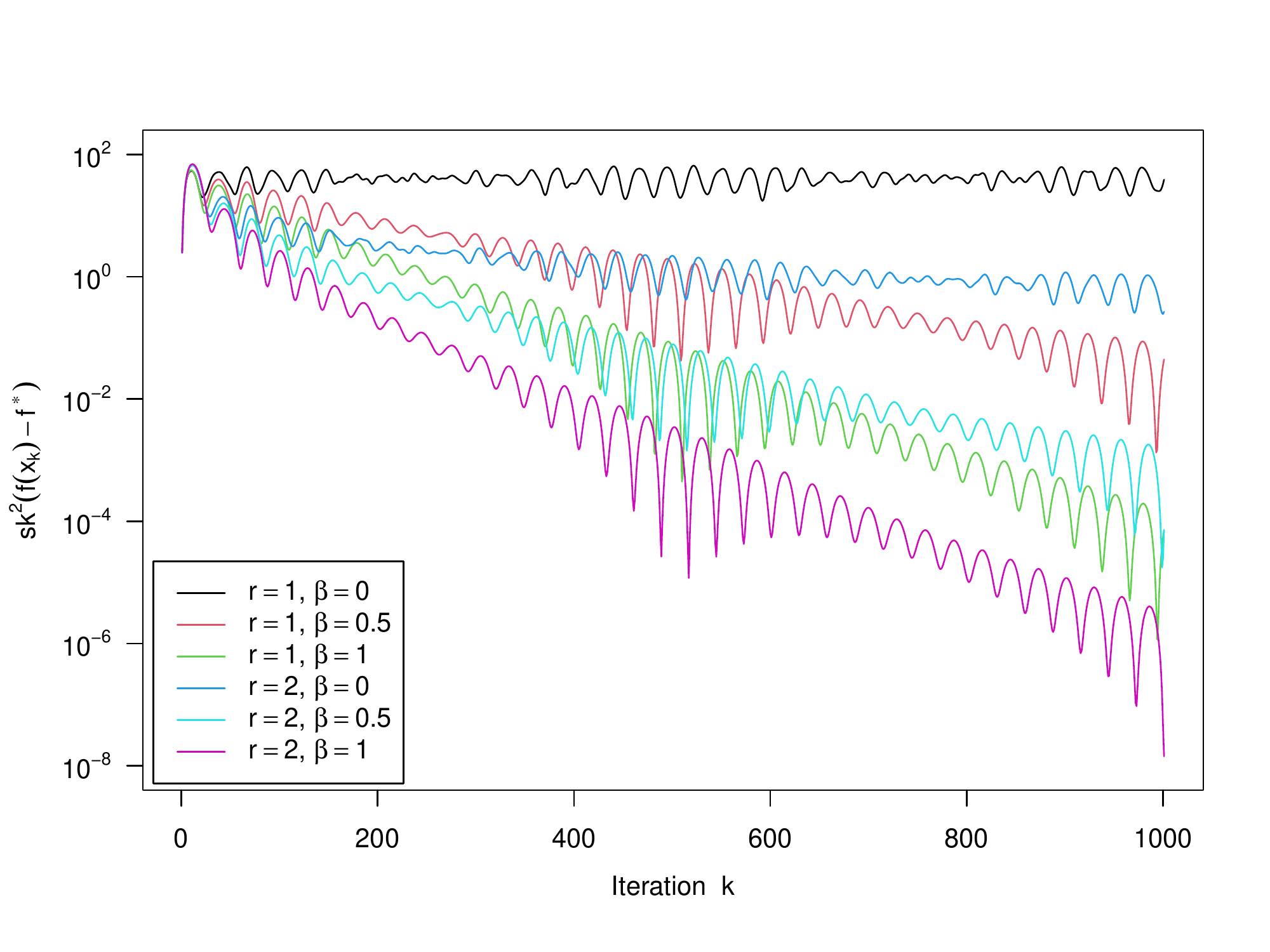}
       \vspace{-0.45in}
       \caption{$s=1$}
    \end{subfigure}
    \hfill
    \begin{subfigure}{0.49\textwidth}
       \includegraphics[width=\textwidth]{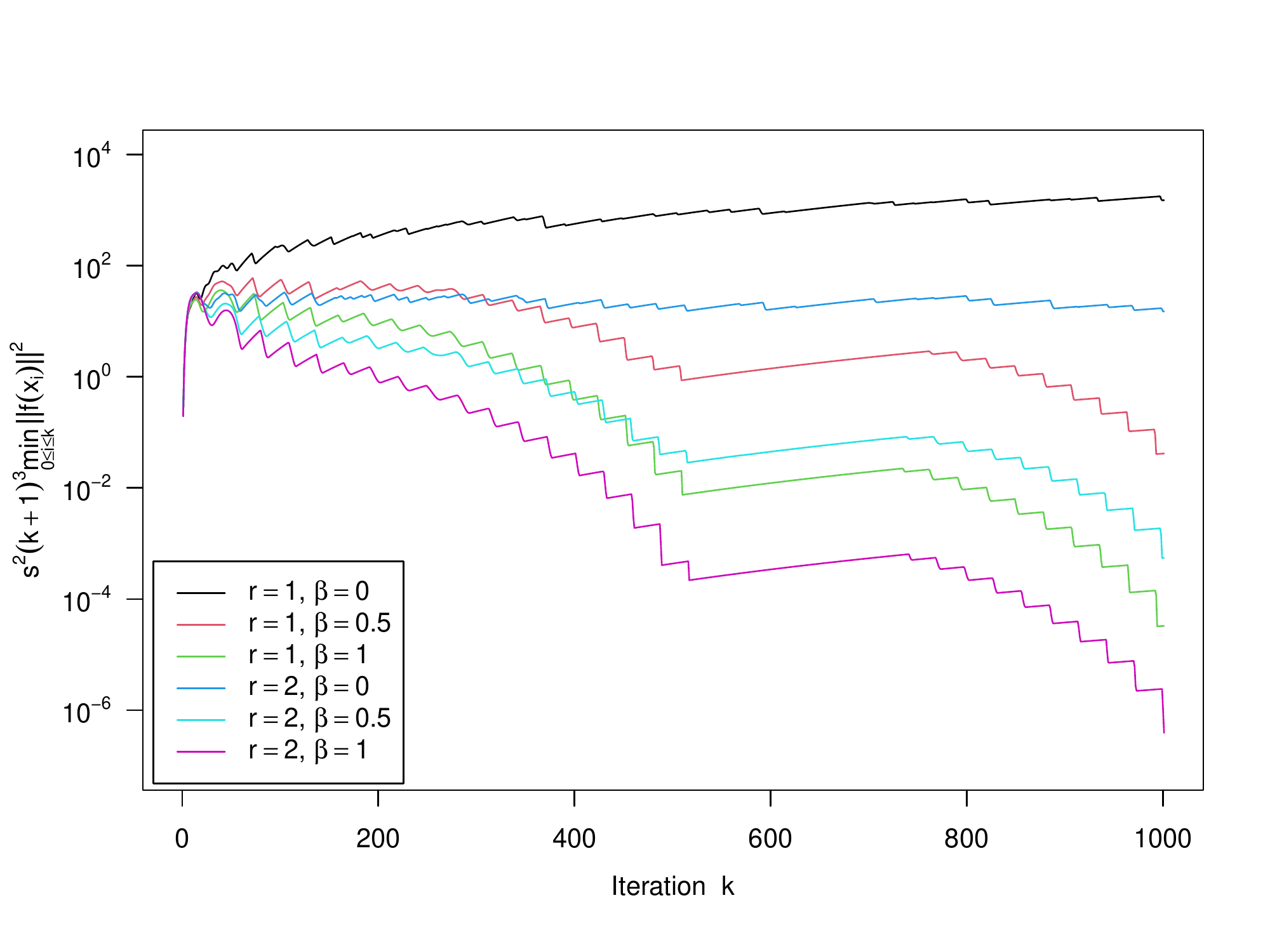}
       \vspace{-0.45in}
       \caption{$s=1$}
    \end{subfigure}
    \vspace{-0.2in}

    \begin{subfigure}{0.49\textwidth}
       \includegraphics[width=\textwidth]{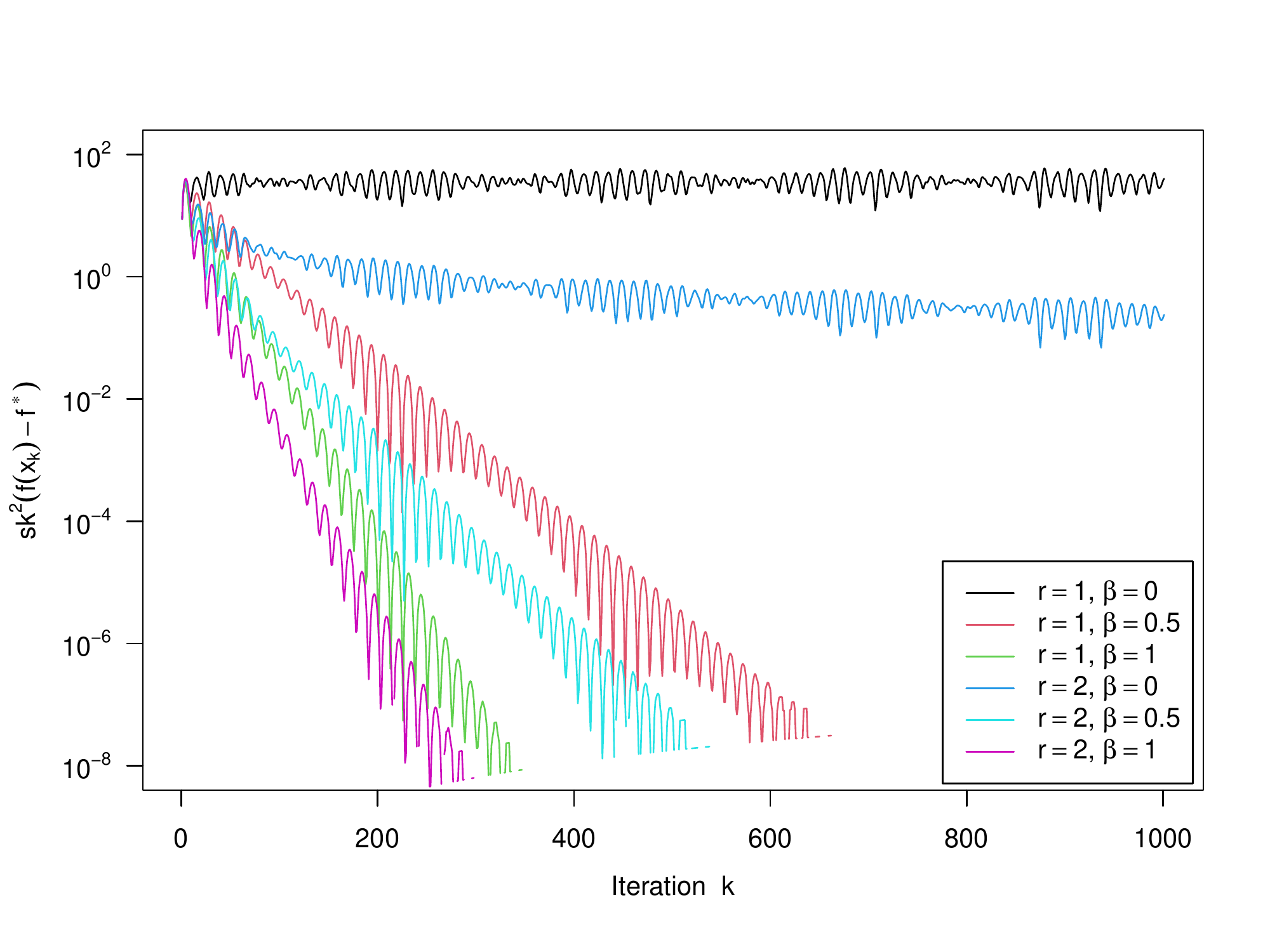}
       \vspace{-0.45in}
       \caption{$s=5$}
    \end{subfigure}
    \hfill
    \begin{subfigure}{0.49\textwidth}
       \includegraphics[width=\textwidth]{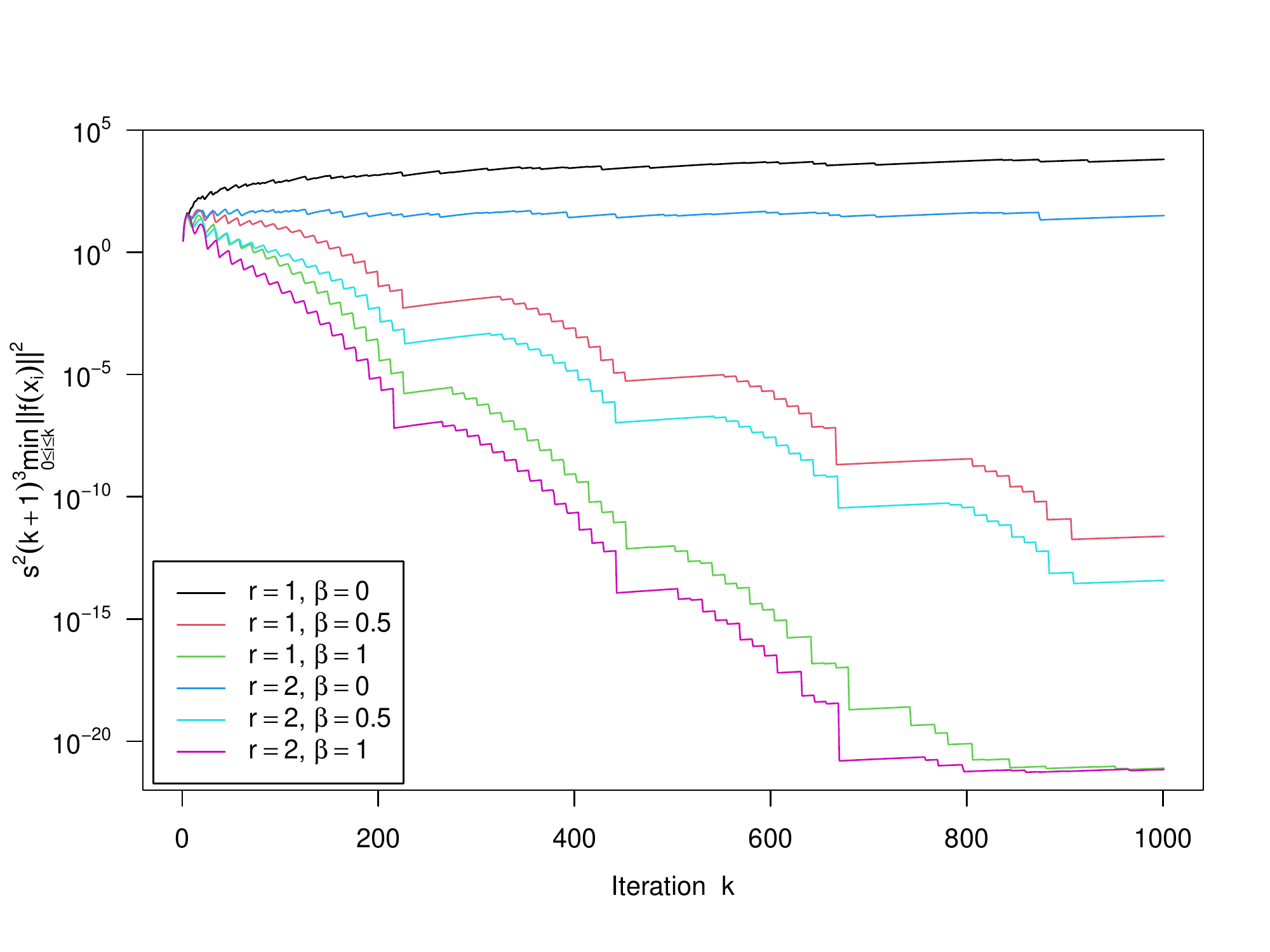}
       \vspace{-0.45in}
       \caption{$s=5$}
    \end{subfigure}

    \caption{Scaled errors and squared gradient norms in minimizing $f(x)=\rho \log \sum_{i=1}^{200} \me^{\frac{a_i^\T x - b_i}{\rho}}$ by (\ref{eq:extended-NAG-C-shi}) under different step sizes, where $A=[a_1,\ldots, a_{200}]\in\bbR^{50\times 200}$, $b\in\bbR^{200}$, and $\rho=20$. All entries in $A$ and $b$ are i.i.d. draws from $\N(0,1)$.  We take $\gamma = 1$, $\sigma_{k+1} = \frac{k}{k+r+1}$ for $r=1,2$ and $\beta=0, 0.5, 1$ in (\ref{eq:extended-NAG-C-shi}).}
    \label{fig: convex log-sum-exp}
\end{figure}

\textbf{General convex setting.}
Following \cite{su2016differential}, we consider a quadratic objective $f(x)=\frac{1}{2}x^\T Ax + b^\T x$ and the log-sum-exp objective $f(x)=\rho \log \sum_{i=1}^{200} \me^{\frac{a_i^\T x - b_i}{\rho}}$. The quadratic function is strongly convex but with $\mu\approx 0$, and the log-sum-exp function is not strongly convex. Figure \ref{fig: convex quadratic} and \ref{fig: convex log-sum-exp} present the numerical results of applying (\ref{eq:extended-NAG-C-shi}) to these two objectives respectively. We set $\gamma\equiv 1$, $\sigma_{k+1}=\frac{k}{k+r+1}$ for $r=1,2$ (boundary) and $\beta = 0, 0.5$ (boundary), $1$. In particular, NAG-C corresponds to $r=2$ and $\beta=1$. In the log-sum-exp case, the minimizer has no closed form. We obtain an approximation of $f^*$ by running NAG-C, which converges fast in hundreds of iterations using a relatively large step size. See Figure \ref{fig: convex log-sum-exp} (e).

In both Figure \ref{fig: convex quadratic} and Figure \ref{fig: convex log-sum-exp}, for fixed $\beta$, increasing $r=1$ to $r=2$ significantly accelerates the convergence of $f(x_k)-f^*$ and $\min_{0\leq i \leq k}\|\nabla f(x_i)\|^2$. In particular, it can be observed in both figures that the bound $\min_{0\leq i \leq k}\|\nabla f(x_i)\|^2 = O(1/k^3)$ may fail when $r=1$. When $r$ is fixed, increasing $\beta$ markedly improves the performance in decreasing $f(x_k)-f^*$ and $\min_{0\leq i \leq k}\|\nabla f(x_i)\|^2$ in Figure \ref{fig: convex log-sum-exp}, especially when the step size is large. A similar beneficial effect of a larger $\beta$ can be observed from Figure \ref{fig: convex quadratic} in decreasing $\min_{0\leq i \leq k}\|\nabla f(x_i)\|^2$, but the effect there is less obvious in decreasing $f(x_k)-f^*$, where all algorithms exhibit oscillations during iterations, and the oscillations become stronger as the step size increases.

\end{document}